# Coupled $G_2$-instantons

Agnaldo A. da Silva Jr.
Universidade Estadual de Campinas (UNICAMP)
agnjr1998@gmail.com

Mario Garcia-Fernandez
Instituto de Ciencias Matemáticas (CSIC-UAM-UC3M-UCM)
mario.garcia@icmat.es

Jason D. Lotay
University of Oxford
jason.lotay@maths.ox.ac.uk

Henrique N. Sá Earp
Universidade Estadual de Campinas (UNICAMP)
henrique.saearp@ime.unicamp.br

**Abstract**

We introduce the coupled instanton equations for a metric, a spinor, a three-form, and a connection on a bundle, over a spin manifold. Special solutions in dimensions 6 and 7 arise, respectively, from the Hull–Strominger and the heterotic $G_2$ system. The equations are motivated by recent developments in theoretical physics and can be recast using generalized geometry; we investigate how coupled instantons relate to generalized Ricci-flat metrics and also to Killing spinors on a Courant algebroid. We present two open questions regarding how these different geometric conditions are intertwined, for which a positive answer is expected from recent developments in the physics literature by De la Ossa, Larfors and Svanes, and in the mathematics literature on Calabi–Yau manifolds, in recent work by the second-named author with González Molina.

We give a complete solution to the first of these problems, providing a new method for the construction of instantons in arbitrary dimensions. For $G_2$-structures with torsion coupled to $G_2$-instantons, in dimension 7, we also establish results around the second problem. The last part of the present work carefully studies the approximate solutions to the heterotic $G_2$-system constructed by the third and fourth authors on contact Calabi–Yau 7-manifolds, for which we prove the existence of approximate coupled $G_2$-instantons and generalized Ricci-flat metrics.

## 1 Introduction

The study of instantons in higher dimensional manifolds is driven by the hope of defining invariants [DT98], emulating the Donaldson invariants obtained from the moduli space of anti-self-dual connections in four dimensions (cf. [DK90]). This programme has motivated substantial activity in recent decades, particularly in the realm of manifolds with special holonomy, such as Calabi–Yau manifolds, see e.g. [SE15, Ste23], and 7-manifolds endowed with a torsion-free $G_2$-structure, see e.g. [Wal13, SEW15, Wal16, LO18, MNE21, ST24]. A new phenomenon that has recently emerged in the physics literature is the construction of *instantons from instantons* in the study of compactifications and domain wall solutions of heterotic string theory on a 7-dimensional manifold $M$ [dlOLS18a, dlOLS18b]. The basic idea is to start with an integrable $G_2$-structure $\varphi$ and a $G_2$-instanton connection $\theta$ on a bundle $P \to M$, that is,

$$F_\theta \wedge \psi = 0,$$

for $\psi = *\varphi$, solving the so-called *heterotic Bianchi identity*

$$dH = \langle F_\theta \wedge F_\theta \rangle. \tag{1.1}$$

The four-form $\langle F_\theta \wedge F_\theta \rangle \in \Omega^4(M)$ depends on a choice of invariant bilinear form on the Lie algebra of the structure group of $P$, often dependent on a real (non-zero) constant denoted $\alpha'$, and represents a multiple of the *first Pontryagin class* of $P$ via Chern–Weil theory. The three-form $H \in \Omega^3(M)$ corresponds to the torsion of the characteristic connection constructed by Friedrich and Ivanov [FI03], and can be written in terms of the torsion forms of $\varphi$ as

$$H = \frac{1}{6}\tau_0 \varphi - \tau_1 \lrcorner \psi - \tau_3.$$

From this data, one obtains an instanton $D$ on the extended bundle $Q = TM \oplus \mathrm{ad}P$ (Theorem 4.6), with respect to the same $G_2$-structure $\varphi$:

$$F_D \wedge \psi = 0. \tag{1.2}$$

The connection $D$ has an intricate yet explicit dependence on $\varphi$, which motivates branding a solution of (1.2) a *coupled $G_2$-instanton*.

Analogues of this interesting new phenomenon have recently been exploited in complex geometry, leading to obstructions to the existence of canonical metrics on a compact non-Kähler Calabi–Yau manifold endowed with solutions of the Hull–Strominger system and pluriclosed Hermitian metrics with vanishing Bismut–Ricci form (see [GFGM23, GFJS23]). Coupled instantons on Calabi–Yau manifolds also play an important role in recent developments in non-Kähler mirror symmetry [ACDAdLHGF24], and they are expected to provide *instantons in higher gauge theory*, where the structure group of the bundle is replaced by a mild category [GFRT20a, TD24].

Motivated by these instances, we introduce and study a general notion of *coupled G-instanton*, for a class of $G$-structures arising from parallel spinors on an oriented spin manifold $M$ endowed with a principal bundle $P$. The unknowns for the equations are tuples $(g, H, \theta, \eta)$, where $g$ is a Riemannian metric on $M$, $(H, \theta)$ is a solution of the heterotic Bianchi identity (1.1), and $\eta$ is a non-zero real spinor. We say that such a tuple satisfies the *coupled instanton equation* if

$$F_D \cdot \eta = 0, \tag{1.3}$$

where $D$ is the connection on the extended tangent bundle $TM \oplus \mathrm{ad}P$ defined by

$$D = \begin{pmatrix} \nabla^- & \mathbb{F}^\dagger \\ -\mathbb{F} & d^\theta \end{pmatrix}. \tag{1.4}$$

Here, $\nabla^-$ is the spin connection on $M$ with skew-symmetric torsion $-H$, given by

$$\nabla^- = \nabla^g - \frac{1}{2} g^{-1} H,$$

and $\mathbb{F} \in \Omega^1(\mathrm{Hom}(T, \mathrm{ad}P))$ is a suitable tensor determined by the curvature of $\theta$ with formal adjoint $\mathbb{F}^\dagger$ (see (2.36) for a more precise definition).

Coupled $G$-instantons are motivated by the unifying language of generalized geometry [Hit03]. As a matter of fact, the basic building blocks for the construction of coupled instantons in six and seven real dimensions correspond to solutions of the *gravitino equation* on a Courant algebroid of string type [GF19]. More explicitly, those tuples $(g, H, \theta, \eta)$ satisfying the heterotic Bianchi identity (1.1) and the *gravitino constraints*

$$\nabla^+ \eta = 0, \qquad F_\theta \cdot \eta = 0, \tag{1.5}$$

where $\nabla^+ = \nabla^g + \frac{1}{2} g^{-1} H$ is the spin connection on $M$ with skew-symmetric torsion $H$. When supplemented with the so-called *dilatino equation*, formulated in terms of a Dirac-type operator acting on $\eta$ (see Lemma 2.12), these equations correspond to the *Killing spinor equations on a Courant algebroid* [GFRT16, GF19] (see Definition 2.11), and hence they are closely related to the Hull–Strominger system [Hul86, Str86] on Calabi–Yau manifolds and also to the heterotic $G_2$ system [FIUV11, dlOLS18a, CGFT22].

In the present work we investigate how coupled instantons and solutions of the gravitino equation are intertwined, and ask whether any solution of (1.5) and the heterotic Bianchi identity (1.1) induces a coupled instanton: see Problem 1. We also study how coupled instantons are related to generalized Ricci-flat metrics, given by solutions of the system of equations

$$\begin{aligned} \mathrm{Rc} - \frac{1}{4} H^2 + \langle i_{v_i} F_\theta, i_{v_i} F_\theta \rangle + \frac{1}{2} L_{\zeta^\sharp} g &= 0, \\ d^* H - d\zeta + i_{\zeta^\sharp} H &= 0, \\ d_\theta^* F_\theta + *(F_\theta \wedge *H) + i_{\zeta^\sharp} F_\theta &= 0, \end{aligned} \tag{1.6}$$

for a tuple $(g, H, \theta)$, as before, and a one-form $\zeta \in \Omega^1$. In this respect, we propose as an open question in Problem 2 to characterise the precise conditions a coupled instanton needs to satisfy, in terms of the $G$-structure determined by the spinor, in order to solve (1.6) for a suitable choice of one-form $\zeta$.

When $\zeta = d\phi$ for a smooth function $\phi$, the equations (1.6) match the equations of motion of the heterotic supergravity for the metric, the three-form flux, and the gauge field, in the mathematical physics literature, see e.g. [GF14, Mol24]. The parallelism with the physics setup (see Remark 2.16) leads naturally to ask whether the equations (1.6) imply that the analogue of the equation of motion for the *dilaton field* is satisfied. This is a scalar equation, given by the vanishing of the function

$$R_g - \frac{1}{2} |H|^2 + |F_\theta|^2 - 2 d^* \zeta - |\zeta|^2 \in C^\infty(M). \tag{1.7}$$

This scalar quantity plays an important role in the theory of generalized Ricci flow, being closely related to the volume density of the generalized Perelman energy functional, see [AMP24, GFS20, ŠV20, SSCV24]. We explore this interesting aspect in Section 2.4.

We prove three main results in the present work. Theorem 2.32 provides a precise answer to Problem 1 for oriented spin manifolds of arbitrary dimension, from which we derive in particular Theorem 4.6 for $G_2$-structures with torsion coupled to $G_2$-instantons. Specialising henceforth to dimension 7, Theorem 4.9 investigates the relation between coupled instantons and generalized Ricci-flat metrics, in relation to Problem 2. More precisely, the coupled $G_2$-instantons constructed in Theorem 4.6 are, in fact, generalized Ricci-flat, and we reach the strong conclusion that the gravitino equations (1.5) together with the heterotic Bianchi identity (1.1) actually imply solving the full heterotic $G_2$ system (3.12), as known to be the case for coupled $\mathrm{Spin}(7)$-instantons, see Remark 4.10. These two results rely on some aspects of the theory of parallel spinors and $G$-structures for connections with skew-symmetric torsion, initiated in the seminal paper [FI03] and further developed in [FKMS97, AF03, AF10, ACFH15]. Third, Theorem 5.15 is devoted to the study of examples arising from the analysis of the heterotic $G_2$ system on certain 7-manifolds which are circle bundles over Calabi–Yau 3-orbifolds and carry natural integrable $G_2$-structures [CARDSE20, LSE23]. On these contact Calabi–Yau 7-manifolds, we find approximate coupled $G_2$-instantons and generalized Ricci curvature which is approximately zero in a precise quantitative sense, relative to the (string) constant $\alpha' > 0$ appearing in the heterotic $G_2$ system. In these examples, as $\alpha' \to 0$, the circle fibres necessarily shrink to zero size.

This paper is organized as follows. In Section 2 we recall the background in generalized geometry and introduce our equations of interest: the gravitino equation, the dilatino equation with parameter $\lambda \in \mathbb{R}$ (see Theorem 2.11) and the coupled instanton equations (see Theorem 2.26). Collectively, we call a solution of the first two sets of equations a *Killing spinor* with parameter $\lambda$. In Section 2.5 we propose two open questions which relate solutions of the gravitino equation (1.5) with coupled instantons (see Problem 1), and also to generalized Ricci-flat metrics (see Problem 2). In Section 3 we specify to the case where the base manifold is seven-dimensional. After briefly reviewing the necessary background on $G_2$-structures, we characterise the Killing spinors with parameter $\lambda$ in terms of $G_2$-structures with torsion coupled to $G_2$-instantons, via a mild generalisation of the *heterotic $G_2$ system* (where the torsion class $\tau_1$ may not be necessarily exact) (see Proposition 3.6). The rest of this section is devoted to prove, by means of generalized geometry techniques [GF19, CSCW11], that any solution of this system gives a generalized Ricci-flat metric with constant generalized scalar curvature, for a suitable choice of divergence operator. Section 4 is devoted to give a precise answer to Problem 1 and to make progress in Problem 2 in the case of a seven-dimensional oriented spin manifold. Finally, Section 5 introduces and studies a notion of approximate instanton, based on [LSE23], and provides conditions for approximate generalized Ricci-flatness, which are shown to hold in the contact Calabi–Yau examples from [LSE23].

**Acknowledgements:** The authors thank Lucía Martín-Merchán and Jeffrey Streets for valuable discussions and Raúl González Molina for helpful comments and a careful reading of the manuscript. MGF was partially supported by the Spanish Ministry of Science and Innovation, through the *'Severo Ochoa Programme for Centres of Excellence in R&D'* (CEX2019-000904-S) and under grants PID2022-141387NB-C22 and CNS2022-135784. JDL was partially supported by the *Simons Collaboration on Special Holonomy in Geometry, Analysis, and Physics* (#24071 Jason Lotay). HSE was supported by the São Paulo Research Foundation (Fapesp) [2021/04065-6] *BRIDGES collaboration* and a Brazilian National Council for Scientific and Technological Development (CNPq) [311128/2020-3] Productivity Grant level 1D. ASJr was supported by the São Paulo Research Foundation (Fapesp) [2021/11603-4]. JDL and HSE benefited from a *Newton Mobility* bilateral collaboration (2019-2023), granted by the UK Royal Society [NMG\R1\191068].

## 2 Coupled instantons and Killing spinors in generalized geometry

### 2.1 Background on string algebroids

We recall the necessary background material on Courant algebroids of string type following [GFRT20b] (see also [BH15, GF14]).

**Definition 2.1.** *A* Courant algebroid $(E, \langle \cdot, \cdot \rangle, [\cdot, \cdot], \pi)$ *over a manifold $M$ consists of a vector bundle $E \to M$ endowed with a non-degenerate symmetric bilinear form $\langle \cdot, \cdot \rangle$ and a (Dorfman) bracket $[\cdot, \cdot]$ on $\Omega^0(E)$, and a bundle map $\pi \colon E \to T$, called an* anchor map*, such that the following axioms are satisfied, for all $a, b, c \in \Omega^0(E)$ and $f \in C^\infty(M)$:*

*(1)* $[a, [b, c]] = [[a, b], c] + [b, [a, c]]$,

*(2)* $\pi [a, b] = [\pi(a), \pi(b)]$,

*(3)* $[a, fb] = f[a, b] + \pi(a)(f) b$,

*(4)* $\pi(a) \langle b, c \rangle = \langle [a, b], c \rangle + \langle b, [a, c] \rangle$,

*(5)* $[a, b] + [b, a] = \mathcal{D} \langle a, b \rangle$.

Here, $\mathcal{D} \colon C^\infty(M) \to \Omega^0(E)$ denotes the composition of the exterior differential $d \colon C^\infty(M) \to \Omega^1(M)$, the dual map $\pi^* \colon T^* \to E^*$ and the isomorphism $E^* \cong E$ provided by $\langle \cdot, \cdot \rangle$.

We will denote a Courant algebroid $(E, \langle \cdot, \cdot \rangle, [\cdot, \cdot], \pi)$ simply by $E$. Using the isomorphism $\langle \cdot, \cdot \rangle : E \to E^*$, we obtain a complex of vector bundles

$$T^* \xrightarrow{\pi^*} E \xrightarrow{\pi} T. \tag{2.1}$$

We will say that $E$ is *transitive* if the anchor map $\pi$ in (2.1) is surjective. Given a transitive Courant algebroid $E$ over $M$, there is an associated Lie algebroid

$$A_E := E/(\operatorname{Ker} \pi)^\perp.$$

Furthermore, the subbundle

$$\operatorname{ad}_E := \operatorname{Ker} \pi/(\operatorname{Ker} \pi)^\perp \subset A_E$$

inherits the structure of a bundle of quadratic Lie algebras. Therefore, the bundle $E$ fits into a double extension of vector bundles

$$\begin{aligned} 0 &\longrightarrow T^* \xrightarrow{\pi^*} E \longrightarrow A_E \longrightarrow 0, \\ 0 &\longrightarrow \operatorname{ad}_E \xrightarrow{\pi^*} A_E \xrightarrow{\pi} T \longrightarrow 0. \end{aligned} \tag{2.2}$$

A classification of transitive Courant algebroids has been obtained in [GFRT20b, Proposition A.6] for the special case in which $A_E$ is isomorphic to the Atiyah algebroid of a smooth principal bundle; such a Courant algebroid is said to be of *string type*. In order to state a more precise definition, we briefly discuss the basic example which we will need.

**Example 2.2.** Let $K$ be a real Lie group, whose Lie algebra $\mathfrak{k}$ is endowed with a non-degenerate bi-invariant symmetric bilinear form

$$\langle \cdot, \cdot \rangle : \mathfrak{k} \otimes \mathfrak{k} \longrightarrow \mathbb{R}.$$

Let $p \colon P \to M$ be a smooth principal $K$-bundle, and consider the Atiyah Lie algebroid $A_P := TP/K$. The smooth bundle of Lie algebras $\operatorname{ad} P := \operatorname{Ker} dp \subset A_P$ fits into the short exact sequence of Lie algebroids

$$0 \longrightarrow \operatorname{ad} P \longrightarrow A_P \longrightarrow T \longrightarrow 0.$$

We construct next a transitive Courant algebroid such that the second sequence in (2.2) is canonically isomorphic to the exact sequence of Lie algebroids above. Assume that the *first Pontryagin class* of $P$ associated to $\langle \cdot, \cdot \rangle$ via Chern–Weil theory is trivial:

$$p_1(P) = 0 \in H^4_{dR}(M, \mathbb{R}). \tag{2.3}$$

Then, given a choice of principal connection $\theta$ on $P$, with curvature $F_\theta$, there exists a smooth real three-form $H \in \Omega^3$ such that the *Bianchi identity* holds:

$$dH - \langle F_\theta \wedge F_\theta \rangle = 0. \tag{2.4}$$

Given such a pair $(H, \theta)$, we define a Courant algebroid $E_{P,H,\theta}$ with underlying vector bundle

$$T \oplus \operatorname{ad} P \oplus T^*,$$

non-degenerate symmetric bilinear form

$$\langle X + r + \zeta, X + r + \zeta \rangle = \zeta(X) + \langle r, r \rangle, \tag{2.5}$$

bracket given by

$$\begin{aligned} [X + r + \zeta, Y + t + \eta] &= [X, Y] - F_\theta(X, Y) + d^\theta_X t - d^\theta_Y r - [r, t] \\ &\quad + L_X \eta - \iota_Y d\zeta + \iota_Y \iota_X H \\ &\quad + 2\langle d^\theta r, t \rangle + 2\langle \iota_X F_\theta, t \rangle - 2\langle \iota_Y F_\theta, r \rangle, \end{aligned} \tag{2.6}$$

and anchor map $\pi(X + r + \zeta) = X$. For simplicity, we will abuse notation and use the same symbol for the pairing on the Lie algebra $\mathfrak{k}$ and $E$, which should lead to no confusion. Under the previous hypothesis, $E_{P,H,\theta}$ defines a smooth transitive Courant algebroid over $M$. △

Transitive Courant algebroids as in Example 2.2 fit into the category of *string algebroids* [GFRT20b], which motivates the following definition. The notion of isomorphism we use is the standard one for Courant algebroids, given by (base-preserving) smooth orthogonal bundle morphisms which preserve the bracket and the anchor map (cf. Remark 2.4):

**Definition 2.3.** *A Courant algebroid $E$ over $M$ is of* string type *if it is isomorphic to a Courant algebroid $E_{P,H,\theta}$, as in Example 2.2, for some triple $(P, H, \theta)$ satisfying (2.4). In this case, we will refer to $E$ simply as a* string algebroid.

*Remark* 2.4. Note that the Courant algebroid $E$ in Example 2.2 carries a natural bracket-preserving map $E \to A_P$, induced by the identification $A_P \cong T \oplus \operatorname{ad} P$ provided by the connection $\theta$. This data is often regarded as part of the structure of a string algebroid, and morphisms in the string algebroid category are compatible with this map (cf. [GFRT20b, Definition 2.3]). ◯

## 2.2 Weak Koszul formula and generalized Ricci tensor

Let $M$ be an oriented manifold endowed with a string algebroid $E$. In this section, we recall basic aspects of generalized Riemannian geometry, following [GF19, GFS20].

**Definition 2.5.** *A generalized metric on a string algebroid $E$ is an orthogonal decomposition $E = V_+ \oplus V_-$, so that the restriction of $\langle \cdot, \cdot \rangle$ to $V_+$ is positive definite and that $\pi_{|V_+} : V_+ \to T$ is an isomorphism.*

A generalized metric $V_+ \subset E$ is equivalent to a pair $(g, \sigma)$, where $g$ is a Riemannian metric on $M$ and $\sigma \colon T \to E$ is an isotropic splitting, see e.g. [GF14]. Alternatively, a generalized metric can be encoded in an orthogonal endomorphism $\mathbf{G} \colon E \to E$ such that $\mathbf{G}^2 = \mathrm{Id}$. The orthogonal decomposition $E = V_+ \oplus V_-$ is then recovered from the eigenbundles

$$V_\pm = \mathrm{Ker}\,(\mathbf{G} \mp \mathrm{Id}).$$

We will use the following notation for the induced orthogonal projections

$$\pi_\pm := \frac{1}{2}(\mathbf{G} \pm \mathrm{Id}) \colon E \longrightarrow V_\pm \colon a \longmapsto a_\pm.$$

More explicitly, in our case of interest, the isotropic splitting $\sigma \colon T \to E$ determined by $\mathbf{G}$ induces an isomorphism

$$E \cong T \oplus \mathrm{ad}P \oplus T^*,$$

and hence an explicit string algebroid structure as in Example 2.2, with bracket (2.6), for a uniquely determined $H \in \Omega^3(M)$ and principal connection $\theta$ on $P$ satisfying (2.4), cf. [GF14, Proposition 3.4]. Furthermore, via this identification we have

$$V_+ = \{X + gX : X \in T\}, \quad V_- = \{X + r - gX : X \in T, r \in \mathrm{ad}P\}, \tag{2.7}$$

as well as

$$\mathbf{G} = \begin{pmatrix} 0 & 0 & g^{-1} \\ 0 & -\mathrm{Id} & 0 \\ g & 0 & 0 \end{pmatrix},$$

with orthogonal projections

$$\pi_+(X + r + \zeta) = \frac{1}{2}(X + gX + g^{-1}\zeta + \zeta), \qquad \pi_-(X + r + \zeta) = \frac{1}{2}(X - gX - g^{-1}\zeta + \zeta) + r.$$

It will be useful to introduce the notation

$$\begin{aligned} \sigma_\pm \colon T &\to V_\pm \\ X &\mapsto \sigma_\pm(X) := X \pm gX. \end{aligned} \tag{2.8}$$

*Remark* 2.6. If we take the string algebroid of the form $E = E_{P,H_0,\theta_0}$, as in Example 2.2, then a generalized metric $\mathbf{G}$ on $E$ is equivalent to a triple $(g, b, \theta)$, for $b \in \Omega^2$, and the three-form $H$ above is determined by

$$H = H_0 + 2\langle a \wedge F_{\theta_0} \rangle + \langle a \wedge d^{\theta_0} a \rangle + \tfrac{1}{3}\langle a \wedge [a \wedge a] \rangle + db, \tag{2.9}$$

where $a = \theta_0 - \theta \in \Omega^1(\mathrm{ad}P)$. Condition (2.9) can be expressed more invariantly by the following equivalent condition on the equivariant cohomology of the principal bundle:

$$[p^* H_0 - CS(\theta_0)] = [p^* H - CS(\theta)] \in H^3(P, \mathbb{R})^K,$$

where $p \colon P \to M$ is the canonical projection and $CS(\theta) \in \Omega^3(P)^K$ is the Chern–Simons 3-form of $\theta$. The class $[p^* H_0 - CS(\theta_0)] \in H^3(P, \mathbb{R})^K$ can be regarded as the isomorphism class of a $K$-equivariant (exact) Courant algebroid over the total space of $P$, from which the (transitive) Courant algebroid $E$ is obtained by reduction [BH15, GF14]. ○

In order to introduce natural curvature quantities associated to a generalized metric $\mathbf{G}$, the main difficulty is that there is no uniquely determined analogue of the Levi-Civita connection [CSCW11, GF19]. Instead, there is a weak version of Koszul's formula: a generalized metric $\mathbf{G}$ on $E$ determines a pair of differential operators

$$D_-^+ \colon \Gamma(V_+) \to \Gamma(V_-^* \otimes V_+) \quad \text{and} \quad D_+^- \colon \Gamma(V_-) \to \Gamma(V_+^* \otimes V_-), \tag{2.10}$$

defined on sections $a_- \in \Gamma(V_-)$ and $b_+ \in \Gamma(V_+)$ respectively by

$$\langle a_-, D_-^+ b_+ \rangle = \pi_+[a_-, b_+] \quad \text{and} \quad \langle b_+, D_+^- a_- \rangle = \pi_-[b_+, a_-]. \tag{2.11}$$

In the sequel, we will abuse notation and write simply $D_{a_-}b_+ := \langle a_-, D_-^+ b_+\rangle$, and similarly for $D_+^-$. They satisfy natural Leibniz rules, with respect to the anchor map, for any smooth function $f \in C^\infty(M)$:

$$D_{a_-}(fb_+) = \pi(a_-)(f)b_+ + fD_{a_-}b_+,$$
$$D_{b_+}(fa_-) = \pi(b_+)(f)a_- + fD_{a_+}a_-.$$

To give an explicit formula for the operators in (2.10), we fix a generalized metric $\mathbf{G}$ on $E$ and consider the associated isomorphism $E \cong T \oplus \mathrm{ad}P \oplus T^*$ and pair $(H,\theta)$ satisfying (2.4). Define a pair of metric connections $\nabla^\pm$ on $(T,g)$ with totally skew-symmetric torsion by

$$\nabla_X^+ Y = \nabla_X^g Y + \frac{1}{2} g^{-1} H(X,Y,\cdot), \quad \text{and} \quad \nabla_X^- Y = \nabla_X^g Y - \frac{1}{2} g^{-1} H(X,Y,\cdot), \tag{2.12}$$

where $\nabla^g$ is the Levi-Civita connection of $g$.

**Lemma 2.7** ([GF14, GFRT16]). *Let $\mathbf{G}$ be a generalized metric on a string algebroid $E$ and let $(H,\theta)$ be the uniquely determined solution of (2.4). Then, the differential operators (2.10) and the metric connections (2.12) on $(T,g)$ are related by:*

$$\begin{aligned} D_{b_+}a_- &= \sigma_-(\nabla_Y^- X - g^{-1}\langle i_Y F_\theta, r\rangle) + d_Y^\theta r - F_\theta(Y,X), \\ D_{a_-}b_+ &= \sigma_+(\nabla_X^+ Y - g^{-1}\langle i_Y F_\theta, r\rangle), \end{aligned} \tag{2.13}$$

*where $F_\theta$ is the curvature of $\theta$ and*

$$\begin{aligned} a_- &= \sigma_-(X) + r = X + r - gX, \\ b_+ &= \sigma_+(Y) = Y + gY. \end{aligned} \tag{2.14}$$

Even though the "right" notion of curvature tensor in generalized geometry is still unknown, one can construct a pair of generalized Ricci tensors associated to a generalized metric. For this, it is customary to consider *divergence operators* on the string algebroid $E$ which keep track of the conformal geometry of $E$, cf. [GF19].

**Definition 2.8.** *A divergence operator on a string algebroid $E$ is a map $\mathrm{div} \colon \Omega^0(E) \to C^\infty(M)$ satisfying*

$$\mathrm{div}(fa) = f\,\mathrm{div}(a) + \pi(a)(f), \quad \text{for } f \in C^\infty(M).$$

By definition, the space of divergence operators on $E$ is affine and modelled on the space of sections of $E^* \cong E$. Note that the generalized metric has an associated *Riemannian divergence* defined by

$$\mathrm{div}^\mathbf{G}(X + r + \zeta) = \frac{L_X \mathrm{vol}_M}{\mathrm{vol}_M}, \tag{2.15}$$

where $\mathrm{vol}_M$ is the volume element of $g$. As a result, any divergence $\mathrm{div}$ on $E$ can be expressed in the form

$$\mathrm{div} = \mathrm{div}^\mathbf{G} - \langle \varepsilon, \cdot\rangle$$

for a uniquely determined section $\varepsilon \in \Gamma(E)$.

Naturally associated to a pair $(\mathbf{G}, \mathrm{div})$, given by a generalized metric $\mathbf{G}$ and a divergence operator $\mathrm{div}$, one can construct a pair of generalized Ricci tensors (see Proposition A.2):

$$\mathrm{Rc}^+_{\mathbf{G},\mathrm{div}} \in V_- \otimes V_+ \quad \text{and} \quad \mathrm{Rc}^-_{\mathbf{G},\mathrm{div}} \in V_+ \otimes V_-.$$

The first general definition of a Ricci tensor in generalized geometry was provided in [GF19], using torsion-free generalized connections. A simpler definition has been introduced in [ŠV20], which makes explicit use of the divergence operator. Both definitions agree on string algebroids, as recently proved by the second author jointly with R. Gonzalez Molina and J. Streets [GFS20, GFGM23]. We now recall an explicit characterisation of generalized Ricci tensors on string algebroids, as originally computed in [GF14] using Levi-Civita connections. A brief account of the motivation for the following formula will be given in Appendix A using the method from [GF14, GF19], for completeness. In the sequel, we will identify $V_\pm \cong V_\pm^*$ via the natural isomorphism provided by the pairing on $E$.

**Proposition 2.9.** *Let $(\mathbf{G}, \mathrm{div})$ be a pair given by a generalized metric $\mathbf{G}$ and a divergence operator $\mathrm{div}$ on a string algebroid $E$ over an $n$-manifold $M$. Define $\varepsilon \in \Gamma(E)$ by $\langle \varepsilon, \cdot\rangle := \mathrm{div}^\mathbf{G} - \mathrm{div}$. Via the isomorphism $E \cong T \oplus \mathrm{ad}P \oplus T^*$ provided by $\mathbf{G}$, we can uniquely write*

$$\varepsilon = \sigma_+(\zeta_+^\sharp) + z + \sigma_-(\zeta_-^\sharp)$$

*for $\zeta_\pm \in \Gamma(T^*)$ and $z \in \Gamma(\mathrm{ad}P)$. Then, one has*

$$\mathrm{Rc}^+_{\mathbf{G},\mathrm{div}}(a_-, b_+) = i_Y i_X \Big( \mathrm{Rc}_{\nabla^+} + \sum_j \langle i_{v_j} F_\theta, i_{v_j} F_\theta \rangle + \nabla^+ \zeta_+ \Big) \\ - i_Y \Big\langle d_\theta^* F_\theta + (-1)^n * (F_\theta \wedge *H) + i_{\zeta_+^\sharp} F_\theta, r \Big\rangle, \quad (2.16)$$

$$\mathrm{Rc}^-_{\mathbf{G},\mathrm{div}}(b_+, a_-) = i_X i_Y \Big( \mathrm{Rc}_{\nabla^-} + \sum_j \langle i_{v_j} F_\theta, i_{v_j} F_\theta \rangle - \nabla^- \zeta_- + \langle F_\theta, z \rangle \Big) \\ - i_Y \Big\langle d_\theta^* F_\theta + (-1)^n * (F_\theta \wedge *H) - d^\theta z - i_{\zeta_-^\sharp} F_\theta, r \Big\rangle. \quad (2.17)$$

*where $a_-, b_+$ are as in* (2.14) *and $\{v_j\}$ is an orthonormal frame for $g$.*

## 2.3 Killing spinors with parameter $\lambda$

In this section, we introduce a natural system of coupled equations, which provides a mild generalisation of the Killing spinor equations in generalized geometry [GFRT16, GF19]. In particular, these equations accommodate the Hull–Strominger system and the heterotic $G_2$ system as particular cases (Lemma 2.12). As we will see in Proposition 2.15, an important motivation for their study is that their solutions give special examples of generalized Ricci-flat metrics.

To introduce the equations, we fix a string algebroid $E$ over a spin manifold $M$ of dimension $n$. Given a pair $(\mathbf{G}, \mathrm{div})$, where $\mathbf{G}$ and $\mathrm{div}$ denote, respectively, a generalized metric and a divergence operator on $E$, the fixed spin structure on $M$, combined with the isometry

$$\sigma_+ \colon (T, g) \longrightarrow V_+ \\ X \longmapsto X + gX, \quad (2.18)$$

determines a real spinor bundle $S$ for $V_+$ (upon a choice of irreducible representation of the real Clifford algebra $\mathrm{Cl}(n, \mathbb{R})$). Associated to the pair $(\mathbf{G}, \mathrm{div})$, there are canonical first-order differential operators [GFRT16, GF19] (see also [CSCW11]):

$$D^S_- \colon \Omega^0(S) \to \Omega^0(V_-^* \otimes S), \quad \text{and} \quad \slashed{D}^+ \colon \Omega^0(S) \to \Omega^0(S). \quad (2.19)$$

The operator $D^S_-$ corresponds to the unique lift to $S$ of the metric-preserving operator $D^+_-$ in (2.10). The Dirac-type operator $\slashed{D}^+$ is more difficult to construct, as it involves torsion-free generalized connections (see proof of Proposition 2.18 and Appendix A). In our situation of interest, both operators can be described explicitly in terms of the affine metric connections with totally skew-symmetric torsion in (2.12). The formula for $D^S_-$ in the next result is, for instance, a direct consequence of the second equation in (2.13).

**Lemma 2.10** ([GFRT16]). *Let $(\mathbf{G}, \mathrm{div})$ be a generalized metric and a divergence operator on a string algebroid $E$, and let $(H, \theta)$ be the unique pair satisfying* (2.4) *determined by $\mathbf{G}$, where $H \in \Omega^3(M)$ and $\theta$ is a principal connection on $P$. Denote $\mathrm{div}^{\mathbf{G}} - \mathrm{div} = \langle \varepsilon, \cdot \rangle$, set $\zeta = g(\pi \varepsilon_+, \cdot) \in T^*$, and identify $S$ with a spinor bundle for $(T, g)$, via the isometry* (2.18). *Then, for any spinor $\eta \in \Omega^0(S)$ and $a_- = X + r - g(X) \in \Gamma(V_-)$, one has*

$$D^S_{a_-} \eta = \nabla^+_X \eta - \langle F_\theta, r \rangle \cdot \eta,$$
$$\slashed{D}^+ \eta = \slashed{\nabla}^{1/3} \eta - \tfrac{1}{2} \zeta \cdot \eta,$$

*where $\nabla^{1/3}_X Y = \nabla^g_X Y + \tfrac{1}{6} g^{-1} H(X, Y, \cdot)$ and $\slashed{\nabla}^{1/3}$ is the associated Dirac operator.*

We are ready to introduce our first system of equations of interest.

**Definition 2.11.** *Let $E$ be a string algebroid over a spin manifold $M$, and fix a constant $\lambda \in \mathbb{R}$. A triple $(\mathbf{G}, \mathrm{div}, \eta)$, given by a generalized metric $\mathbf{G}$, a divergence operator $\mathrm{div}$, and a spinor $\eta \in \Omega^0(S)$, is a solution of the* Killing spinor equations with parameter $\lambda$, *if*

$$D^S_- \eta = 0, \\ \slashed{D}^+ \eta = \lambda \eta. \quad (2.20)$$

*From their origins in theoretical physics, we will refer to the first equation in* (2.20) *as the* gravitino equation, *and to the second equation in* (2.20) *as the* dilatino equation.

The previous definition can be generalized, in the obvious way, to complex spinors. In even dimensions, the system (2.20) for a complex pure spinor $\eta$ forces $\lambda = 0$, since $\eta$ is necessarily chiral and the Dirac operator $\slashed{D}^+$ changes chirality. In odd dimensions, or for non-pure spinors, equations (2.20) are more general than the ones considered in [GFRT16, GF19]. The next result provides an explicit characterisation of the Killing spinor equations with parameter $\lambda$, in terms of classical tensors. The proof is a straightforward consequence of Lemma 2.10 and the classification of string algebroids in [GFRT20b, Proposition A.6], almost identical to the proof of [ACDAdLHGF24, Lemma 3.8], and it is therefore omitted.

**Lemma 2.12.** *Let $E$ be a string algebroid over a spin manifold $M$. Let $(\mathbf{G}, \mathrm{div}, \eta)$ be as in Lemma 2.10, and consider the associated data $(g, H, \theta, \eta, \zeta)$. Then the following holds.*

1. *$(\mathbf{G}, \eta)$ solves the gravitino equation in (2.20) if and only if $(g, H, \theta, \eta)$ solves*

$$\nabla^+ \eta = 0, \qquad F_\theta \cdot \eta = 0. \tag{2.21}$$

2. *$(\mathbf{G}, \mathrm{div}, \eta)$ solves the dilatino equation with parameter $\lambda$ in (2.20) if and only if $(g, H, \theta, \eta, \zeta)$ solves*

$$\left(\slashed{\nabla}^{1/3} - \frac{1}{2}\zeta\right) \cdot \eta = \lambda \eta. \tag{2.22}$$

*Conversely, any solution $(g, H, \theta, \eta)$ of (2.21) (resp. $(g, H, \theta, \eta, \zeta)$ of (2.22)) satisfying the heterotic Bianchi identity (2.4) determines a string algebroid as in Definition 2.2, endowed with a solution of the gravitino equation (resp. dilatino equation) in (2.20).*

*Remark* 2.13. In the mathematical physics literature, the name *gravitino equation* is often reserved for the first equation in (2.21), while the second receives the name of *gaugino equation*, referring to the superpartners of the graviton field and the gauge field, respectively (see e.g. [II05]). The unified treatment of the equations (2.21) is motivated by the way they appear in generalized geometry via the weak Koszul formula (2.11) and is borrowed from [ACDAdLHGF24]. ◯

We next prove a first structural property of the Killing spinor equations (2.20) with parameter $\lambda$ in relation to generalized Ricci-flat metrics, which motivates Definition 2.11. This is based on an interesting formula for the generalized Ricci tensor in terms of operators $D^S_-$ and $\slashed{D}^+$, discovered in the physics literature [CSCW11] (without proof) and first established in [GF19, Lemma 4.7].

**Lemma 2.14** ([GF19]). *Let $(\mathbf{G}, \mathrm{div})$ be a pair given by a generalized metric and a divergence operator on a string algebroid $E$ over a spin manifold $M$. Then, for any $a_- \in \Omega^0(V_-)$ and any spinor $\eta \in \Omega^0(S)$, the generalized Ricci tensor $\mathrm{Rc}^+_{\mathbf{G}, \mathrm{div}} \in V_- \otimes V_+$ satisfies*

$$\langle a_-, \mathrm{Rc}^+_{\mathbf{G}, \mathrm{div}} \rangle \cdot \eta = 4\Big(\slashed{D}^+ D^S_{a_-} - D^S_{a_-} \slashed{D}^+ - \sum_{j=1}^n v_j \cdot D^S_{\pi_-[v_j, a_-]}\Big)\eta,$$

*where $\{v_j\}$ is any choice of local orthonormal frame for $V_+$.*

As a direct consequence of the previous formula, any solution of the Killing spinor equations with parameter $\lambda$ is generalized Ricci-flat.

**Proposition 2.15.** *Let $(\mathbf{G}, \mathrm{div}, \eta)$ be a solution of the Killing spinor equations (2.20) with parameter $\lambda \in \mathbb{R}$, on a string algebroid $E$ over a spin $n$-manifold $M$. Provided that $\eta$ is nowhere-vanishing on $M$, the pair $(\mathbf{G}, \mathrm{div})$ satisfies*

$$\mathrm{Rc}^+_{\mathbf{G}, \mathrm{div}} = 0.$$

*More explicitly, in terms of the tuple $(g, H, \theta, \eta, \zeta)$ determined by $(\mathbf{G}, \mathrm{div})$, cf. Lemma 2.10, and a local orthonormal frame $\{v_j\}$ for $g$, one has*

$$\begin{aligned} \mathrm{Rc} - \frac{1}{4}H^2 + \sum_j \langle i_{v_j} F_\theta, i_{v_j} F_\theta \rangle + \frac{1}{2} L_{\zeta^\sharp} g &= 0, \\ d^* H - d\zeta + i_{\zeta^\sharp} H &= 0, \\ d^*_\theta F_\theta + (-1)^n * (F_\theta \wedge *H) + i_{\zeta^\sharp} F_\theta &= 0. \end{aligned} \tag{2.23}$$

*Proof.* Applying Lemma 2.14 to a solution of the Killing spinor equations (2.20), we have

$$\langle a_-, \mathrm{Rc}_{\mathbf{G},\mathrm{div}}^+\rangle \cdot \eta = -4 D_{a_-}^S \slashed{D}^+ \eta = -4\lambda D_{a_-}^S \eta = 0,$$

for every $a_- \in \Omega^0(V_-)$. Consequently,

$$|\langle a_-, \mathrm{Rc}_{\mathbf{G},\mathrm{div}}^+\rangle|^2 \eta = \langle a_-, \mathrm{Rc}_{\mathbf{G},\mathrm{div}}^+\rangle \cdot \langle a_-, \mathrm{Rc}_{\mathbf{G},\mathrm{div}}^+\rangle \cdot \eta = 0,$$

and therefore $|\langle a_-, \mathrm{Rc}_{\mathbf{G},\mathrm{div}}^+\rangle|^2 = 0$ for every section $a_-$, since $\eta$ is nowhere-vanishing. The first part of the proof follows now from the fact that the pairing on $V_+$ is positive-definite. As for the second part of the statement, equations (2.23) follow from the explicit formula for the generalized Ricci tensor (2.16), together with the unique decomposition of $\mathrm{Rc}_{\nabla^+}$ and $\nabla^+\zeta$ into symmetric and skew-symmetric 2-tensors, see e.g. [GFS20, IP01]:

$$\begin{aligned}
\mathrm{Rc}_{\nabla^+} &= \mathrm{Rc} - \frac{1}{4}H^2 - \frac{1}{2}d^*H, \\
\nabla^+\zeta &= \frac{1}{2}L_{\zeta^\sharp}g + \frac{1}{2}d\zeta - \frac{1}{2}i_{\zeta^\sharp}H.
\end{aligned} \tag{2.24}$$

$\square$

*Remark* 2.16. When $\zeta = d\phi$, for a smooth function $\phi$, equations (2.23) match the heterotic supergravity equations of motion for the metric, the 3-form flux, and the gauge field, in the mathematical physics literature, see e.g. [GF14, Mol24]. This suggests that solutions of (2.20) with closed, or even exact, one-form $\zeta$, play a special role; we explore this interesting aspect in Section 2.4. $\bigcirc$

## 2.4 Generalized scalar curvature

A comparison with the physics setup (see Remark 2.16) leads naturally to asking whether equations (2.20) imply the analogue of the equation of motion for the *dilaton field*. This is a scalar equation, given by the vanishing of the function

$$R_g - \frac{1}{2}|H|^2 + |F_\theta|^2 - 2d^*\zeta - |\zeta|^2 \in C^\infty(M), \tag{2.25}$$

where $(g, H, \theta, \eta, \zeta)$ is the associated data as in Lemma 2.10. This scalar quantity plays an important role in the theory of generalized Ricci flow, being closely related to the volume density of the generalized Perelman energy functional, see [GFS20, GFGMS24]. In (2.25), and in the sequel, we use the Hodge norm on differential forms, given by

$$|\beta|^2 \mathrm{Vol}_g = \beta \wedge *\beta = \frac{1}{k!} \sum_{i_1\ldots i_n=1}^n \beta_{i_1\ldots i_k} \beta^{i_1\ldots i_k}$$

for $\beta \in \Omega^k$. Note further that the summand $|F_\theta|^2$ in (2.25) is computed using the bilinear form $\langle \cdot, \cdot \rangle : \mathfrak{k} \otimes \mathfrak{k} \longrightarrow \mathbb{R}$, via

$$|F_\theta|^2 \mathrm{Vol}_g = \langle F_A \wedge *F_A \rangle,$$

and hence it might not be non-negative as a function on $M$.

To describe the dynamics of the dilaton field and provide an answer to the above question, we start by giving an interpretation of the scalar (2.25) in generalized geometry, using the operators (2.19). We build on a Lichnerowicz-type formula for the *cubic Dirac operator* $\slashed{\nabla}^{1/3}$ due to Bismut [Bis89], see also [AF03, Theorem 6.2]. We follow closely [GFS20, Proposition 3.39], see also [CSCW11] (alternative approaches can be found in [AMP24, ŠV20, SSCV24]). Given a pair $(\mathbf{G}, \mathrm{div})$, we can define a *rough Laplacian* operator

$$\Delta_-^S : \Omega^0(S) \to \Omega^0(S)$$

by the formula

$$\Delta_-^S \eta := \mathrm{tr}_{V_-}(D_-^- \otimes D_-^S)(D_-^S \eta),$$

where we recall that $D_-^S \eta \in \Omega^0(V_-^* \otimes S)$ and $D_-^-$ is the operator defined in Lemma A.3. It is not difficult to see that $\Delta_-^S$ is actually independent of the choice of $\mathbf{G}$-compatible torsion-free generalized connection with divergence div (see Appendix A), similarly as for $\slashed{D}^+$ (see [GF19, Lemma 3.4]), and hence it is a natural quantity associated canonically to the pair $(\mathbf{G}, \mathrm{div})$. We give next an explicit formula for $\Delta_-^S$.

**Lemma 2.17.** *Let $(\mathbf{G}, \operatorname{div})$ be given by a generalized metric and a divergence operator on a string algebroid $E$. Let $E \cong T \oplus \operatorname{ad}P \oplus T^*$ and $(H, \theta)$ be the isomorphism and the solution of (2.4) uniquely determined by $\mathbf{G}$, respectively. Define $\varepsilon \in \Omega^0(E)$ by $\langle \varepsilon, \cdot \rangle := \operatorname{div}^{\mathbf{G}} - \operatorname{div}$. We can uniquely write*

$$\varepsilon = \sigma_+(\zeta_+^\sharp) + z + \sigma_-(\zeta_-^\sharp)$$

*for $\zeta_\pm \in \Omega^1$ and $z \in \Omega^0(\operatorname{ad}P)$. Then, for any spinor $\eta \in \Omega^0(S)$, one has*

$$\Delta_-^S \eta = (\nabla^+)^* \nabla^+ \eta + \frac{1}{4} \langle F_\theta \wedge F_\theta \rangle \cdot \eta - \frac{1}{4} |F_\theta|^2 \eta - \nabla^+_{\zeta_-^\sharp} \eta + \langle F_\theta, z \rangle \cdot \eta, \qquad (2.26)$$

*where $\nabla^+$ is as in (2.12).*

*Proof.* Taking $a_-, c_- \in \Omega^0(V_-)$, we have

$$\langle (D_-^- \otimes D_-^S)(D_-^S \eta), a_- \otimes c_- \rangle := D_{a_-}^S D_{c_-}^S \eta - D_{D_{a_-} c_-}^S \eta.$$

To calculate the different elements in this formula explicitly, consider the natural isometries, cf. (2.7),

$$(T, g) \to V_+ : X \to X + gX, \qquad (T, -g) \oplus (\operatorname{ad}P, \langle \cdot, \cdot \rangle) \to V_- : X + r \to X + r - gX. \qquad (2.27)$$

Via the identification $T \oplus \operatorname{ad}P \cong V_-$, we have, cf. (A.2),

$$\begin{aligned}
D^-_{X+r}(Z+t) &= \nabla^{-1/3}_X Z - \tfrac{2}{3} g^{-1} \langle i_X F_\theta, t \rangle - \tfrac{1}{3} g^{-1} \langle i_Z F_\theta, r \rangle \\
&\quad + d^\theta_X t - \tfrac{2}{3} F_\theta(X, Z) - \tfrac{1}{3}[r, t] \\
&\quad + \tfrac{1}{\dim \mathfrak{k} + n - 1} \big((\langle r, t \rangle - g(X, Z))(z + \zeta_-^\sharp) - (\langle z, t \rangle - \zeta_-(Z))(X + r)\big).
\end{aligned}$$

Identifying $S$ with a spinor bundle for $(T, g)$, the Clifford bundle $\operatorname{Cl}(T)$ is defined via the relation (we follow [LM90])

$$X \cdot X = -g(X, X)$$

and, consequently, in a local orthonormal frame $\{e^j\}$ of $T$, the 2-forms $e^i \wedge e^j \in \mathfrak{so}(T) = \Lambda^2 T^*$ embed as $\frac{1}{2} e^i \cdot e^j$ into $\operatorname{Cl}(T)$, cf. [LM90, Proposition 6.2]. Hence, we have an identification of $D_-^S$ with the operator, cf. Lemma 2.7,

$$D^S_{X+r} \eta = \nabla^+_X \eta - \langle F_\theta, r \rangle \cdot \eta,$$

for any local spinor $\eta$. We choose a local orthogonal frame $\{v_\mu\}$ of $V_-$ and let $v^\mu$ denote the corresponding metric dual frame, so that $\langle v_\mu, v^\nu \rangle = \delta_{\mu\nu}$, which we assume without loss of generality to be of the form

$$v_\mu = \begin{cases} X_\mu & \text{if } 1 \leq \mu \leq n, \\ r_{\mu-n} & \text{if } n < \mu \leq \dim \mathfrak{k} + n, \end{cases} \qquad v^\mu = \begin{cases} -X_\mu & \text{if } 1 \leq \mu \leq n, \\ r^{\mu-n} & \text{if } n < \mu \leq \dim \mathfrak{k} + n, \end{cases}$$

where $X_\mu$ lie in $T$ and $r_\mu$ lie in $\operatorname{ad}P$, and analogously for their metric duals. Using this, we calculate

$$\begin{aligned}
\sum_{\mu=1}^n D^S_{v_\mu} D^S_{v^\mu} \eta &= \sum_{\mu=1}^n D^S_{X_\mu} \left(-\nabla^+_{X_\mu} \eta\right) + \sum_{\mu=1}^{\dim \mathfrak{k}} D^S_{r_\mu}(-\langle F_\theta, r^\mu \rangle \cdot \eta) \\
&= -\sum_{\mu=1}^n \nabla^+_{X_\mu} \nabla^+_{X_\mu} \eta + \sum_{\mu=1}^{\dim \mathfrak{k}} \langle F_\theta, r_\mu \rangle \langle F_\theta, r^\mu \rangle \cdot \eta \\
&= -\sum_{\mu=1}^n \nabla^+_{X_\mu} \nabla^+_{X_\mu} \eta + \frac{1}{16} \sum_{i,j,k,l=1}^n \langle F_{ij}, F_{kl} \rangle X_i X_j X_k X_l \cdot \eta \\
&= -\sum_{\mu=1}^n \nabla^+_{X_\mu} \nabla^+_{X_\mu} \eta + \frac{1}{4} \langle F_\theta \wedge F_\theta \rangle \cdot \eta - \frac{1}{4} |F_\theta|^2 \eta.
\end{aligned}$$

Note that the element

$$\Omega_\mathfrak{k} = \sum_{\mu=1}^{\dim \mathfrak{k}} [r_\mu, r^\mu] \in \mathfrak{k}$$

is independent of the choice of local frame $\{r_\mu\}$ and therefore

$$\Omega_\mathfrak{k} = \sum_{\mu=1}^{\dim \mathfrak{k}} [r^\mu, r_\mu] = -\Omega_\mathfrak{k} = 0.$$

Using this fact, we also have

$$\sum_{\mu=1}^{\dim \mathfrak{k}+n} D^-_{v_\mu} v_\mu = \sum_{\mu=1}^{n} D_{X_\mu}(-X_\mu) + \sum_{\mu=1}^{\dim \mathfrak{k}} D_{r_\mu} r^\mu$$

$$= \sum_{\mu=1}^{n} \left( -\nabla^{-1/3}_{X_\mu} X_\mu + \tfrac{2}{3} F_\theta(X_\mu, X_\mu) + \tfrac{1}{\dim \mathfrak{k}+n-1}(\delta_{\mu\mu}(z+\zeta^\sharp_-) - (\zeta_-(X_\mu))X_\mu) \right)$$

$$+ \sum_{\mu=1}^{\dim \mathfrak{k}} \left( -\tfrac{1}{3}[r_\mu, r^\mu] + \tfrac{1}{\dim \mathfrak{k}+n-1}(\delta_{\mu\mu}(z+\zeta^\sharp_-) - (\langle z, r^\mu\rangle)r_\mu) \right)$$

$$= -\sum_{\mu=1}^{n} \left( \nabla^g_{X_\mu} X_\mu \right) + \tfrac{1}{\dim \mathfrak{k}+n-1}(nz + (n-1)\zeta_-) + \tfrac{1}{\dim \mathfrak{k}+n-1}((\dim \mathfrak{k} - 1)z + \dim \mathfrak{k} \zeta^\sharp_-)$$

$$= -\sum_{\mu=1}^{n} \left( \nabla^g_{X_\mu} X_\mu \right) + \zeta^\sharp_- + z.$$

From the last formula

$$\sum_{\mu=1}^{n} D^S_{D_{v_\mu} v_\mu} \eta = -\sum_{\mu=1}^{n} \nabla^+_{\nabla^g_{X_\mu} X_\mu} \eta + \nabla^+_{\zeta^\sharp_-} \eta - \langle F_\theta, z \rangle \cdot \eta.$$

Using now [AF03, Theorem 6.1], which states

$$(\nabla^+)^* \nabla^+ \eta = -\sum_{\mu=1}^{n} \left( \nabla^+_{X_\mu} \nabla^+_{X_\mu} \eta + \nabla^+_{\nabla^g_{X_\mu} X_\mu} \eta \right),$$

we conclude

$$\Delta^S_- \eta = (\nabla^+)^* \nabla^+ \eta + \tfrac{1}{4} \langle F_\theta \wedge F_\theta \rangle \cdot \eta - \tfrac{1}{4} |F_\theta|^2 \eta - \nabla^+_{\zeta^\sharp_-} \eta + \langle F_\theta, z \rangle \cdot \eta. \qquad \square$$

We are ready to prove the main technical result of this section.

**Proposition 2.18.** *Let $(\mathbf{G}, \mathrm{div})$ be given by a generalized metric and a divergence operator on a string algebroid $E$. Consider the pair $(H, \theta)$ satisfying (2.4) uniquely determined by $\mathbf{G}$. Define $\varepsilon \in \Omega^0(E)$ by $\langle \varepsilon, \cdot \rangle := \mathrm{div}^{\mathbf{G}} - \mathrm{div}$. Via the isomorphism $E \cong T \oplus \mathrm{ad}P \oplus T^*$ provided by $\mathbf{G}$, we can uniquely write*

$$\varepsilon = \sigma_+(\zeta^\sharp) + z + \sigma_-(\zeta^\sharp_-)$$

*for $\zeta, \zeta_- \in \Omega^0(T^*)$ and $z \in \Omega^0(\mathrm{ad}P)$. Then, for any spinor $\eta \in \Omega^0(S)$, one has*

$$\left( \left( \slashed{D}^+ \right)^2 - \Delta^S_- - D^S_{\tilde{\varepsilon}_-} \right) \eta = \tfrac{1}{4}(\mathcal{S}^+ - 2d\zeta) \cdot \eta, \qquad (2.28)$$

*where $\tilde{\varepsilon}_- = \sigma_-(\zeta^\sharp + \zeta^\sharp_-) + z$ and*

$$\mathcal{S}^+ = R_g - \tfrac{1}{2}|H|^2 + |F_\theta|^2 - 2d^*\zeta - |\zeta|^2. \qquad (2.29)$$

*Proof.* For the proof, we use the same notation as in Lemma 2.17 and its proof. Via the natural isometries (2.27) we identify $S$ with a spinor bundle for $(T, g)$. Then, the operator $D^+_+$ in Lemma A.3 is given by the following affine metric connection in the tangent bundle

$$\tilde{\nabla}_X Y = \nabla^{1/3}_X Y + \tfrac{1}{n-1}(g(X,Y)\zeta^\sharp - \zeta(Y)X)$$

$$= \nabla^{1/3}_X Y + \tfrac{1}{n-1} g^{-1}(X^\flat \wedge \zeta)(Y)$$

and the operator $\slashed{D}^+$, defined as the Dirac operator for $D^+_+$ (see [GF19, Lemma 3.4]), is therefore $\slashed{D}^+ \eta = \tilde{\slashed{\nabla}} \eta$. Hence, given a local spinor $\eta$, we have

$$\tilde{\nabla}_X \eta = \nabla^{1/3}_X \eta + \tfrac{1}{4(n-1)}(-\zeta \cdot X + X \cdot \zeta) \cdot \eta.$$

Moreover, writing $\zeta = \sum_k \zeta_k e^k$, we have

$$\tilde{\nabla}\eta = \nabla^{1/3}\eta + \frac{1}{4(n-1)}\sum_{j,k}\zeta_k e^j \cdot (-e^k \cdot e^j + e^j \cdot e^k) \cdot \eta$$

$$= \nabla^{1/3}\eta + \frac{1}{4(n-1)}\sum_{j,k}\zeta_k(2\delta_{jk}e^j - 2e^k) \cdot \eta$$

$$= \nabla^{1/3}\eta - \frac{1}{2}\zeta \cdot \eta.$$

With these preliminaries, following the proof of [GFS20, Proposition 3.39], we compute

$$\tilde{\nabla}^2\eta = \left(\nabla^{1/3}\right)^2\eta - \frac{1}{2}\sum_j\left(e_j \cdot \nabla^g_{e_j}\zeta + \frac{1}{6}e_j \cdot H(e_j, \zeta^\sharp, \cdot)\right) \cdot \eta + \nabla^+_{\zeta^\sharp}\eta - \frac{1}{3}i_{\zeta^\sharp}H \cdot \eta - \frac{1}{4}|\zeta|^2\eta,$$

where we have used that $e^j \cdot \alpha + \alpha \cdot e^j = -2\alpha_j$, for any $\alpha = \sum_j \alpha_j e^j \in T^*$. Now, for any $\alpha, \beta \in T^*$, we have (cf. [LM90, Proposition 3.9])

$$(\alpha \wedge \beta) \cdot \eta = (\alpha \cdot \beta) \cdot \eta + (\alpha \lrcorner \beta) \cdot \eta,$$

and hence

$$\sum_j e_j \cdot \nabla^g_{e_j}\zeta \cdot \eta = \sum_j(e_j \wedge \nabla^g_{e_j}\zeta) \cdot \eta - \sum_j(e_j \lrcorner \nabla^g_{e_j}\zeta) \cdot \eta = (d\zeta + d^*\zeta) \cdot \eta.$$

Moreover,

$$\frac{1}{3}i_{\zeta^\sharp}H \cdot \eta = \frac{1}{6}\sum_{j,k,l}\zeta_l H_{jkl}e^j \wedge e^k \cdot \eta = \frac{1}{12}\sum_{j,k,l}\zeta_l H_{jkl}e^j \cdot e^k \cdot \eta = -\frac{1}{12}\sum_k e_j \cdot H(e_j, \zeta, \cdot) \cdot \eta.$$

We deduce that

$$\tilde{\nabla}^2\eta = \left(\nabla^{1/3}\right)^2\eta - \frac{1}{2}d\zeta \cdot \eta - \frac{1}{2}d^*\zeta\eta + \nabla^+_{\zeta^\sharp}\eta - \frac{1}{4}|\zeta|^2\eta.$$

Applying now Lemma 2.17 and using the Lichnerowicz-type formula [Bis89], cf. [AF03, Theorem 6.2],

$$\left(\nabla^{1/3}\right)^2 - (\nabla^+)^*\nabla^+ = \frac{1}{4}R_g - \frac{1}{8}|H|^2 + \frac{1}{4}dH$$

we conclude, applying the Bianchi identity (2.4),

$$\left((\slashed{D}^+)^2 - \Delta^S_-\right) \cdot \eta = \frac{1}{4}\left(R_g - \frac{1}{2}|H|^2 + |F_\theta|^2 - 2d^*\zeta - |\zeta|^2\right)\eta + \frac{1}{4}dH \cdot \eta - \frac{1}{2}d\zeta \cdot \eta + \nabla^+_{\zeta^\sharp}\eta$$

$$- \frac{1}{4}\langle F_\theta \wedge F_\theta\rangle \cdot \eta + \nabla^+_{\zeta^\sharp_-}\eta - \langle F_\theta, z\rangle \cdot \eta$$

$$= \frac{1}{4}\left(R_g - \frac{1}{2}|H|^2 + |F_\theta|^2 - 2d^*\zeta - |\zeta|^2\right)\eta - \frac{1}{2}d\zeta \cdot \eta + \nabla^+_{\zeta^\sharp + \zeta^\sharp_-}\eta - \langle F_\theta, z\rangle \cdot \eta. \quad \square$$

The previous result motivates the following definition of the *generalized scalar curvature* of a pair $(\mathbf{G}, \mathrm{div})$ on a string algebroid $E$.

**Definition 2.19.** *Let* $(\mathbf{G}, \mathrm{div})$ *be given by a generalized metric and a divergence operator on a string algebroid $E$. Define $\varepsilon \in \Omega^0(E)$ by $\langle \varepsilon, \cdot \rangle := \mathrm{div}^{\mathbf{G}} - \mathrm{div}$. Via the isomorphism $E \cong T \oplus \mathrm{ad}P \oplus T^*$ provided by $\mathbf{G}$, we can uniquely write*

$$\varepsilon = \sigma_+(\zeta^\sharp) + z + \sigma_-(\zeta^\sharp_-)$$

*for $\zeta, \zeta_- \in \Omega^0(T^*)$ and $z \in \Omega^0(\mathrm{ad}P)$. The* generalized scalar curvature

$$\mathcal{S}^+ = \mathcal{S}^+_{\mathbf{G}, \mathrm{div}} \in C^\infty(M)$$

*of the pair $(\mathbf{G}, \mathrm{div})$ is defined by*

$$\mathcal{S}^+\eta = 4\left(\left(\slashed{D}^+\right)^2 - \Delta^S_- - D^S_{\tilde{\varepsilon}_-} + \frac{1}{2}d\zeta\right) \cdot \eta, \tag{2.30}$$

*for any spinor $\eta \in \Omega^0(S)$, where $\tilde{\varepsilon}_- = \sigma_-(\zeta^\sharp + \zeta^\sharp_-) + z$. This is well-defined, and explicitly given by (2.29), by Proposition 2.18.*

*Remark* 2.20. As we will see shortly, the generalized scalar curvature plays a distinguished role when $\varepsilon \in \Omega^0(\ker \pi)$ and $[\varepsilon, \cdot] = 0$; in other words, when $\varepsilon$ gives a symmetry of the Dorfman bracket lying in the kernel of the anchor map. In this case, one has $\zeta = -\zeta_-$ and consequently $\tilde{\varepsilon}_- = z$. Furthermore, the condition $[\varepsilon, \cdot] = 0$ implies

$$d\zeta + 2\langle F_\theta, z \rangle = 0, \qquad d^\theta z = 0, \qquad [z, \cdot] = 0.$$

A particularly interesting instance arises when $z = 0$, that is, when $\varepsilon$ lies on the cotangent subbundle $T^* \subset E$, as in this case $d\zeta = 0$ and $\tilde{\varepsilon}_- = 0$. An important fact about $\mathcal{S}^+_{\mathbf{G},\mathrm{div}}$, which we will see in the proof of Proposition 2.23, is that it does *not* always coincide with the trace of the symmetric part of the generalized Ricci tensor $\mathrm{Rc}^+_{\mathbf{G},\mathrm{div}}$ (cf. [GFS20, Remark 3.42]). ◯

*Remark* 2.21. We can regard (2.30) as a local formula on $M$, so that there is no obstruction for the existence of the spinor bundle. Therefore, we can define the generalized scalar curvature of a pair $(\mathbf{G}, \mathrm{div})$ by (2.29), for a string algebroid over an arbitrary smooth manifold. ◯

We finish this section establishing the desired relation between the generalized Ricci-flat condition, the Killing spinor equations (2.20), and the generalized scalar curvature. The next general result about generalized Ricci metrics was proved by Gonzalez Molina in [Mol24, Proposition 6.4.5].

**Proposition 2.22** ([Mol24])**.** *Let $(\mathbf{G}, \mathrm{div})$ be a pair given by a generalized metric and a divergence operator on a string algebroid $E$ over an $n$-manifold $M$. Denote $\mathrm{div}^\mathbf{G} - \mathrm{div} = \langle \varepsilon, \cdot \rangle$, and set $\zeta = g(\pi \varepsilon_+, \cdot) \in T^*$. Then, assuming that $\mathrm{Rc}^+_{\mathbf{G},\mathrm{div}} = 0$, cf. (2.23), one has*

$$d\mathcal{S}^+ = (-1)^n * (d\zeta \wedge *H). \tag{2.31}$$

*In particular, if $d\zeta = 0$, the generalized scalar curvature of $(\mathbf{G}, \mathrm{div})$ is constant and furthermore one has*

$$d\left(|H|^2 - |F_\theta|^2 - d^*\zeta - |\zeta|^2\right) = 0 \tag{2.32}$$

Formula (2.31) follows from an explicit calculation in local coordinates while the proof of (2.32) follows by subtracting $\mathcal{S}^+$ minus the trace of the symmetric tensor in the generalized Ricci-flat equations (2.23), cf. Remark 2.20. In the last result of this section we establish that solutions of the Killing spinor equations have constant generalized scalar curvature.

**Proposition 2.23.** *Provided that $(\mathbf{G}, \mathrm{div}, \eta)$ is a solution of the Killing spinor equations (2.20) with parameter $\lambda \in \mathbb{R}$, then the generalized scalar curvature satisfies*

$$(\mathcal{S}^+ - 4\lambda^2 - 2d\zeta) \cdot \eta = 0,$$

*where $\mathrm{div}^\mathbf{G} - \mathrm{div} = \langle \varepsilon, \cdot \rangle$ and $\zeta = g(\pi \varepsilon_+, \cdot) \in T^*$. In particular, if $d\zeta = 0$ and $\eta$ is nowhere-vanishing, one has*

$$\mathcal{S}^+ = |H|^2 - |F_\theta|^2 - d^*\zeta - |\zeta|^2 = 4\lambda^2.$$

*Proof.* Combining the Killing spinor equations (2.20) with (2.28), we conclude immediately

$$(\mathcal{S}^+ - 2d\zeta)\eta = 4\left(\left(\slashed{D}^+\right)^2 - \Delta^S_- - D^S_{\tilde{\varepsilon}_-}\right)\eta = 4\lambda^2 \eta.$$

The last part of the statement follows, as in the proof of Proposition 2.22, by subtracting $\mathcal{S}^+$ minus the trace of the first equation in (2.23). □

*Remark* 2.24. Equation (2.31) implies that, for a Ricci-flat pair $(\mathbf{G}, \mathrm{div})$ with $d\zeta = 0$, the *dilaton equation of motion* $\mathcal{S}^+ = 0$ is satisfied, up to an overall constant on the manifold. Furthermore, for any solution of the Killing spinor equations (2.20) with $d\zeta = 0$ and parameter $\lambda \neq 0$, one has $\mathcal{S}^+ > 0$. ◯

*Remark* 2.25. For a solution $(\mathbf{G}, \mathrm{div}, \eta)$ of the Killing spinor equations in dimension 6 with $\eta \neq 0$, one has that $\lambda = 0$ (because $\eta$ is pure). Imposing further that $\zeta = d\phi$, one has that $(\mathbf{G}, \mathrm{div}, \eta)$ is equivalent to a solution of the Hull-Strominger system [GFRT16] and the previous result implies that $\mathcal{S}^+ = 0$ in this case. ◯

## 2.5 Coupled instantons

We now introduce a second natural system of coupled equations on a string algebroid $E$ over a spin manifold $M$ of arbitrary dimensions. This is inspired by recent developments on *coupled instantons*, both in the physics [dlOLS18a, dlOLS18b] and mathematical literature [GFGM23, GFJS23]. These equations are closely related to Killing spinors and generalized Ricci-flat metrics (see Problem 1 and Problem 2 below), and play an important role in recent developments around the Hull–Strominger system and non-Kähler mirror symmetry [ACDAdLHGF24].

To introduce our equations of interest, we fix a string algebroid $E$ over an oriented spin manifold $M$. Given a generalized metric $\mathbf{G}$ on $E$, we note that the operator $D_+^-$ in (2.10) has a natural curvature endomorphism

$$F_{D_+^-} \in \Lambda^2 V_+ \otimes \Lambda^2 V_-,$$

defined by

$$F_{D_+^-}(a_+, b_+)c_- = D_{a_+} D_{b_+} c_- - D_{b_+} D_{a_+} c_- - D_{\sigma_+ \pi[a_+, b_+]} c_-$$

for any sections $a_+, b_+$ of $V_+$ and $c_-$ of $V_-$.

**Definition 2.26.** *Let $E$ be a string algebroid over an oriented spin manifold $M^n$. A pair $(\mathbf{G}, \eta)$, given by a generalized metric $\mathbf{G}$ and a spinor $\eta \in \Omega^0(S)$, is a solution of the* coupled instanton equation, *if*

$$F_{D_+^-} \cdot \eta = 0. \tag{2.33}$$

*When $\eta$ is nowhere-vanishing, denoting by $G$ the stabilizer of $\eta$ in $\mathrm{Spin}(n)$, we will refer to a solution of (2.33) as a* coupled $G$-instanton.

We present next two problems, which relate the coupled instanton equation (2.33) to the Killing spinor equations (2.20) and generalized Ricci-flat metrics, providing an important motivation for their study. We first observe that, if $M^n$ is even-dimensional with $n = 2m$, then any solution of the gravitino equation in (2.20) (see also Lemma 2.12), with complex pure spinor $\eta$ and integrable complex structure is, in fact, a coupled $\mathrm{SU}(m)$-instanton, in the sense of Definition 2.26, cf. [GFGM23, Lemma 5.4]. In this case, solutions of the Killing spinor equations (2.20) with $\eta$ pure and $\zeta$ exact are in correspondence with solutions of the Hull–Strominger system on complex Calabi-Yau manifolds, with the Hermitian Yang–Mills Ansatz, see [GFGM23]. This motivates the following.

**Problem 1.** *Let $E$ be a string algebroid over an oriented spin manifold $M^n$. Let $(\mathbf{G}, \eta)$ be a solution of the gravitino equation in (2.20), i.e.*

$$D_-^S \eta = 0.$$

*Then, $(\mathbf{G}, \eta)$ satisfies the coupled instanton equation (2.33).*

More explicitly, the data $(E, \mathbf{G}, \eta)$ in the hypothesis of the previous Problem is equivalent to a tuple $(g, H, \theta, \eta)$ solving the equations, cf. Lemma 2.12,

$$\nabla^+ \eta = 0, \qquad F_\theta \cdot \eta = 0, \qquad dH - \langle F_\theta \wedge F_\theta \rangle = 0. \tag{2.34}$$

Furthermore, coupled $\mathrm{SU}(m)$-instantons which satisfy the gravitino equation in (2.20) with integrable complex structure are generalized Ricci-flat by [GFGM23, Proposition 5.6], for a suitable choice of divergence operator canonically determined by the Lee form of the $\mathrm{SU}(m)$-structure. This motivates our second problem.

**Problem 2.** *Let $E$ be a string algebroid over an oriented spin manifold $M^n$. Let $(\mathbf{G}, \eta)$ be a solution of the coupled instanton equation (2.33). Find the precise conditions, in terms of the $G$-structure determined by $\eta$, which imply that*

$$\mathrm{Rc}^+_{\mathbf{G}, \mathrm{div}_0} = 0,$$

*for a canonical choice of divergence operator $\mathrm{div}_0 = \mathrm{div}(\mathbf{G}, \eta)$ uniquely determined by the pair $(\mathbf{G}, \eta)$.*

*Remark* 2.27. Notice that the generalized Ricci tensor $\mathrm{Rc}^+_{\mathbf{G}, \mathrm{div}}$ only depends on $\zeta_+$ in the decomposition

$$\varepsilon = \sigma_+(\zeta_+^\sharp) + z + \sigma_-(\zeta_-^\sharp),$$

where $\langle \varepsilon, \cdot \rangle := \mathrm{div}^{\mathbf{G}} - \mathrm{div}$, for $\zeta_\pm \in \Gamma(T^*)$ and $z \in \Gamma(\mathrm{ad}P)$, see Proposition 2.9. We expect a solution of (2.21) to determine the divergence uniquely, by imposing the constraint $\varepsilon \in T^*$. ○

In Section 2.6 we present a complete answer to Problem 1 in arbitrary dimensions, via a spinorial proof, cf. Theorem 2.32. Section 3 and Section 4 are devoted to study those two Problems in the 7-dimensional case, where any nowhere-vanishing spinor determines a $G_2$-structure. In particular, by direct application of Theorem 2.32, in Theorem 4.6 we extend a result by De la Ossa, Larfors and Svanes in seven dimensions in the physics literature [dlOLS18a, dlOLS18b]. Before going into those specifics, we conclude this section by giving an explicit characterisation of the coupled instanton equation (2.33) in arbitrary dimensions, which we will use for the proof of Theorem 2.32. To that end, we consider a generalized metric $\mathbf{G}$ as in Definition 2.26, and observe that there exists an anti-isometry

$$\sigma_- \colon (T \oplus \mathrm{ad}P, \langle \cdot, \cdot \rangle^0) \longrightarrow V_- \tag{2.35}$$
$$X + r \longmapsto X + r - gX,$$

where $\langle X + r, X + r \rangle^0 = g(X, X) - \langle r, r \rangle$. Using this, and the isometry $\sigma_+ \colon T \to V_+$ from (2.8), the operator $D_-^+$ in (2.10) can be identified with the ordinary connection (see Lemma 2.7)

$$D = \begin{pmatrix} \nabla^- & \mathbb{F}^\dagger \\ -\mathbb{F} & d^\theta \end{pmatrix}, \tag{2.36}$$

where $\mathbb{F} \in \Omega^1(\mathrm{Hom}(T, \mathrm{ad}P))$ is the $\mathrm{Hom}(T, \mathrm{ad}P)$-valued 1-form

$$(i_X \mathbb{F})(Y) := F_\theta(X, Y) \tag{2.37}$$

and $\mathbb{F}^\dagger \in \Omega^1(\mathrm{Hom}(\mathrm{ad}P, T))$ is the corresponding $\langle \cdot, \cdot \rangle^0$-adjoint

$$(i_X \mathbb{F}^\dagger)(r) = -g^{-1} \langle i_X F_\theta, r \rangle^0.$$

We will use the standard notation $R_{\nabla^-}$ for the curvature of $\nabla^-$ and also $\nabla^{\theta,-}$ for the covariant derivative induced by $\theta$ and $\nabla^-$ on $\Lambda^2 T^* \otimes \mathrm{ad}P$. In particular,

$$(\nabla_Z^{\theta,-} F_\theta)(X, Y) = d_Z^\theta(F_\theta(X, Y)) - F_\theta(\nabla_Z^- X, Y) - F_\theta(X, \nabla_Z^- Y)$$

for any triple of vector fields $X, Y, Z$ on $M$. An explicit formula for the curvature of $D$ has been computed in [GFGM23, Lemma 4.7] as follows.

**Lemma 2.28** ([GFGM23]). *The curvature of $D$ is given by*

$$F_D = \begin{pmatrix} R_{\nabla^-} - \mathbb{F}^\dagger \wedge \mathbb{F} & -\mathbb{I}^\dagger \\ \mathbb{I} & [F_\theta, \cdot] - \mathbb{F} \wedge \mathbb{F}^\dagger \end{pmatrix}, \tag{2.38}$$

*where*

$$i_Y i_X \mathbb{F}^\dagger \wedge \mathbb{F}(Z) = g^{-1}\langle i_Y F_\theta, F_\theta(X, Z)\rangle - g^{-1}\langle i_X F_\theta, F_\theta(Y, Z)\rangle,$$
$$i_Y i_X \mathbb{I}(Z) = (\nabla_Z^{\theta,-} F_\theta)(X, Y) + F_\theta(X, g^{-1} i_Z i_Y H) - F_\theta(Y, g^{-1} i_Z i_X H),$$
$$i_Y i_X \mathbb{F} \wedge \mathbb{F}^\dagger(r) = F_\theta(Y, g^{-1}\langle i_X F_\theta, r\rangle) - F_\theta(X, g^{-1}\langle i_Y F_\theta, r\rangle).$$

The desired explicit characterisation of the coupled instanton equation (2.33) follows now from the previous formula.

**Proposition 2.29.** *Let $E$ be a string algebroid over an oriented spin manifold $M^n$. A solution $(\mathbf{G}, \eta)$ of the coupled instanton equation (2.33) on $E$ determines a tuple $(g, H, \theta, \eta)$, as in Lemma 2.10, solving the coupled instanton system*

$$(R_{\nabla^-} - \mathbb{F}^\dagger \wedge \mathbb{F}) \cdot \eta = 0,$$
$$\nabla^{\theta,+} F_\theta \cdot \eta = 0,$$
$$[F_\theta \cdot \eta, \cdot] - \mathbb{F} \wedge \mathbb{F}^\dagger \cdot \eta = 0, \tag{2.39}$$
$$dH - \langle F_\theta \wedge F_\theta \rangle = 0.$$

*Conversely, any solution of (2.39) determines a string algebroid as in Example 2.2, endowed with a solution of the coupled instanton equation (2.33).*

*Proof.* Note that

$$i_W i_V \mathbb{I}(Z) = (\nabla_Z^{\theta,-} F_\theta)(V, W) - F_\theta(V, g^{-1} H(Z, W)) + F_\theta(W, g^{-1} H(Z, V))$$
$$= d_Z^\theta(F_\theta(V, W)) - F_\theta(V, \nabla_Z^- W) - F_\theta(\nabla_Z^- V, W) - F_\theta(V, g^{-1} H(Z, W)) + F_\theta(W, g^{-1} H(Z, V))$$
$$= d_Z^\theta(F_\theta(V, W)) - F_\theta\left(V, \nabla_Z^g W - \frac{1}{2} g^{-1} H(Z, W)\right) - F_\theta\left(\nabla_Z^g V - \frac{1}{2} g^{-1} H(Z, V), W\right)$$
$$\quad - F_\theta(V, g^{-1} H(Z, W)) + F_\theta(W, g^{-1} H(Z, V))$$
$$= d_Z^\theta(F_\theta(V, W)) - F_\theta\left(V, \nabla_Z^g W + \frac{1}{2} g^{-1} H(Z, W)\right) - F_\theta\left(\nabla_Z^g V + \frac{1}{2} g^{-1} H(Z, V), W\right)$$
$$= (\nabla_Z^{\theta,+} F_\theta)(V, W).$$

The statement follows now directly from (2.38). □

*Remark* 2.30. More explicitly, choosing an orthonormal frame $\{v_j\}$ for $(T, g)$, the first three equations in (2.39) can be written as (cf. proof of Proposition 2.18)

$$\begin{aligned}(R_{\nabla^-}(v_i, v_j) - g^{-1}(\langle i_{v_j} F_\theta, F_\theta(v_i, \cdot)\rangle - \langle i_{v_i} F_\theta, F_\theta(v_j, \cdot)\rangle))v_i v_j \cdot \eta &= 0,\\ (\nabla^{\theta,+} F_\theta)(v_i, v_j) v_i v_j \cdot \eta &= 0,\\ ([F_\theta(v_i, v_j)), \cdot] - (F_\theta(v_j, g^{-1}\langle i_{v_i} F_\theta, \cdot\rangle) - F_\theta(v_i, g^{-1}\langle i_{v_j} F_\theta, \cdot\rangle))v_i v_j \cdot \eta &= 0.\end{aligned} \quad (2.40)$$

○

## 2.6 Gravitino solutions and instanton towers

In this section, we present a complete solution to Problem 1 in arbitrary dimensions. We also discuss a curious phenomenon which creates infinite numbers of instantons, with increasing rank, from solutions of the gravitino equation (2.21) with $dH = 0$. Concrete examples of these *instanton towers* are discussed in Section 4.3 in the seven-dimensional case, by application of Theorem 4.6. Key to our development is the following symmetry, originally due to Bismut [Bis89] (see also [GFS20, Proposition 3.21]), between the curvatures of the metric connections $\nabla^{\pm}$ with totally skew-symmetric torsion $\pm H \in \Omega^3$ defined in (2.12):

$$g(R_{\nabla^+}(X, Y)Z, W) = g(R_{\nabla^-}(Z, W)X, Y) + \frac{1}{2}dH(X, Y, Z, W). \quad (2.41)$$

We will also need the following basic technical lemma, which provides a generalization of [dlOLS18a, Lemma 5], originally proved in dimension 7, to arbitrary dimensions.

**Lemma 2.31.** *Let $(M^n, g)$ be an oriented spin manifold Riemannian manifold of dimension $n$. Let $\alpha, \beta \in \Omega^2$ a pair of two forms on $M$ and $\eta \in \Gamma(S)$ an arbitrary spinor. Then, we have*

$$[\alpha \cdot, \beta \cdot]\eta = \gamma \cdot \eta \quad (2.42)$$

*where, for any choice of local orthonormal frame $\{e_j\}$ for $g$, $\gamma \in \Omega^2(M)$ is defined by*

$$\gamma = \sum_{j=1}^n i_{e_j}\alpha \wedge i_{e_j}\beta.$$

*Proof.* Writing the 2-forms in the standard way $\alpha = \frac{1}{2!}\alpha_{ij}e^{ij}$ and $\beta = \frac{1}{2!}\beta_{ij}e^{ij}$ and using the canonical embedding of $\Omega^2$ into the Clifford algebra $e^j \wedge e^k = e^{jk} \mapsto \frac{1}{2}e_j e_k$, we have

$$\begin{aligned}\alpha \cdot (\beta \cdot \eta) &= \frac{1}{16}\alpha_{ij}\beta_{kl}e_i e_j e_k e_l \cdot \eta \\ \beta \cdot (\alpha \cdot \eta) &= \frac{1}{16}\alpha_{ij}\beta_{kl}e_k e_l e_i e_j \cdot \eta.\end{aligned} \quad (2.43)$$

Furthermore, $\gamma \in \Omega^2$ is given by $\gamma = i_{e_j}\alpha \wedge i_{e_j}\beta = \alpha_{jk}\beta_{jl}e^{kl} = \frac{1}{2}(\alpha_{jk}\beta_{jl} - \alpha_{jl}\beta_{jk})e^{kl}$ and, consequently,

$$\gamma \cdot \eta = \frac{1}{4}(\alpha_{jk}\beta_{jl} - \alpha_{jl}\beta_{jk})e_k e_l \cdot \eta. \quad (2.44)$$

The basic Clifford identity $e_i e_j = -e_j e_i - 2\delta_{ij}$ implies

$$\begin{aligned}e_k e_l e_i e_j &= -e_k e_i e_l e_j - 2\delta_{il}e_k e_j \\ &= e_k e_i e_j e_l + 2\delta_{jl}e_k e_i - 2\delta_{il}e_k e_j \\ &= -e_i e_k e_j e_l - 2\delta_{ik}e_j e_l + 2\delta_{jl}e_k e_i - 2\delta_{il}e_k e_j \\ &= e_i e_j e_k e_l + 2\delta_{jk}e_i e_l - 2\delta_{ik}e_j e_l + 2\delta_{jl}e_k e_i - 2\delta_{il}e_k e_j.\end{aligned}$$

Using the previous expression and substituting in the first equation in (2.43), we then have

$$\begin{aligned}\beta \cdot (\alpha \cdot \eta) &= \frac{1}{16}\alpha_{ij}\beta_{kl}\left(e_i e_j e_k e_l + 2\delta_{jk}e_i e_l - 2\delta_{ik}e_j e_l + 2\delta_{jl}e_k e_i - 2\delta_{il}e_k e_j\right) \cdot \eta \\ &= \alpha \cdot (\beta \cdot \eta) + \left(\frac{1}{8}\alpha_{ij}\beta_{jl}e_i e_l - \frac{1}{8}\alpha_{ij}\beta_{il}e_j e_l + \frac{1}{8}\alpha_{ij}\beta_{kj}\ e_k e_i - \frac{1}{8}\alpha_{ij}\beta_{ki}e_k e_j\right) \cdot \eta \\ &= \alpha \cdot (\beta \cdot \eta) + \left(-\frac{1}{4}\alpha_{ij}\beta_{il}e_j e_l + \frac{1}{4}\beta_{jk}\alpha_{ji}e_k e_i\right) \cdot \eta \\ &= \alpha \cdot (\beta \cdot \eta) - \gamma \cdot \eta,\end{aligned}$$

as desired. □

As in the previous section, we consider a string algebroid $E$ over an oriented spin manifold $M^n$. By application of Lemma 2.12, a solution $(\mathbf{G}, \eta)$ of the gravitino equation in (2.20) is equivalent to a tuple $(g, H, \theta, \eta)$ solving the equations (2.21) and the heterotic Bianchi identity (2.4). Similarly, by Proposition 2.29, the coupled instanton equation (2.33) on $E$ is equivalent to a tuple $(g, H, \theta, \eta)$ solving the coupled instanton system (2.39). With these preliminaries, our next result provides a complete solution to Problem 1 in arbitrary dimensions.

**Theorem 2.32.** *Let $P \to M$ be a principal $K$-bundle over an oriented spin manifold of arbitrary dimension. Then any solution $(g, H, \theta, \eta)$ of the equations*

$$\nabla^+ \eta = 0, \qquad F_\theta \cdot \eta = 0, \qquad dH - \langle F_\theta \wedge F_\theta \rangle = 0 \tag{2.45}$$

*solves the coupled instanton system (2.39), and consequently the connection on $T \oplus \mathrm{ad} P$ defined in (2.36) by*

$$D = \begin{pmatrix} \nabla^- & \mathbb{F}^\dagger \\ -\mathbb{F} & d^\theta \end{pmatrix}$$

*is an instanton with respect to $\eta$, i.e.*

$$F_D \cdot \eta = 0.$$

*In particular, given any solution $(g, H, \theta, \eta, \zeta)$ of the gravitino equation (2.21) and the dilatino equation (2.22), satisfying the heterotic Bianchi identity (2.4), the tuple $(g, H, \theta, \eta)$ solves the coupled instanton system (2.39).*

*Proof.* To prove that the first equation in (2.39) is satisfied with our hypothesis, we start by writing (using the summation convention, as we shall use throughout the proof):

$$\mathbb{F}^\dagger \wedge \mathbb{F} := \frac{1}{2} f^l{}_{kij} e^{ij} \otimes e^k \otimes e_l \in \Omega^2(M, \mathrm{End}(TM))$$

in a local orthonormal frame $\{e_j\}$ of $T$. The coefficients are computed as follows:

$$\begin{aligned}
f^l{}_{kij} &:= e^l(i_{e_j} i_{e_i} \mathbb{F}^\dagger \wedge \mathbb{F}(e_k)) \\
&= e^l\left(g^{-1}\langle i_{e_j} F_\theta, F_\theta(e_i, e_k)\rangle - g^{-1}\langle i_{e_i} F_\theta, F_\theta(e_j, e_k)\rangle\right) \\
&= e^l\left(g^{-1}\langle F_{jp} e^p, F_{ik}\rangle - g^{-1}\langle F_{ip} e^p, F_{jk}\rangle\right) = e^l\left((\langle F_j{}^p, F_{ik}\rangle - \langle F_i{}^p, F_{jk}\rangle) e_p\right) \\
&= \langle F_j{}^l, F_{ik}\rangle - \langle F_i{}^l, F_{jk}\rangle,
\end{aligned}$$

where $F_\theta := \frac{1}{2} F_{\mu\nu} e^{\mu\nu}$ and $F_{\mu\nu} := F_\theta(e_\mu, e_\nu) \in \Omega^0(\mathrm{ad} P)$. Note that

$$\langle F^l \wedge F_k \rangle = \langle F^l{}_i, F_{kj}\rangle e^{ij} = \frac{1}{2}\left(\langle F^l{}_i, F_{kj}\rangle - \langle F^l{}_j, F_{ki}\rangle\right) e^{ij},$$

i.e.

$$\langle F^l{}_j, F_{ki}\rangle - \langle F^l{}_i, F_{kj}\rangle = f^l{}_{kij} = -\langle F^l \wedge F_k\rangle_{ij}.$$

On the other hand,

$$\begin{aligned}
\langle F_\theta \wedge F_\theta \rangle &= \frac{1}{4}\langle F_{li}, F_{kj}\rangle e^{likj} \\
&= \frac{1}{12}\left(\langle F_{li}, F_{kj}\rangle - \langle F_{lk}, F_{ij}\rangle + \langle F_{lj}, F_{ik}\rangle\right) e^{likj} \\
&= \frac{1}{12}\left(\langle F_{li}, F_{kj}\rangle - \langle F_{lj}, F_{ki}\rangle - \langle F_{lk}, F_{ij}\rangle\right) e^{likj} \\
&= -\frac{1}{12}\left(-\langle F_l \wedge F_k\rangle_{ij} + \langle F_{lk}, F_{ij}\rangle\right) e^{likj},
\end{aligned}$$

which gives the components

$$-\langle F_l \wedge F_k\rangle_{ij} + \langle F_{lk}, F_{ij}\rangle = -\frac{1}{2}\langle F_\theta \wedge F_\theta\rangle_{likj}.$$

Combining the above and introducing the heterotic Bianchi identity $dH = \langle F_\theta \wedge F_\theta \rangle$, we therefore have

$$\begin{aligned}
f^l{}_{kij} &= -\langle F^l \wedge F_k\rangle_{ij} \\
&= \frac{1}{2}\langle F_\theta \wedge F_\theta\rangle^l{}_{kij} - \langle F_\theta, F^l{}_k\rangle_{ij} \\
&= -\frac{1}{2}(dH)^l{}_{kij} - \langle F_\theta, F^l{}_k\rangle_{ij}.
\end{aligned}$$

From formula (2.41), relating the curvatures of $\nabla^+$ and $\nabla^- = \nabla^g - \frac{1}{2}g^{-1}H$, we deduce:

$$\begin{aligned}
R_{\nabla^-} - \mathbb{F}^\dagger \wedge \mathbb{F} &= -\frac{1}{2}\left(f^l{}_{kij} - (R_{\nabla^-})^l{}_{kij}\right) e^{ij} \otimes e^k \otimes e_l \\
&= -\frac{1}{2}\left(-\frac{1}{2}\cancel{(dH)^l{}_{kij}} - \langle F_\theta, F^l{}_k\rangle_{ij} + \frac{1}{2}\cancel{(dH)^l{}_{kij}} - (R_{\nabla^+})_{ij}{}^l{}_k\right) e^{ij} \otimes e^k \otimes e_l \\
&= \frac{1}{2}\left(\langle F_\theta, F^l{}_k\rangle_{ij} + (R_{\nabla^+})_{ij}{}^l{}_k\right) e^{ij} \otimes e^k \otimes e_l.
\end{aligned}$$

Applying now the first equation in (2.45), we also have

$$g(R_{\nabla^+}(X,Y)e_i, e_j)e_ie_j \cdot \eta = \nabla^+_X\nabla^+_Y\eta - \nabla^+_Y\nabla^+_X\eta - \nabla^+_{[X,Y]}\eta = 0,$$

for arbitrary vector field $X, Y \in \Gamma(T)$. Hence, the first equation in (2.39) follows from

$$\left(R_{\nabla^-} - \mathbb{F}^\dagger \wedge \mathbb{F}\right) \cdot \eta = \frac{1}{4}\left(g(R_{\nabla^+}(e_l, e_k)e_i, e_j)e_ie_j \cdot \eta + 2\langle F_\theta \cdot \eta, F^l{}_k\rangle\right) \otimes e^k \otimes e_l = 0.$$

We next prove that the second and third equations in (2.39) are satisfied. Since $\nabla^+\eta = 0$ and $F_\theta \cdot \eta = 0$ by hypothesis, we have

$$(\nabla^{\theta,+}F_\theta) \cdot \eta = \nabla^{\theta,+}(F_\theta \cdot \eta) = 0, \quad \text{and} \quad [F_\theta \cdot \eta, \cdot] = 0,$$

so, it only remains to show that $\mathbb{F} \wedge \mathbb{F}^\dagger \cdot \eta = 0$. To see this, taking an orthonormal basis $\{\zeta_j\}$ for the Lie algebra $\mathfrak{k}$, write

$$\mathbb{F} \wedge \mathbb{F}^\dagger = \frac{1}{2}h^l{}_{kij}e^{ij} \otimes \zeta^k \otimes \zeta_l \in \Omega^2(M, \text{End}(\text{ad}P)),$$

where the coefficients are given by

$$\begin{aligned}
h^l{}_{kij} &= \zeta^l\left(i_{e_j}i_{e_i}\mathbb{F} \wedge \mathbb{F}^\dagger(\zeta_k)\right) = \zeta^l\left(F_\theta(e_j, g^{-1}\langle i_{e_i}F_\theta, \zeta_k\rangle) - F_\theta(e_i, g^{-1}\langle i_{e_j}F_\theta, \zeta_k\rangle)\right) \\
&= \zeta^l\left(F_\theta(e_j, \langle F_i{}^a e_a, \zeta_k\rangle) - F_\theta(e_i, \langle F_j{}^a e_a, \zeta_k\rangle)\right) \\
&= \zeta^l\left(F_\theta(e_j, \langle F^\alpha{}_i{}^a e_a \otimes \zeta_\alpha, \zeta_k\rangle) - F_\theta(e_i, \langle F^\alpha{}_j{}^a e_a \otimes \zeta_\alpha, \zeta_k\rangle)\right) \\
&= \zeta^l\left(F_\theta(e_j, F_{ki}{}^a e_a) - F_\theta(e_i, F_{kj}{}^a e_a)\right) = \zeta^l\left(F^\alpha{}_{ja}F_{ki}{}^a\zeta_\alpha - F^\alpha{}_{ia}F_{kj}{}^a\zeta_\alpha\right) \\
&= F^l{}_{ja}F_{ki}{}^a - F^l{}_{ia}F_{kj}{}^a.
\end{aligned}$$

Those are precisely the coefficients of the 2-form

$$\gamma^l{}_k := \sum_j (e_j \lrcorner F^l) \wedge (e_j \lrcorner F_k) \tag{2.46}$$

and hence the proof follows from $F_\theta \cdot \eta = 0$ by direct application of Lemma 2.31. $\square$

We discuss next the salient implications of the previous result for *instanton engineering* on the oriented spin $n$-manifold $M$. To start the iteration scheme, consider $M$ endowed with a metric $g$, a spinor $\eta$, a three-form $H \in \Omega^3$, and a connection $\theta$ on a principal $K$-bundle $P \to M$, satisfying

$$\nabla^+\eta = 0, \qquad F_\theta \cdot \eta = 0, \qquad dH = 0.$$

We consider the Lie algebra $\mathfrak{k}$ of the structure group $K$ of $P$ endowed with a biinvariant pairing $\langle,\rangle$. Using that $(d^\theta)^2 = [F_\theta, \cdot]$, it follows that the covariant derivative $\nabla^1 = d^\theta$ on the orthogonal bundle $V_1 = (\text{ad}P, \langle,\rangle)$ is an instanton with respect to $\eta$. Let $P_1$ be the principal bundle of split orthogonal frames of $V_1 \oplus V_1$, with structure group $K_1 = SO(r_1) \times SO(r_1)$, for $r_1 = \dim \mathfrak{k}$ (here we abuse notation for the special orthogonal group, since $\langle,\rangle$ may have arbitrary signature), and Lie algebra $\mathfrak{k}_1$ endowed with the neutral pairing

$$\langle \cdot, \cdot \rangle_1 = \text{tr}_{\mathfrak{so}(r_1)} - \text{tr}_{\mathfrak{so}(r_1)}.$$

The product connection $D^1 = \nabla^1 \times \nabla^1$ provides a new solution $(g, H, D^1, \eta)$ of the gravitino equation (2.21),

$$\nabla^+\eta = 0, \qquad F_{D^1} \cdot \eta = 0,$$

together with a trivial split solution of the Bianchi identity (cf. (2.4)):

$$dH = 0, \qquad \langle F_{D^1} \wedge F_{D^1}\rangle_1 = \text{tr}_{\mathfrak{so}(r_1)}(F_{\nabla^1} \wedge F_{\nabla^1}) - \text{tr}_{\mathfrak{so}(r_1)}(F_{\nabla^1} \wedge F_{\nabla^1}) = 0. \tag{2.47}$$

This type of ansatz for solving the (supersymmetry) equations is known in the supergravity literature as the *standard embedding*, cf. [GFRST22].

Applying now Theorem 2.32, from the data $(g, H, D^1, \eta)$ we can construct an instanton $\nabla^2$ on the bundle

$$V_2 = T \oplus \mathrm{ad} P_1$$

for the same metric $g$ and spinor $\eta$, explicitly given by (2.36) (with $\theta$ replaced by $D^1$). Note that $V_2$ is an orthogonal bundle with metric

$$\langle X + r, X + r \rangle_2 = g(X, X) - \langle r, r \rangle_1,$$

and the connection $\nabla^2$ is compatible with this metric. As before, let $P_2$ be the principal bundle of split orthogonal frames of $V_2 \oplus V_2$, with structure group $K_2 = SO(r_2) \times SO(r_2)$, for $r_2 = n + r_1(r_1 - 1)$, and Lie algebra $\mathfrak{k}_2$ endowed with the neutral pairing

$$\langle \cdot, \cdot \rangle_2 = \mathrm{tr}_{\mathfrak{so}(r_2)} - \mathrm{tr}_{\mathfrak{so}(r_2)}.$$

The product connection $D^2 = \nabla^2 \times \nabla^2$ provides a new solution $(g, H, D^2, \eta)$ of the gravitino equation (2.21) and a trivial split solution of the Bianchi identity (2.47). Iterating this scheme we obtain, by induction, an infinite tower of instantons with rank going to infinity on a fixed manifold. We summarise the construction in the next result.

**Proposition 2.33.** *Let $M$ be an oriented spin manifold of dimension $n$ endowed with a metric $g$, a three-form $H \in \Omega^3$, a spinor $\eta$, and a connection $\theta$ on a principal $K$-bundle $P \to M$, solving the equations, cf. (2.12),*

$$\nabla^+ \eta = 0, \qquad F_\theta \cdot \eta = 0, \qquad dH = 0.$$

*Then, there exists an infinite sequence of instantons $\{(V_k, \nabla^k)\}_{k \in \mathbb{N}}$ on $M$, for the same metric $g$ and spinor $\eta$, where $V_k$ is a real orthogonal bundle of rank*

$$r_k = n + r_{k-1}(r_{k-1} - 1), \qquad r_1 = \dim K,$$

*and $\nabla^k$ is a linear orthogonal connection on $V_k$.*

*Remark* 2.34. Given a solution $(g, H, \eta)$ of the equations

$$\nabla^+ \eta = 0, \qquad dH = 0,$$

the proof of Theorem 2.32 implies that the connection $\nabla^-$ on the orthogonal vector bundle $(T, g)$ is an instanton. This is a direct consequence of the identity (2.41). Hence, in this setup one can always choose $\theta = \nabla^-$ to start the iteration scheme. ◯

*Remark* 2.35. Note that the proof of Theorem 2.32 is local in nature and hence it applies, in particular, to solutions of the Hermitian-Yang-Mills equations on a Hermitian bundle $V$ over a complex manifold, provided that $c_1(V) = 0$. Consequently, by [GFJS23, Proposition 4.4] and the proof of Proposition 2.33 any pluriclosed Hermitian metric with vanishing Bismut-Ricci form leads to an instanton tower, cf. Remark 2.34. ◯

*Remark* 2.36. We speculate that this curious phenomenon holds for more general coupled instantons, not necessarily obtained via the standard embedding ansatz. ◯

# 3 Killing spinors and the heterotic $G_2$ system

## 3.1 $G_2$-structures

A $G_2$-structure always exists on an oriented spin manifold $M^7$, cf. [Gra69], and it is equivalent to a 3-form $\varphi \in \Omega^3(M)$ pointwise modelled on (we are following the conventions in [FKMS97] to define the form)

$$\varphi_0 = e^{127} + e^{347} + e^{567} + e^{135} - e^{146} - e^{236} - e^{245}. \tag{3.1}$$

Given that $G_2$ is a subgroup of $SO(7)$, $\varphi$ induces a metric via

$$g_\varphi(X, Y) = \frac{1}{6\mathrm{vol}_M} i_X \varphi \wedge i_Y \varphi \wedge \varphi, \tag{3.2}$$

where $X$ and $Y$ are vector fields on $M$, and $\mathrm{vol}_M$ is the volume form of $M$. We denote its Hodge dual under $* = *_\varphi$ by

$$\psi := *\varphi \in \Omega^4(M).$$

The pointwise model of $\psi$ is $\psi_0 \in \Lambda^4 \mathbb{R}^{7*}$ which is expressed as:
$$\psi_0 = e^{3456} + e^{1256} + e^{1234} - e^{2467} + e^{2357} + e^{1457} + e^{1367}. \tag{3.3}$$

With our conventions, we have $|\varphi|^2 = |\psi|^2 = 7$. Following [Kar08], we have the following natural decomposition of the space of forms on a manifold with $G_2$-structure $(M, \varphi)$.

**Proposition 3.1** (Decomposition of forms). *On a 7-manifold $M^7$ with $G_2$-structure $\varphi$, the spaces of differential forms $\Omega^k$ decompose orthogonally into irreducible $G_2$-representations as follows, where the $l$ in $\Omega_l^k$ denotes the rank of the subbundle in question.*

1. *The spaces of 0-forms and 1-forms are irreducible over $G_2$:*
$$\Omega^0 \cong \Omega_1^0, \qquad \Omega^1 \cong \Omega_7^1.$$

2. *The space of 2-forms decomposes as*
$$\Omega^2 = \Omega_7^2 \oplus \Omega_{14}^2,$$
*where*
$$\Omega_7^2 = \{X \lrcorner \varphi : X \in TM\} = \{\beta \in \Omega^2 : (\beta \lrcorner \varphi) \lrcorner \varphi = 3\beta\} = \{\beta \in \Omega^2 : \beta \lrcorner \psi = 2\beta\}, \tag{3.4}$$
$$\Omega_{14}^2 = \{\beta \in \Omega^2 : \beta \wedge \psi = 0\} = \{\beta \in \Omega^2 : \beta \lrcorner \varphi = 0\} = \{\beta \in \Omega^2 : \beta \lrcorner \psi = -\beta\}. \tag{3.5}$$

*Furthermore, we have the projection formulas for these spaces*
$$\pi_7(\beta) = \frac{1}{3}(\beta \lrcorner \varphi) \lrcorner \varphi = \frac{1}{3}(\beta + \beta \lrcorner \psi), \quad \pi_{14}(\beta) = \frac{1}{3}(2\beta - \beta \lrcorner \psi). \tag{3.6}$$

3. *The space of 3-forms decomposes as*
$$\Omega^3(M) = \Omega_1^3 \oplus \Omega_7^3 \oplus \Omega_{27}^3 \cong \Omega^0 \oplus \Omega^1 \oplus S_0^2(\mathbb{R}^7),$$
*where*
$$\Omega_1^3 = \{f\varphi : f \in C^\infty(M)\},$$
$$\Omega_7^3 = \{X \lrcorner \psi : X \in \Omega^1\},$$
$$\Omega_{27}^3 = \{\gamma \in \Omega^3 : \gamma \wedge \varphi = 0, \gamma \wedge \psi = 0\}.$$

4. *For $4 \leq k \leq 7$, the space $\Omega^k$ breaks into irreducible components using the isomorphism induced by the Hodge star operator $*\colon \Omega^k \cong \Omega^{7-k}$, as follows:*
$$\Omega^4 = \Omega_1^4 \oplus \Omega_7^4 \oplus \Omega_{27}^4, \qquad \Omega^5 = \Omega_7^5 \oplus \Omega_{14}^5; \quad \Omega^6 = \Omega_7^6, \qquad \Omega^7 = \Omega_1^7.$$

As first observed in [MCMS94], $G_2$-structures are classified by their torsion forms, cf. [Kar08, Theorem 2.23].

**Lemma 3.2** ([MCMS94]). *Let $(M,^7 \varphi)$ be a seven-dimensional manifold with $G_2$-structure $\varphi$. Then there are unique differential forms $\tau_0 \in \Omega^0$, $\tau_1 \in \Omega^1$, $\tau_2 \in \Omega_{14}^2$ and $\tau_3 \in \Omega_{27}^3$, called the torsion forms of $\varphi$, satisfying*
$$d\varphi = \tau_0 \psi + 3\tau_1 \wedge \varphi + *\tau_3, \tag{3.7}$$
$$d\psi = 4\tau_1 \wedge \psi + \tau_2 \wedge \varphi. \tag{3.8}$$

In the present work, *integrable* $G_2$-structures (such that $\tau_2 = 0$) play a distinguished role, as they admit compatible connections with totally skew-symmetric torsion.

**Proposition 3.3** ([FI03]). *Let $(M^7, \varphi)$ be a 7-manifold endowed with $G_2$-structure $\varphi$. Then, there exists an affine connection $\nabla$ on $T$ with totally skew-symmetric torsion preserving $\varphi$, that is, $\nabla\varphi = 0$, if and only if $\tau_2 = 0$. In this case, we say $\varphi$ is integrable, and the connection $\nabla$ is unique and given by*
$$\nabla = \nabla^g + \frac{1}{2}g^{-1}H_\varphi, \qquad H_\varphi = \frac{1}{6}\tau_0 \varphi - \tau_1 \lrcorner \psi - \tau_3, \tag{3.9}$$
*where $\nabla^g$ is the Levi-Civita connection of $g = g_\varphi$ and $H_\varphi \in \Omega^3(M)$ is the torsion of $\nabla$. The connection $\nabla$ is called the* characteristic connection *of $\varphi$.*

In particular, the compatibility of $\nabla$ with the $G_2$-structure implies that the endomorphism part of its curvature $R_\nabla$ lives in $\Omega_{14}^2 \subset \Omega^2$, that is, for any pair of vector fields $X, Y$ on $M$,
$$g(R_\nabla(X,Y)\cdot,\cdot) \in \Omega_{14}^2. \tag{3.10}$$

## 3.2 The heterotic $G_2$ system

In this section, we establish a relation between the Killing spinor equations (2.11) in seven dimensions and the heterotic $G_2$ system. We define the heterotic $G_2$ system following closely [dIOLMS20], on a fixed oriented and spin manifold $M^7$ endowed with a principal $K$-bundle $P$. We assume that $\mathfrak{k} = \mathrm{Lie}(K)$ is endowed with a non-degenerate bi-invariant symmetric bilinear form

$$\langle \cdot, \cdot \rangle : \mathfrak{k} \otimes \mathfrak{k} \longrightarrow \mathbb{R}.$$

We also assume that the first Pontryagin class of $P$, associated to $\langle \cdot, \cdot \rangle$ via Chern–Weil theory, is trivial:

$$p_1(P) = 0 \in H^4_{dR}(M, \mathbb{R}). \tag{3.11}$$

**Definition 3.4.** *Let $M^7$ be an oriented spin manifold, endowed with a principal $K$-bundle $P$, as before. A pair $(\varphi, \theta)$, where $\varphi$ is a $G_2$-structure on $M$ and $\theta$ is a principal connection on $P$, satisfies the* heterotic $G_2$ system *if there exists a constant $\lambda \in \mathbb{R}$ such that*

$$\begin{aligned} F_\theta \wedge \psi &= 0, \\ d\varphi &= \frac{12}{7}\lambda\psi + 3\tau_1 \wedge \varphi + *\tau_3, \\ d\psi &= 4\tau_1 \wedge \psi, \\ dH_\varphi - \langle F_\theta \wedge F_\theta \rangle &= 0, \end{aligned} \tag{3.12}$$

*where*

$$H_\varphi := \frac{2}{7}\lambda\varphi - \tau_1 \lrcorner \psi - \tau_3.$$

*Remark* 3.5. It is clear that there are no solutions to the final equation in (3.12) if $p_1(P) \neq 0$, so the assumption (3.11) is clearly necessary for the theory to be non-empty. ◯

In particular, the system (3.12) implies the following conditions on the torsion components:

$$\tau_0 = \frac{12}{7}\lambda \in \mathbb{R}, \qquad \tau_2 = 0.$$

Hence, by Proposition 3.3, the characteristic connection $\nabla$ of $\varphi$ exists, and it coincides with the metric connection with totally skew-symmetric torsion $H_\varphi$, that is,

$$\nabla = \nabla^+ := \nabla^g + \frac{1}{2}g^{-1}H_\varphi.$$

When the torsion component $\tau_1 \in \Omega^1(M)$ is exact, solutions of (3.12) provide compactifications of ten-dimensional heterotic string theory to three dimensions with $N = 1$ supersymmetry, in the supergravity limit. In this case, supersymmetry constrains the 3-dimensional geometry to be Anti de Sitter (AdS$_3$) or Minkowski and the natural scale for the 3-dimensional cosmological constant is set by $-\lambda^2$, cf. [dIOLMS20].

We now show that the heterotic $G_2$ system (3.12) is equivalent to the Killing spinor equations with parameter $\lambda$ on $M^7$, cf. (2.20). To see this, recall that a $G_2$-structure on $M$ is equivalent to a nowhere-vanishing spinor field $\eta \in \Omega^0(S)$, where $S = P_M \times_{G_2} \Delta_7$, for $P_M \to M$ the principal Spin(7)-bundle and $\Delta_7$ the irreducible real spinor representation, see Appendix B.1. More precisely, any nowhere-vanishing spinor $\eta$ defines a $G_2$-structure via

$$\varphi(X, Y, Z) := \langle X \cdot Y \cdot Z \cdot \eta, \eta \rangle. \tag{3.13}$$

Conversely, any positive three-form $\varphi$ determines a non-vanishing spinor via the identification

$$S \cong \langle \varphi \rangle \oplus \Omega^3_7 \tag{3.14}$$

provided by (B.6). To state the next result, we follow the notation in Lemma 2.12.

**Proposition 3.6.** *Let $E$ be a string algebroid over an oriented spin manifold $M^7$. Let $(\mathbf{G}, \mathrm{div}, \eta)$ be a solution of the Killing spinor equations with parameter $\lambda \in \mathbb{R}$ on $E$, cf. (2.20). Assume that the real spinor $\eta$ is non-zero and consider the tuple $(g, H, \theta, \zeta)$ determined by $(\mathbf{G}, \mathrm{div})$, cf. Lemma 2.10, as well as the $G_2$-structure $\varphi$ defined by $\eta$ via (3.13). Then $(\varphi, \theta)$ satisfies the heterotic $G_2$ system (3.12) and*

$$g = g_\varphi, \qquad H = H_\varphi := \frac{1}{6}\tau_0\varphi - \tau_1 \lrcorner \psi - \tau_3, \qquad \zeta = 4\tau_1, \qquad \lambda = \frac{7}{12}\tau_0. \tag{3.15}$$

*Furthermore, provided that $E = E_{P,H_0,\theta_0}$ as in Definition 2.2, one has*

$$\frac{1}{6}\tau_0\varphi - \tau_1 \lrcorner\, \psi - \tau_3 = H_0 + 2\langle a \wedge F_{\theta_0}\rangle + \langle a \wedge d^{\theta_0} a\rangle + \tfrac{1}{3}\langle a \wedge [a \wedge a]\rangle + db,$$

*for a uniquely determined $b \in \Omega^2$, where $a = \theta_0 - \theta \in \Omega^1(\mathrm{ad}P)$.*

*Conversely, any solution $(\varphi,\theta)$ of the heterotic $G_2$ system (3.12) determines a string algebroid $E_{P,H_\varphi,\theta}$ as in Definition 2.2, endowed with a solution $(\mathbf{G}, \mathrm{div}, \eta)$ of the Killing spinor equations with parameter $\lambda$, as in (2.20), and nowhere-vanishing spinor. The tuple $(g, H, \theta, \zeta)$ determined by $(\mathbf{G}, \mathrm{div})$ satisfies (3.15) and the spinor $\eta$ is given by (3.14), where we identify $S$, the spinor bundle for $V_+$, with a spinor bundle for $(T,g)$ via (2.18).*

*Proof.* By Lemma 2.12 and Remark 2.6, it suffices to prove the equivalence between solutions of the heterotic $G_2$ system (3.12) and solutions $(g, H, \theta, \zeta, \eta)$ of the coupled system defined by (2.21), (2.22), and the heterotic Bianchi identity (2.4). Given such a tuple $(g, H, \theta, \zeta)$, consider the $G_2$-structure $\varphi$ defined by the real spinor $\eta$ via (3.13). Note that, as $G_2 \leq SO(7)$, we have $g = g_\varphi$. Then, $\nabla^+ \eta = 0$ implies that $\nabla^+ \varphi = 0$ and hence, applying Proposition 3.3,

$$\tau_2 = 0, \qquad H = H_\varphi := \frac{1}{6}\tau_0\varphi - \tau_1 \lrcorner\, \psi - \tau_3.$$

Furthermore, by the identity (B.5), the condition $F_\theta \cdot \eta = 0$ is equivalent to $\theta$ being a $G_2$-instanton, that is,

$$F_\theta \wedge \psi = 0.$$

Using $\nabla^{1/3} = \nabla^+ - \frac{1}{3}H$, we have

$$\lambda\eta = \left(\slashed{\nabla}^{1/3} - \frac{1}{2}\zeta\right)\cdot\eta = \slashed{\nabla}^+\eta - \frac{1}{3}\slashed{H}\cdot\eta - \frac{1}{2}\zeta\cdot\eta = -\frac{1}{3}\slashed{H}\cdot\eta - \frac{1}{2}\zeta\cdot\eta.$$

On the other hand, applying (B.7) in Lemma B.2, it follows that

$$\slashed{H}\cdot\eta = \left(\frac{1}{6}\tau_0\varphi - \tau_1 \lrcorner\, \psi - \tau_3\right)\cdot\eta = -\frac{21}{12}\tau_0\cdot\eta - 6\tau_1\cdot\eta \tag{3.16}$$

and thus

$$\lambda\eta = -\frac{1}{3}\slashed{H}\cdot\eta - \frac{1}{2}\zeta\cdot\eta = \frac{7}{12}\tau_0\cdot\eta + 2\tau_1\cdot\eta - \frac{1}{2}\zeta\cdot\eta \Rightarrow \left(2\tau_1 - \frac{1}{2}\zeta\right)\cdot\eta + \left(\frac{7}{12}\tau_0 - \lambda\right)\cdot\eta = 0.$$

Using now that $(2\tau_1 - \frac{1}{2}\zeta)\cdot\eta \in \langle\eta\rangle^\perp \cong \Omega^1 \subset S$, cf. Appendix B.2, it follows that:

$$\zeta = 4\tau_1, \qquad \lambda = \frac{7}{12}\tau_0.$$

Conversely, given a solution $(\varphi,\theta)$ of the heterotic $G_2$ system (3.12), consider the real nowhere-vanishing spinor $\eta$ defined by (3.14). Then, by the third equation in (3.12), $\tau_2 = 0$ and hence

$$\nabla^+\varphi = 0 \Rightarrow \nabla^+\eta = 0,$$

where $\nabla^+ = \nabla^g + \frac{1}{2}g^{-1}H$ and $H = H_\varphi$, see (3.15). As before, $F_\theta \wedge \psi = 0$ implies $F_\theta \cdot \eta = 0$, while the second equation in (3.12) implies $\lambda = \frac{7}{12}\tau_0$, by Lemma 3.2. Finally, setting $\zeta = 4\tau_1$, we have

$$\left(\slashed{\nabla}^{1/3} - \frac{1}{2}\zeta\right)\cdot\eta = \left(\slashed{\nabla}^+ - \frac{1}{3}\slashed{H} - \frac{1}{2}\zeta\right)\cdot\eta = \frac{7}{12}\tau_0\cdot\eta + 2\tau_1\cdot\eta - \frac{1}{2}\zeta\cdot\eta = \frac{7}{12}\tau_0\cdot\eta = \lambda\cdot\eta. \qquad \square$$

*Remark* 3.7. Given a solution $(\varphi,\theta)$ of the heterotic $G_2$ system, the associated string algebroid is $E_{P,H_\varphi,\theta}$ as in Definition 2.2, and the corresponding generalized metric is

$$\mathbf{G}_\varphi = \begin{pmatrix} 0 & 0 & g_\varphi^{-1} \\ 0 & -\mathrm{Id} & 0 \\ g_\varphi & 0 & 0 \end{pmatrix},$$

with eigenbundles

$$V_+ = \{X + g_\varphi X : X \in T\}, \quad V_- = \{X + r - g_\varphi X : X \in T, r \in \mathrm{ad}P\}.$$

The divergence operator associated to a solution, given by

$$\mathrm{div} = \mathrm{div}^{\mathbf{G}_\varphi} - \langle\varepsilon, \cdot\rangle,$$

is uniquely determined provided that we impose the natural condition $\varepsilon \in T^*$. In this case, $\varepsilon = 8\tau_1$; that is, it coincides up to a constant multiple with the Lee form of the $G_2$-structure. $\bigcirc$

## 3.3 Curvature constraints on the heterotic $G_2$ system

We will derive various curvature constraints for solutions of the heterotic $G_2$ system. Our results follow from the characterisation of the system using generalized geometry, in Proposition 3.6, combined with Proposition 2.15 and Proposition 2.18.

We keep the notation of the previous section. In particular, we fix an oriented spin manifold $M^7$ endowed with a principal $K$-bundle $P$ with trivial first Pontryagin class, see (2.3). Our first result gives an interpretation of the heterotic $G_2$ system as a special class of generalized Ricci-flat metrics.

**Theorem 3.8.** *Given a solution $(\varphi, \theta)$ of the heterotic $G_2$-system (3.12) on $(M, P)$, the associated Riemannian metric $g = g_\varphi$ satisfies, for a local orthonormal frame $\{v_j\}$:*

$$\mathrm{Rc} - \frac{1}{4} H^2 + \sum_j \langle i_{v_j} F_\theta, i_{v_j} F_\theta \rangle + \frac{1}{2} L_{\zeta^\sharp} g = 0,$$
$$d^* H - d\zeta + i_{\zeta^\sharp} H = 0, \qquad (3.17)$$
$$d^*_\theta F_\theta - *(F_\theta \wedge *H) + i_{\zeta^\sharp} F_\theta = 0,$$

*where $H$ and $\zeta$ are uniquely determined by the torsion components of $\varphi$, via*

$$H = \frac{2}{7} \lambda \varphi - \tau_1 \lrcorner \psi - \tau_3, \qquad \zeta = 4\tau_1.$$

*Proof.* By Proposition 3.6, $(\varphi, \theta)$ determines a solution $(\mathbf{G}_\varphi, \mathrm{div}^\varphi, \eta_\varphi)$ of the Killing spinor equations with parameter $\lambda$ on the string algebroid $E_{P,H,\theta}$. More explicitly, the generalized metric $\mathbf{G}_\varphi$ is as in Remark 3.7 and

$$\mathrm{div}^\varphi = \mathrm{div}^{\mathbf{G}_\varphi} - \langle 8\tau_1, \cdot \rangle.$$

Applying now Proposition 2.15, we have $\mathrm{Rc}^+_{\mathbf{G}_\varphi, \mathrm{div}^\varphi} = 0$, and the result follows from (2.23). □

By the proof of the previous result, a solution $(\varphi, \theta)$ of the heterotic $G_2$ system determines a generalized Ricci-flat metric. Alternatively, we can think of $(\varphi, \theta)$ as solving the equations of motion of heterotic supergravity for the metric, the three-form flux, and the gauge field, in the mathematical physics literature, see Remark 2.16. In our next result, the analogue of the equation of motion for the *dilaton field* is satisfied up to an overall constant on the manifold, explicitly given in terms of the parameter $\lambda$ in (3.12). In other words, solutions of the heterotic $G_2$ system have constant generalized scalar curvature, proportional to the square of the torsion component $\tau_0$.

**Theorem 3.9.** *Given $(\varphi, \theta)$ a solution of heterotic $G_2$-system (3.12) on $(M, P)$, one has*

$$d\tau_1 \wedge *H = 0,$$
$$\mathcal{S}^+ = R_g - \frac{1}{2} |H|^2 + |F_\theta|^2 - 8 d^* \tau_1 - 16 |\tau_1|^2 = \frac{49}{36} \tau_0^2, \qquad (3.18)$$

*where $H = \frac{2}{7} \lambda \varphi - \tau_1 \lrcorner \psi - \tau_3$.*

*Proof.* As in the proof of Theorem 3.8, $(\varphi, \theta)$ determines a solution $(\mathbf{G}_\varphi, \mathrm{div}^\varphi, \eta_\varphi)$ of the Killing spinor equations with parameter $\lambda = \frac{7}{12} \tau_0$ on the string algebroid $E_{P,H,\theta}$. Applying now Proposition 2.18, we have

$$(\mathcal{S}^+ - 8 d\tau_1) \cdot \eta = 4 \left( \left( \slashed{D}^+ \right)^2 - \Delta^S_- - D^S_{\tilde{\varepsilon}_-} \right) \eta = 4 \lambda^2 \eta = \frac{49}{36} \tau_0^2 \eta,$$

where

$$\mathcal{S}^+ = R_g - \frac{1}{2} |H|^2 + |F_\theta|^2 - 8 d^* \tau_1 - 16 |\tau_1|^2. \qquad (3.19)$$

The proof follows from the orthogonal decomposition $S = \langle \eta \rangle \oplus \Omega^1$, combined with $d\tau_1 \cdot \eta \in \Omega^1$ (see Appendix B.2), which implies $\mathcal{S}^+ = \frac{49}{36} \tau_0^2$ and $d\tau_1 \cdot \eta = 0$, combined with formula (2.31) in Proposition 2.22. Note that the condition $d\tau_1 \cdot \eta = 0$ is equivalent to $d\tau_1 \wedge \psi = 0$, which is trivially satisfied since $d\psi = 4\tau_1 \wedge \psi$. □

We conclude this section with an alternative form of the scalar equation in (3.18).

**Corollary 3.10.** *Given $(\varphi, \theta)$ a solution of heterotic $G_2$-system (3.12) on $(M, P)$, one has*

$$\frac{7}{6} \tau_0^2 + 12 |\tau_1|^2 + 4 d^* \tau_1 - |\tau_3|^2 + |F_\theta|^2 = 0. \qquad (3.20)$$

*Proof.* Applying the result in Theorem 3.9 and the fact that

$$|H|^2 = H \lrcorner H = \frac{1}{6}\tau_0\varphi \lrcorner \frac{1}{6}\tau_0\varphi + (\tau_1 \lrcorner \psi) \lrcorner (\tau_1 \lrcorner \psi) + \tau_3 \lrcorner \tau_3 = \frac{7}{36}\tau_0^2 + 4|\tau_1|^2 + |\tau_3|^2,$$

where we have used $|\varphi|^2 = 7$ and $(\tau_1 \lrcorner \psi) \lrcorner \psi = -4\tau_1$. On the other hand, we have (cf. [Bry03, Equation (4.2)])

$$R_g = \frac{21}{8}\tau_0^2 + 30|\tau_1|^2 - \frac{1}{2}|\tau_3|^2 + 12d^*\tau_1.$$

which gives us

$$\begin{aligned}
\frac{49}{36}\tau_0^2 &= R_g - \frac{1}{2}|H|^2 - 8d^*\tau_1 - 16|\tau_1|^2 + |F_\theta|^2 \\
&= \frac{21}{8}\tau_0^2 + 30|\tau_1|^2 - \frac{1}{2}|\tau_3|^2 + 12d^*\tau_1 - \frac{1}{2}|H|^2 - 8d^*\tau_1 - 16|\tau_1|^2 + |F_\theta|^2 \\
&= \frac{21}{8}\tau_0^2 + 14|\tau_1|^2 - \frac{1}{2}|\tau_3|^2 + 4d^*\tau_1 - \frac{1}{2}|H|^2 + |F_\theta|^2 \\
&= \frac{21}{8}\tau_0^2 + 14|\tau_1|^2 - \frac{1}{2}|\tau_3|^2 + 4d^*\tau_1 - \frac{1}{2}\left(\frac{7}{36}\tau_0^2 + 4|\tau_1|^2 + |\tau_3|^2\right) + |F_\theta|^2 \\
&= \left(\frac{21}{8} - \frac{7}{72}\right)\tau_0^2 + 12|\tau_1|^2 + 4d^*\tau_1 - |\tau_3|^2 + |F_\theta|^2 \\
&= \frac{91}{36}\tau_0^2 + 12|\tau_1|^2 + 4d^*\tau_1 - |\tau_3|^2 + |F_\theta|^2
\end{aligned}$$

and the result follows. $\square$

# 4 The coupled $G_2$-instanton equations

## 4.1 Coupled instantons and the gravitino equation

We introduce the coupled $G_2$-instanton equations, a particular instance of the system (2.29) in seven dimensions. We will also establish the relation to the gravitino equation (2.21), by application of Theorem 2.32 in the present setup.

As in the previous section, we fix an oriented spin manifold $M^7$. Given a $G_2$-structure $\varphi$ on $M$ and a three-form $H \in \Omega^3$, we introduce the following quantity, which plays a similar role to the Bismut–Ricci form in the theory of coupled $SU(n)$-instantons, see [GFGM23, GFJS23]. Recall that the vector cross product $\times \colon T \otimes T \to T$ associated to $\varphi$ is defined by

$$g_\varphi(X \times Y, Z) = \varphi(X, Y, Z)$$

for any $X, Y \in T$.

**Definition 4.1** (Bismut–Ricci form)**.** *The* Bismut–Ricci form *associated to a pair* $(\varphi, H)$, *where $\varphi$ is a $G_2$-structure and $H \in \Omega^3(M)$, is the vector-valued 2-form*

$$\rho = \rho(\varphi, H) \in \Omega^2(M, T)$$

*defined by*

$$\rho(X, Y) := \frac{1}{2}\sum_j (R_{\nabla^+}(X,Y)e_j) \times e_j \tag{4.1}$$

*in terms of the vector cross product, where $\{e_j\}$ is a local orthonormal frame on $T$, and $\nabla^+$ is the metric connection with skew-symmetric torsion*

$$\nabla^+ = \nabla^g + \frac{1}{2}g^{-1}H.$$

*Remark* 4.2. An interesting special case of the previous definition follows when we take

$$H = H_\varphi := \frac{1}{6}\tau_0\varphi - \tau_1 \lrcorner \psi - \tau_3.$$

In this case, we say that $\rho_\varphi = \rho(\varphi, H_\varphi)$ is the *Bismut–Ricci form* of the $G_2$-structure. $\bigcirc$

As we observe in the next result, the Bismut–Ricci form is an obstruction to the integrability of the $G_2$-structure.

**Lemma 4.3.** *Assume that $(\varphi, H)$ satisfies $\nabla^+ \varphi = 0$. Then, $(\varphi, H)$ has vanishing Bismut–Ricci form:*

$$\rho(\varphi, H) = 0.$$

*Proof.* Assuming $\nabla^+ \varphi = 0$, the endomorphism part of the curvature tensor $R_{\nabla^+}$ lives in $\Omega^2_{14} \subset \Omega^2$, i.e. for any vector fields $X, Y$ on $M$,

$$g(R_{\nabla^+}(X,Y)\cdot, \cdot) \in \Omega^2_{14}.$$

The proof follows from the identity

$$g(\rho(X,Y), e_l) = \frac{1}{2} \sum_j g((R_{\nabla^+}(X,Y)e_j) \times e_j, e_l) = \frac{1}{2} \sum_{j,k} R_{\nabla^+}(X,Y)^k_j g(e_k \times e_j, e_l) = \frac{1}{2} \sum_{j,k} \varphi_{jkl} R_{\nabla^+}(X,Y)^k_j$$

and Proposition 3.1. Now, by Proposition 3.3, the assumption $\nabla^+ \varphi = 0$ implies that $\varphi$ is integrable and furthermore (3.9) holds, so

$$H = H_\varphi := \frac{1}{6} \tau_0 \varphi - \tau_1 \lrcorner \psi - \tau_3. \qquad \square$$

To introduce our equations of interest, we fix a principal $K$-bundle $P \to M$. The Lie algebra $\mathfrak{k} = \operatorname{Lie}(K)$ is endowed with a non-degenerate bi-invariant symmetric bilinear form $\langle \cdot, \cdot \rangle$, and we assume that $P$ has vanishing first Pontryagin class, cf. (3.11).

**Definition 4.4.** *Let $P \to M^7$ be a principal $K$-bundle over an oriented spin 7-manifold. For a triple $(\varphi, H, \theta)$, where $\varphi$ is a $G_2$-structure on $M$, $H \in \Omega^3(M)$, and $\theta$ is a principal connection on $P$, the* coupled $G_2$-instanton equation *is*

$$\begin{aligned}
\rho(\varphi, H) + \langle F_\theta, (F_\theta \lrcorner \varphi)^\sharp \rangle &= 0, \\
(\nabla^{\theta,+} F_\theta) \lrcorner \varphi &= 0, \\
[F_\theta \lrcorner \varphi, \cdot\,] - \mathbb{F} \wedge \mathbb{F}^\dagger \lrcorner \varphi &= 0, \\
dH - \langle F_\theta \wedge F_\theta \rangle &= 0,
\end{aligned} \qquad (4.2)$$

*where $\mathbb{F} \wedge \mathbb{F}^\dagger \in \Omega^2(\operatorname{End}(\operatorname{ad}P))$ is defined as in Lemma 2.28 by*

$$i_Y i_X \mathbb{F} \wedge \mathbb{F}^\dagger(r) = F_\theta(Y, g^{-1}\langle i_X F_\theta, r \rangle) - F_\theta(X, g^{-1}\langle i_Y F_\theta, r \rangle).$$

In the next result we establish a bijection between the solutions of the coupled $G_2$-instanton equation (4.2) and the coupled instanton equations formulated in terms of spinors, in Proposition 2.29. Recall that a $G_2$-structure on $M$ is equivalent to a nowhere-vanishing spinor field $\eta \in \Omega^0(S)$, via (3.13) and (3.14). Note that the system (2.39), introduced in Proposition 2.29, can be regarded as a system for tuples $(g, H, \theta, \eta)$.

**Proposition 4.5.** *Let $P \to M^7$ be a principal $K$-bundle over an oriented spin 7-manifold. Then, any solution of the coupled $G_2$-instanton equations (4.2) determines a solution $(g_\varphi, H, \theta, \eta_\varphi)$ of (2.39). Conversely, any solution $(g, H, \theta, \eta)$ of (2.39) determines a solution of the coupled $G_2$-instanton equations (4.2) of the form $(\varphi_\eta, H, \theta)$, where $\varphi_\eta$ is defined by (3.13).*

*Proof.* The equivalence between the second and third equations in (2.39) and (4.2) follows easily from (B.5). It remains therefore to prove the equivalence between the first equation in (4.2) and the first equation in (2.39), so long as the Bianchi identity $dH = \langle F_\theta \wedge F_\theta \rangle$ is satisfied. Arguing as in the last part of the proof of Theorem 2.32, the desired equivalence now follows from (B.5):

$$\begin{aligned}
(R_{\nabla^-} - \mathbb{F}^\dagger \wedge \mathbb{F}) \cdot \eta = 0 &\iff R_{\nabla^-} - \mathbb{F}^\dagger \wedge \mathbb{F} \in \Omega^2_{14} \\
&\iff (R_{\nabla^-} - \mathbb{F}^\dagger \wedge \mathbb{F}) \lrcorner \varphi = 0 \\
&\iff \frac{1}{2!1!} \left( \langle F_\theta, F^l{}_k \rangle_{ij} + (R_{\nabla^+})_{ij}{}^l{}_k \right) \varphi_{ijp} e^p \otimes e^k \otimes e_l = 0 \\
&\iff \langle (F_\theta \lrcorner \varphi)^\sharp, F^l{}_k e^k \otimes e_l \rangle + \rho^l{}_k = 0 \\
&\iff \langle (F_\theta \lrcorner \varphi)^\sharp, F_\theta \rangle + \rho = 0. \qquad \square
\end{aligned}$$

As a direct consequence of the previous result and Proposition 2.29, it follows that any solution of the coupled $G_2$-instanton equation (4.2) corresponds to a $G_2$-instanton on $T \oplus \operatorname{ad}P$, given by the connection (2.36).

To finish this section, we prove that any solution of the gravitino equation (2.21) in seven dimensions provides a solution of the coupled $G_2$-instanton equation (4.2), by application of Theorem 2.32. Note that, in the present setup, the gravitino equation is given by (see Lemma 2.12 and Appendix B.2)

$$\nabla^+ \varphi = 0, \qquad F_\theta \wedge \psi = 0, \tag{4.3}$$

where the unknowns are pairs $(\varphi, \theta)$ as before and $\nabla^+ = \nabla^g + \frac{1}{2}g^{-1}H_\varphi$ for (see Proposition 3.3)

$$H_\varphi := \frac{1}{6}\tau_0 \varphi - \tau_1 \lrcorner \psi - \tau_3. \tag{4.4}$$

We are mainly interested in solutions of the gravitino equation solving also the heterotic Bianchi identity (2.4)

$$dH_\varphi = \langle F_\theta \wedge F_\theta \rangle. \tag{4.5}$$

Our result provides an alternative proof of, and is inspired by, some results in [dlOLS18a, dlOLS18b].

**Theorem 4.6.** *Let $P \to M^7$ be a principal $K$-bundle over an oriented spin 7-manifold. Then any solution $(\varphi, \theta)$ of the gravitino equation (4.3) and the Bianchi identity (4.5) determines a solution $(\varphi, H_\varphi, \theta)$ of the coupled $G_2$-instanton equation (4.2), and the connection on $T \oplus \mathrm{ad} P$ defined in (2.36) by*

$$D = \begin{pmatrix} \nabla^- & \mathbb{F}^\dagger \\ -\mathbb{F} & d^\theta \end{pmatrix}$$

*is a $G_2$-instanton with respect to $\varphi =: *\psi$, i.e.*

$$F_D \wedge \psi = 0.$$

*In particular, given a solution $(\varphi, \theta)$ of the heterotic $G_2$-system (3.12), the triple $(\varphi, H_\varphi, \theta)$ solves the coupled $G_2$-instanton equation (4.2).*

*Proof.* The proof follows by direct application of Theorem 2.32, combined with Proposition 4.5. The last part of the statement follows from Proposition 3.6. □

As a straightforward consequence of Theorem 4.6 and Proposition 2.33 we obtain the following:

**Corollary 4.7.** *Let $(M^7, \varphi)$ be a 7-manifold with an integrable $G_2$-structure $\varphi$ and closed torsion, endowed with a $G_2$-instanton connection $\theta$ on a principal $K$-bundle $P \to M$ with respect to $\varphi =: *\psi$, that is, solving the equations,*

$$\tau_2 = 0, \qquad dH_\varphi = 0, \qquad F_\theta \wedge \psi = 0,$$

*cf. Proposition 3.3. Then there exists a sequence of $G_2$-instanton bundles $\{(V_k, \nabla^k)\}_{k \in \mathbb{N}}$ over $M$ with respect to $\varphi$, such that each $V_k$ is a real orthogonal bundle of rank*

$$r_k = 7 + r_{k-1}(r_{k-1} - 1), \qquad r_1 = \dim K,$$

*and $\nabla^k$ is a linear orthogonal connection on $V_k$.*

## 4.2 Gravitino solutions and generalized Ricci-flat metrics

We will provide a partial answer to Problem 2 on any oriented spin manifold $M^7$. The approach consists in considering solutions $(\varphi, \theta)$ of the gravitino equation (4.3) and the heterotic Bianchi identity (4.5), and proving that they induce generalized Ricci-flat metrics for a canonical choice of divergence determined by the Lee form of the $G_2$-structure $\varphi$. In particular, this implies that any solution of the coupled $G_2$-instanton equation (4.2) constructed via Theorem 4.6 induces a generalized Ricci-flat metric, as stated in Problem 2. Incidentally, we reach the strong conclusion that our hypothesis imply solving the full heterotic $G_2$ system (3.12), as known to be the case for coupled $\mathrm{Spin}(7)$-instantons, see Remark 4.10.

We start with a technical Lemma about the failure of a $G_2$-instanton to satisfy the Yang-Mills equations, which is valid for arbitrary $G_2$-structures.

**Lemma 4.8.** *Let $P$ be a principal $K$-bundle over 7-manifold $M^7$ with $G_2$-structure $\varphi$. Given a $G_2$-instanton $\theta$ on $P$, that is, a principal connection $\theta$ satisfying $F_\theta \wedge \psi = 0$, one has*

$$d^*_\theta F_\theta - *(F_\theta \wedge *H) + 4i_{\tau_1^\sharp} F_\theta = 0. \tag{4.6}$$

*Proof.* Recall that the instanton condition for $\theta$ is equivalent to the following equations

$$F_\theta \wedge \psi = 0 \iff F_\theta \lrcorner \varphi = 0 \iff F_\theta \lrcorner \psi = -F_\theta \iff F_\theta \wedge \varphi = -*F_\theta.$$

Taking covariant derivatives in the last expression, and using the usual Bianchi identity $d_\theta F_\theta = 0$, we obtain:

$$d_\theta * F_\theta = d_\theta(-F_\theta \wedge \varphi) = -F_\theta \wedge d\varphi,$$

which implies

$$d_\theta^* F_\theta + F_\theta \lrcorner d^*\psi = 0.$$

Applying (3.7), we have $d^*\psi = *d\varphi = \tau_0 \varphi - 3\tau_1 \lrcorner \psi + \tau_3$, and therefore

$$\begin{aligned}
F_\theta \lrcorner d^*\psi &= \tau_0 F_\theta \lrcorner \varphi - 3F_\theta \lrcorner (\tau_1 \lrcorner \psi) + F_\theta \lrcorner \tau_3 \\
&= \left(-\frac{1}{6} + \frac{7}{6}\right)\tau_0 F_\theta \lrcorner \varphi + (1-4)\tau_1 \lrcorner (F_\theta \lrcorner \psi) + F_\theta \lrcorner \tau_3 \\
&= -F_\theta \lrcorner H + \frac{7}{6}\tau_0 F_\theta \lrcorner \varphi - 4\tau_1 \lrcorner F_\theta \lrcorner \psi \\
&= -F_\theta \lrcorner H + 4\tau_1 \lrcorner F_\theta.
\end{aligned}$$

Recalling that, in seven dimensions, one has

$$X \lrcorner \alpha = (-1)^{7(k-1)} * (X^\flat \wedge *\alpha), \quad \text{for} \quad \alpha \in \Omega^k, X \in T,$$

and the Hodge star operator squares to the identity, the statement now follows from

$$\begin{aligned}
*(F_\theta \wedge *H) &= \frac{1}{2!}F_{ij} * (e^i \wedge e^j \wedge *H) = \frac{1}{2!}F_{ij} * (e^i \wedge *(*(e^j \wedge *H))) \\
&= \frac{(-1)^{14}}{2!}F_{ij} * (e^i \wedge *(e^j \lrcorner H)) = \frac{(-1)^7}{2!}F_{ij}(e_i \lrcorner e^j \lrcorner H) \\
&= F_\theta \lrcorner H,
\end{aligned}$$

where $F_\theta = \frac{1}{2!}\sum_{i,j} F_{ij} e^i \wedge e^j$ in a local orthonormal frame. $\square$

Our next result establishes the desired relation between solutions of the seven-dimensional gravitino equation, generalized Ricci-flat metrics, and the the heterotic $G_2$ system (3.12). Via Theorem 4.6, it can be regarded as a partial answer to Problem 2. For the proof we will use a general formula for the Ricci tensor of the characteristic connection of an integrable $G_2$-structure, from [IS23, Theorem 4.5], combined with the spinorial formula for the generalized Ricci tensor in Lemma 2.14.

**Theorem 4.9.** *Let $P \to M^7$ be a principal $K$-bundle over a connected, oriented, spin 7-manifold, and let $(\varphi, \theta)$ be a solution of the gravitino equation (4.3) and the Bianchi identity (4.5). Then the Riemannian metric $g = g_\varphi$ on $M$ determined by the $G_2$-structure satisfies, for a local orthonormal frame $\{v_j\}$ for $g$,*

$$\begin{aligned}
\mathrm{Rc} - \frac{1}{4}H^2 + \sum_j \langle i_{v_j} F_\theta, i_{v_j} F_\theta \rangle + 2L_{\tau_1^\sharp} g &= 0, \\
d^* H - 4d\tau_1 + 4i_{\tau_1^\sharp} H &= 0, \\
d_\theta^* F_\theta - *(F_\theta \wedge *H) + 4i_{\tau_1^\sharp} F_\theta &= 0,
\end{aligned} \tag{4.7}$$

*where $H = H_\varphi$. Consequently, $(\varphi, \theta)$ is a solution of the heterotic $G_2$ system (3.12) with $\lambda = \frac{7}{12}\tau_0$ and*

$$d\tau_1 \wedge *H = 0,$$
$$\mathcal{S}^+ = R_g - \frac{1}{2}|H|^2 + |F_\theta|^2 - 8d^*\tau_1 - 16|\tau_1|^2 = \frac{49}{36}\tau_0^2.$$

*Proof.* The third equation in (4.7) follows from Lemma 4.8 and the hypothesis (4.3). For the first two equations in (4.7), we use the integrability of $\varphi$, which implies from [IS23, Thrm 4.5] that

$$(\mathrm{Ric}_{\nabla^+})_{ij} - \frac{1}{12}(dH)_{ab\mu i}\psi_{ab\mu j} + 4\nabla_i^+ \tau_{1j} = 0. \tag{4.8}$$

Choosing an orthonormal frame $\{\zeta_\alpha\}$ for the pairing $\langle \cdot, \cdot \rangle$ on the Lie algebra $\mathfrak{k}$, we express the curvature of $\theta$ by

$$F_\theta = \frac{1}{2!}F^\alpha{}_{ij} e^i \wedge e^j \otimes \zeta_\alpha,$$

where the $\{e_j\}$ form a local orthonormal frame for the tangent bundle and here (and throughout the proof) we use summation convention. Using this, we now have

$$\langle F_\theta \wedge F_\theta \rangle = \left\langle \frac{1}{2!} F^\alpha{}_{ab}\, e^{ab} \otimes \zeta_\alpha \wedge \frac{1}{2!} F^\beta{}_{kl}\, e^{kl} \otimes \zeta_\beta \right\rangle = \frac{1}{4} F^\alpha{}_{ab} F^\beta{}_{kl} e^{abkl} \langle \zeta_\alpha, \zeta_\beta \rangle = \frac{1}{4} F^\alpha{}_{ab} F_{\alpha kl} e^{abkl}.$$

By the heterotic Bianchi identity (4.5), we have $dH = \langle F_\theta \wedge F_\theta \rangle = \frac{1}{4} F^\alpha{}_{ab} F_{\alpha kl} e^{abkl}$, which reads in local components:

$$(dH)_i = F^\alpha{}_{i\mu} F_{\alpha\nu\rho} e^{\mu\nu\rho}.$$

Contracting this expression with $\psi_j = \frac{1}{3!} \psi_{j\mu\nu\rho} e^{\mu\nu\rho}$, and using the instanton condition for $\theta$, we conclude:

$$(dH)_i \lrcorner \psi_j = F^\alpha{}_i{}^\mu F_\alpha{}^{\nu\rho} \psi_{j\mu\nu\rho} = -2 F^\alpha{}_i{}^\mu F_{\alpha j\mu}.$$

On the other hand,

$$(dH)_i \lrcorner \psi_j = \frac{1}{6}(dH)_{ab\mu i} \psi_{ab\mu j} = -2 F^\alpha{}_{i\mu} F_{\alpha j\mu} \Rightarrow (dH)_{ab\mu i} \psi_{ab\mu j} = -12 F^\alpha{}_{\mu i} F_{\alpha\mu j},$$

and hence

$$(\mathrm{Ric}_{\nabla^+})_{ij} = -F^\alpha{}_{\mu i} F_{\alpha\mu j} - 4\nabla_i^+ \tau_{1j} = -(\langle i_{e_k} F_\theta, i_{e_k} F_\theta \rangle)_{ij} - 4\nabla_i^+ \tau_{1j}. \tag{4.9}$$

The first and second equations in (4.7) now follow from the unique decomposition of $\mathrm{Rc}_{\nabla^+}$ and $\nabla^+\tau_1$ into symmetric and skew-symmetric 2-tensors given by (2.24), since $\mathrm{Rc}$ and $\langle i_{v_j} F_\theta, i_{v_j} F_\theta \rangle$ are symmetric tensors.

To prove the last part of the statement, we use the explicit formula for the generalized Ricci tensor (2.16), see also (4.7), which implies that $(g_\varphi, H_\varphi, \theta, 4\tau_1)$ determine a solution of the generalized Ricci flat equation

$$\mathrm{Rc}^+_{\mathbf{G}_\varphi, \mathrm{div}^\varphi} = 0 \tag{4.10}$$

on the string algebroid $E_{P, H_\varphi, \theta}$ from Example 2.2. Here, $\mathbf{G}_\varphi$ is obtained as in Remark 3.7 and the divergence operator is uniquely determined by the $G_2$-structure via the explicit formula given by Remark 2.27:

$$\mathrm{div}^\varphi = \mathrm{div}^{\mathbf{G}_\varphi} - \langle 8\tau_1, \cdot \rangle.$$

Consider $\eta_\varphi$, the non-vanishing spinor on $(T, g)$ determined by the $G_2$-structure $\varphi$, cf. (3.14). Applying now Lemma 2.14 and formula (3.16), arguing as in the proof of Proposition 3.6, we have

$$0 = \langle a_-, \mathrm{Rc}^+_{\mathbf{G}, \mathrm{div}} \rangle \cdot \eta = -4 D^S_{a_-} \slashed{D}^+ \eta = -\frac{7}{3} D^S_{a_-} (\tau_0 \eta) = -\frac{7}{3} (\pi(a_-) \tau_0) \eta,$$

for every $a_- \in \Omega^0(V_-)$. Consequently, $d\tau_0 = 0$ and hence $\varphi$ is a solution of the heterotic $G_2$-system (3.12) with $\lambda = \frac{7}{12} \tau_0$. The statement follows now from Theorem 3.9. $\square$

*Remark* 4.10. The proof of the previous theorem boils down to the fact that any solution of the gravitino equation and the Bianchi identity in seven dimensions is also a solution of the heterotic $G_2$-system (3.12) for a suitable choice of $\lambda \in \mathbb{R}$. A similar phenomenon occurs for $\mathrm{Spin}(7)$-structures $\Omega \in \Omega^4$ in dimension 8, where the gravitino equation and the Bianchi identity read

$$\nabla^+ \Omega = 0, \quad *(F_\theta \wedge \Omega) = -F_\theta, \quad dH_\Omega - \langle F_\theta \wedge F_\theta \rangle = 0,$$

cf. (2.34). In this case, we have the torsion forms defined by $d\Omega = \tau_1 \wedge \Omega + *\tau_3$ and $H = H_\Omega = -\frac{1}{6}\tau_1 \lrcorner \Omega - \tau_3$ is the unique torsion three-form determined by $\Omega$ [Iva04] (see also [IP23, MM18]), cf. Proposition 3.3. A pair $(\Omega, \theta)$ solving the previous equations determines uniquely a solution of the dilatino equation (2.22) with $\lambda = 0$ and $\zeta = \frac{7}{6}\tau_1$, by direct application of [Iva04, Theorem 1.1]. In particular, Proposition 2.15, Proposition 2.23, and Theorem 2.32 apply in this situation, giving more examples of generalized Ricci flat metrics and coupled $\mathrm{Spin}(7)$-instantons. $\bigcirc$

To finish this section, in the following result we investigate the failure of generalized Ricci-flatness when we remove the instanton condition on $\theta$ from the hypotheses of Theorem 4.9. We focus on the Yang–Mills-type equation given by the third equation in (4.7), which we relate to the second equation in the coupled $G_2$-instanton equations (4.2). This situation can be then compared to the case of $\mathrm{SU}(n)$-structures $(\omega, \Psi)$ with integrable complex structure studied in [GFGM23, Proposition 4.9], see Remark 4.12. A similar analysis can be adopted for the first and second equations in (4.7), following carefully the proof of Theorem 4.9. This technical result will be key for the proof of the main results in Section 5.2.

**Lemma 4.11.** *Let $(M^7, \varphi)$ be a 7-manifold endowed with an integrable $G_2$-structure $\varphi$. Let $P$ be a principal $K$-bundle over $M$ and $\theta$ an arbitrary principal connection on $P$. Then the following identity holds:*

$$d_\theta^* F_\theta + 4\tau_1 \lrcorner F_\theta - *(F_\theta \wedge *H) = 6\tau_1 \lrcorner \pi_7 F_\theta + \frac{1}{3}\tau_0 \pi_7 F_\theta \lrcorner \varphi - 3\pi_7 F_\theta \lrcorner \tau_3 - 3\sum_j i_{e_j} \pi_7 \nabla_{e_j}^{\theta,+} F_\theta \quad (4.11)$$

*for a local orthonormal frame $\{e_j\}$ on $M$. In particular, if $(\nabla^{\theta,+} F_\theta) \lrcorner \varphi = 0$, we have*

$$d_\theta^* F_\theta + 4\tau_1 \lrcorner F_\theta - *(F_\theta \wedge *H) = 6\tau_1 \lrcorner \pi_7 F_\theta + \frac{1}{3}\tau_0 \pi_7 F_\theta \lrcorner \varphi - 3\pi_7 F_\theta \lrcorner \tau_3. \quad (4.12)$$

*Proof.* Writing the expression for $\nabla_X^{\theta,+} F_\theta$ explicitly, we have

$$\begin{aligned}
\nabla_X^{\theta,+} F_\theta(V,W) &= d_X^\theta(F_\theta(V,W)) - F_\theta(\nabla_X^+ V, W) - F_\theta(V, \nabla_X^+ W) \\
&= d_X^\theta(F_\theta(V,W)) - F_\theta(\nabla_X^g V, W) - F_\theta(V, \nabla_X^g W) - \frac{1}{2}\Big(F_\theta(H(X,V), W) + F_\theta(V, H(X,W))\Big) \\
&= \nabla_X^{\theta,g} F_\theta(V,W) - \frac{1}{2}\Big(F_\theta(H(X,V), W) - F_\theta(H(X,W), V).\Big)
\end{aligned}$$

Define $\mathcal{K} \in \Omega^1(\Lambda^2 T^* \otimes \mathrm{ad}P)$ by $\nabla_X^{\theta,+} F_\theta =: \nabla_X^{\theta,g} F_\theta - \frac{1}{2}\mathcal{K}_X$. Now, using

$$\pi_7 \nabla_X^{\theta,+} F_\theta = \frac{1}{3}\Big(\nabla_X^{\theta,+} F_\theta + \nabla_X^{\theta,+} F_\theta \lrcorner \psi\Big),$$

in a local orthonormal frame $\{e_j\}$, we obtain

$$\begin{aligned}
3\sum_j i_{e_j} \pi_7 \nabla_{e_j}^{\theta,+} F_\theta &= \sum_j i_{e_j} \nabla_{e_j}^{\theta,+} F_\theta + i_{e_j}(\nabla_{e_j}^{\theta,+} F_\theta \lrcorner \psi) \\
&= \sum_j i_{e_j} \nabla_{e_j}^{\theta,+} F_\theta + \nabla_{e_j}^{\theta,+} F_\theta \lrcorner i_{e_j} \psi \\
&= \sum_j i_{e_j} \nabla_{e_j}^{\theta,+} F_\theta - \nabla_{e_j}^{\theta,+} F_\theta \lrcorner *(e_j \wedge \varphi) \\
&= \sum_j i_{e_j} \nabla_{e_j}^{\theta,+} F_\theta - *(\nabla_{e_j}^{\theta,+} F_\theta \wedge e^j \wedge \varphi).
\end{aligned}$$

We compute the first summand in the last expression:

$$\begin{aligned}
\sum_j i_{e_j} \nabla_{e_j}^{\theta,+} F_\theta &= \sum_j i_{e_j} \nabla_{e_j}^{\theta,g} F_\theta - \frac{1}{2}\Big(F_\theta(H(e_j, e_j), \cdot) + F_\theta(e_j, H(e_j, \cdot))\Big) \\
&= \sum_j e_j \lrcorner \nabla_{e_j}^{\theta,g} F_\theta - \frac{1}{2}H_{jkl} F_{jl}\, e^k \\
&= -d_\theta^* F_\theta + F_\theta \lrcorner H \\
&= -d_\theta^* F_\theta + *(F_\theta \wedge *H).
\end{aligned}$$

To compute the second summand, we can write the covariant exterior derivative as $d_\theta = \sum_j e^j \wedge \nabla_{e_j}^{\theta,g}$, so that the Bianchi identity $d_\theta F_\theta = 0$ gives

$$\sum_j e^j \wedge \nabla_{e_j}^{\theta,+} F_\theta = \sum_j e^j \wedge \nabla_{e_j}^{\theta,g} F_\theta - \frac{1}{2} e^j \wedge \mathcal{K}_{e_j} = -\frac{1}{2}\sum_j e^j \wedge \mathcal{K}_{e_j}.$$

Therefore

$$-d_\theta^* F_\theta + F_\theta \lrcorner H + \frac{1}{2}\sum_j *(e^j \wedge \mathcal{K}_{e_j} \wedge \varphi) = 3\sum_j i_{e_j} \pi_7 \nabla_{e_j}^{\theta,+} F_\theta.$$

Setting $\mathcal{K}_j = \mathcal{K}_{e_j}$, and computing directly (now using summation convention for efficiency)

$$\begin{aligned}
e^j \wedge \mathcal{K}_j &= \frac{1}{2}(\mathcal{K}_j)_{\alpha\beta}\, e^{j\alpha\beta} = \frac{1}{2}\big(H_{j\alpha}{}^\gamma F_{\gamma\beta} - H_{j\beta}{}^\gamma F_{\gamma\alpha}\big) e^{j\alpha\beta} \\
&= \frac{1}{2}\big((H_{j\alpha}{}^\gamma\, e^{j\alpha}) \wedge (F_{\gamma\beta}\, e^\beta) + (H_{j\beta}{}^\gamma\, e^{j\beta}) \wedge (F_{\gamma\alpha}\, e^\alpha)\big) \\
&= 2i_{e_\gamma}(H^\gamma \wedge F_\theta),
\end{aligned}$$

we obtain
$$*(e^j \wedge \mathcal{K}_j \wedge \varphi) = 2*(i_{e_j}(H^j \wedge F_\theta) \wedge \varphi) = 2*(i_{e_j}(H^j \wedge F_\theta \wedge \varphi) - H^j \wedge F_\theta \wedge i_{e_j}\varphi).$$

Using now that $\tau_2 = 0$, we have (cf. [dlOLS18a, Section 2.3])
$$d\varphi = \sum_j (e_j \lrcorner H) \wedge (e_j \lrcorner \varphi) = \sum_j H^j \wedge \varphi_j,$$

we deduce
$$\sum_j *(H^j \wedge F_\theta \wedge \varphi_j) = *(F_\theta \wedge d\varphi) = *(F_\theta \wedge **d**\varphi)$$
$$= F_\theta \lrcorner d^*\psi,$$
$$\Rightarrow \quad \sum_j *(i_{e_j}(H^j \wedge F_\theta \wedge \varphi)) = \sum_j *(i_{e_j}(H^j \wedge *(F_\theta \lrcorner \psi))) = \sum_j e^j \wedge *(H^j \wedge *(F_\theta \lrcorner \psi))$$
$$= \sum_j e^j \wedge (H^j \lrcorner (F_\theta \lrcorner \psi)) = \frac{1}{2} \sum_{j,k,l} H^{jkl}(F_\theta \lrcorner \psi)_{kl} e^j$$
$$= \frac{1}{2} \sum_{j,k,l} (F_\theta \lrcorner \psi)_{kl} H_{jkl} e^j$$
$$= (F_\theta \lrcorner \psi) \lrcorner H.$$

From this, we conclude that
$$\sum_j *(e^j \wedge \mathcal{K}_j \wedge \varphi) = 2(F_\theta \lrcorner d^*\psi - (F_\theta \lrcorner \psi) \lrcorner H)$$

and, as desired,
$$d_\theta^* F_\theta + 4\tau_1 \lrcorner F_\theta - *(F_\theta \wedge *H) = 4\tau_1 \lrcorner F_\theta + (F_\theta \lrcorner \psi) \lrcorner H - F_\theta \lrcorner d^*\psi - 3\sum_j i_{e_j} \pi_7 \nabla^{\theta,+}_{e_j} F_\theta$$
$$= 4\tau_1 \lrcorner F_\theta + \frac{1}{3}\tau_0 \pi_7 F_\theta \lrcorner \varphi - \tau_1 \lrcorner ((F_\theta \lrcorner \psi) \lrcorner \psi)$$
$$- (F_\theta \lrcorner \psi) \lrcorner \tau_3 - \tau_0 F_\theta \lrcorner \varphi + 3\tau_1 \lrcorner (F_\theta \lrcorner \psi) - F_\theta \lrcorner \tau_3 - 3\sum_j i_{e_j} \pi_7 \nabla^{\theta,+}_{e_j} F_\theta$$
$$= 6\tau_1 \lrcorner \pi_7 F_\theta + \frac{1}{3}\tau_0 \pi_7 F_\theta \lrcorner \varphi - 3\pi_7 F_\theta \lrcorner \tau_3 - 3\sum_j i_{e_j} \pi_7 \nabla^{\theta,+}_{e_j} F_\theta.$$

The last part of the statement follows by noticing that $\nabla^{\theta,+} F_\theta \lrcorner \varphi = 0 \Leftrightarrow \pi_7 \nabla^{\theta,+} F_\theta = 0$, cf. Proposition 3.1. □

*Remark* 4.12. Equation (4.12) above shows us that the second equation in the coupled $G_2$-instanton equations (4.2), given by
$$\nabla^{\theta,+} F_\theta \lrcorner \varphi = 0,$$
does not imply, in general, the Yang–Mills equation with torsion, given by the third equation in (4.7). An explicit example where this is indeed the case is not known to the authors. This situation stands out in comparison to $SU(n)$-structures $(\omega, \Psi)$ with integrable complex structure studied in [GFGM23, Proposition 4.9], for which the equation $(\nabla^{\theta,+} F_\theta) \wedge \omega^{n-1} = 0$ combined with $F_\theta^{0,2} = 0$ is equivalent to the corresponding Yang–Mills equation with torsion.
○

## 4.3 Examples

In this section we discuss some examples of coupled $G_2$-instantons which can be found scattered in the literature, but which apparently have not been identified as such. The first examples arise from solutions of the heterotic $G_2$ system (3.12), by Theorem 4.6 and Proposition 3.6. Such solutions with exact torsion one-form $\tau_1 = d\phi$ have been constructed in e.g. [GN95, FIUV11, Nol12, FIUV15, dlOG21, CGFT22, GS24], motivated by the concept's origins in heterotic string theory, which require a globally defined dilaton field $\phi$ whose vacuum expectation value determines the string-coupling constant. The approximate solutions constructed in [LSE23] deserve a special treatment, since they do not exactly solve the first equation in (3.12), and we postpone their analysis to Section 5.

Since we are mainly concerned with solving the coupled instanton equation (4.2), we will work with Theorem 4.6 and consider solutions of the gravitino equation (4.3) and the Bianchi identity (4.5). Incidentally, by Theorem 4.9, these conditions are sufficient to imply a solution of the heterotic $G_2$ system, with our relaxed definition (3.12). Note that our equations barely impose any constraint on the torsion one-form $\tau_1$, and therefore are more flexible than the ones usually considered in the mathematical physics literature, yet are still strong enough to prove Theorem 3.8 and Theorem 4.6.

Our first two examples are given by the product of a flat torus with a manifold carrying an $SU(n)$-structure which is integrable and has closed torsion, also known in the literature as *twisted Calabi-Yau* [GFRT20b] or, more generally, *Bismut Hermitian–Einstein metrics* [GFJS23]. The seven-dimensional geometry is given by a strong integrable $G_2$-structure (see Proposition 3.3), i.e. such that

$$\tau_2 = 0 \quad \text{and} \quad dH_\varphi = 0.$$

**Example 4.13.** Let $N^4$ be a four-dimensional manifold endowed with $SU(2)$-structure $(\omega, \Psi)$, with almost complex structure $J$ and Hermitian metric $g = \omega(\cdot, J\cdot)$. Its Lee form $\theta_\omega := -J^*d^*\omega \in \Omega^1(N)$ is defined by

$$d\omega = \theta_\omega \wedge \omega.$$

In this setup, a solution of the gravitino equation is a triple $(\omega, \Psi, H)$ such that, cf. [FI03, Theorem 10.1],

$$H = -d^c\omega + g(N_J, \cdot),$$

where $N_J$ is the Nijenhuis tensor of $J$, which in particular must be skew-symmetric.

Suppose that $(\omega, \Psi)$ satisfies the *twisted Calabi-Yau equation*, introduced in [GFRST22]:

$$d\Psi = \theta_\omega \wedge \Psi, \quad d\theta_\omega = 0, \quad dd^c\omega = 0. \tag{4.13}$$

Then, it was proved in [GFRST22, Lemma 2.2] that $N_J = 0$ and that it determines a solution of the gravitino equation with $H = -d^c\omega$, which also solves the Bianchi identity $dH = 0$. Note that compact solutions of these equations in four dimensions are rather rigid, as they only exist on tori and K3 surfaces, with $H = 0 = \theta_\omega$, and diagonal Hopf surfaces, with $H \neq 0 \neq \theta_\omega$, cf. [GFRST22, Proposition 2.10].

To build the seven-dimensional geometry from a solution of (4.13) we follow closely [FMMR23]. We consider $M = N \times T^3$, where $T^3$ is a three-dimensional flat torus. Denote

$$\psi_+ := \text{Re}(\Psi), \qquad \psi_- := \text{Im}(\Psi).$$

Define a $G_2$-structure on $M$ by

$$\varphi = dx^1 \wedge dx^2 \wedge dx^3 + dx^1 \wedge \omega + dx^2 \wedge \psi_+ - dx^3 \wedge \psi_-, \tag{4.14}$$

where $(x_1, x_2, x_3) \in \mathbb{R}$ are coordinates in the universal cover of $T^3$. Then, $\varphi$ is strong and integrable with

$$\tau_0 = 0, \qquad \theta_\omega = 4\tau_1, \qquad H_\varphi = d^c\omega.$$

For the proof we follow [FMMR23, Proposition 3.5]. For instance, since $\omega$, $\psi_+$ and $\psi_-$ are Hodge self-dual on $N^4$,

$$*\varphi = \tfrac{1}{2}\omega^2 + dx^2 \wedge dx^3 \wedge \omega - dx^1 \wedge dx^3 \wedge \psi_+ - dx^1 \wedge dx^2 \wedge \psi_-.$$

Since $d\omega = \theta_\omega \wedge \omega$ and $d\Psi = \theta_\omega \wedge \Psi$, we obtain $d*\varphi = \theta_\omega \wedge *\varphi$. Thus, $\varphi$ is an integrable $G_2$-structure with Lee form $\theta = \theta_\omega$, and

$$d\varphi \wedge \varphi = \theta_\omega \wedge \left(dx^1 \wedge \omega + dx^2 \wedge \psi_+ - dx^3 \wedge \psi_-\right) \wedge \varphi = 0,$$

since the self-dual forms $\omega$, $\psi_+$ and $\psi_-$ are pairwise orthogonal. The torsion of $\varphi$ is

$$H_\varphi = *(\theta \wedge \varphi - d\varphi) = *(\theta_\omega \wedge dx^1 \wedge dx^2 \wedge dx^3) = -*_4 \theta_\omega = Jd\omega = d^c\omega,$$

since $\theta_\omega = J *_4 d *_4 \omega = *_4 Jd\omega$ and $J\omega = \omega$.

Applying now Theorem 2.32, we obtain a coupled $G_2$-instanton on $TM$ given by the connection

$$\nabla^- = \nabla^{g_7} - \tfrac{1}{2}g_7^{-1}H_\varphi. \tag{4.15}$$

Incidentally, this connection is actually flat [FMMR23], and the tower of coupled $G_2$-instantons over this manifold given by Corollary 4.7 is also flat. $\triangle$

The following example is given by the product of the Calabi–Eckmann 6-manifold $S^3 \times S^3$ with a circle. Similarly as in the previous example, by application of [FMMR23, Proposition 3.5] the seven dimensional torsion classes are inherited from the six dimensional geometry. Consequently, this example is also strong and integrable, but unlike the previous one has $d\tau_1 \neq 0$, which reflects the fact that the Calabi–Eckmann complex threefold does not admit balanced hermitian metrics. Note that this example provides a solution of the heterotic $G_2$ system (3.12) according to our lax Definition 3.4, but it escapes from the orthodoxy for these systems of equations in the literature, precisely due to the fact that $\tau_1$ is non-closed.

**Example 4.14.** Let
$$N^6 = \{\mathbb{C}^2_\times \times \mathbb{C}^2_\times\}/\mathbb{C} \simeq S^3 \times S^3,$$
with its (non-Kähler) Calabi–Eckmann SU(3)-structure $(\omega, \Psi)$. Following [GFS20, Example 8.35], if we let $\pi_j : S^3 \to \mathbb{CP}^1$ denote the Hopf fibration on each of the two factors in $N$, for $j = 1, 2$, and let $\mu_j$ denote the 1-form on $S^3$ such that
$$d\mu_j = \pi_j^* \omega_{\mathbb{CP}^1}$$
for $j = 1, 2$, where $\omega_{\mathbb{CP}^1}$ is the Kähler form for the Fubini–Study metric on $\mathbb{CP}^1$, then we can write $\omega$ explicitly as:
$$\omega = \pi_1^* \omega_{\mathbb{CP}^1} + \pi_2^* \omega_{\mathbb{CP}^1} + \mu_1 \wedge \mu_2.$$
It is straightforward to show that if we let
$$\theta_\omega = \mu_2 - \mu_1$$
then
$$d\Psi = \theta_\omega \wedge \Psi \quad \text{and} \quad dd^c \omega = 0.$$
However, note that
$$d\theta_\omega = \pi_2^* \omega_{\mathbb{CP}^1} - \pi_1^* \omega_{\mathbb{CP}^1} \neq 0$$
and so the second equation (4.13) in the definition of twisted Calabi–Yau is not satisfied, though the rest are.

As in Example 4.13, we now let $\psi_\pm$ denote the real and imaginary parts of $\Psi$ respectively. We may then define a product $G_2$-structure on $M^7 = N^6 \times S^1$ by
$$\varphi = \omega \wedge dt + \psi_+,$$
where $dt$ is the standard nowhere vanishing 1-form on $S^1$. The Hodge dual $\psi$ of $\varphi$ is then given by
$$\psi = \frac{1}{2}\omega^2 + \psi_- \wedge dt.$$
As in [FMMR23, Proposition 3.5], one sees from these formulae that the $G_2$-structure $\varphi$ is integrable with
$$\tau_0 = 0, \quad \tau_1 = \theta_\omega, \quad H_\varphi = d^c \omega = \pi_1^* \omega_{\mathbb{CP}^1} \wedge \mu_1 - \pi_2^* \omega_{\mathbb{CP}^1} \wedge \mu_2.$$
Hence, $dH_\varphi = dd^c \omega = 0$ and thus $\varphi$ is also strong.

Even though $\tau_1$ is not closed, one may still apply Theorem 2.32 and obtain a coupled $G_2$-instanton $\nabla^-$ on $TM$ as in (4.15), which is again flat. The tower of coupled $G_2$-instantons we obtain from Corollary 4.7 are also flat. Notice that the $G_2$-structure presented here is fundamentally different from that obtained on $G = \text{SU}(2)^2 \times S^1$ as a Lie group, in [FMMR23, Proposition 6.2]. △

*Remark* 4.15. Given the observations in Examples 4.13 and 4.14, it would be interesting to find $G_2$-structures which are both strong and integrable but for which the connection $\nabla^-$ in (4.15) is *not flat*, or even irreducible. ○

*Remark* 4.16. No irreducible compact homogeneous spaces admitting invariant $G_2$-structures, up to a covering, admit (invariant) strong integrable $G_2$-structures, cf. [FMMR23, §5]. On the other hand, the same authors find numerous examples of such structures on *reducible* spaces, which according to their preference have closed Lee form, cf. [FMMR23, §6]. Several of those examples can be easily adapted to provide more general solutions of the gravitino equation (4.13). ○

Our final example, originally found in [II05, §6], provides a solution of the heterotic $G_2$-system (3.12) with $\tau_0 \neq 0$ and non-flat instanton $\theta$ in the nearly parallel seven-dimensional sphere. In particular, the coupled $G_2$-instanton obtained from this solution via Theorem 4.6 is non-flat.

**Example 4.17.** Let $M = S^7$ be the standard 7-sphere, viewed as a sphere in the octonions. It is well-known that the embedding of $S^7$ in the octonions induces a natural $\mathrm{Spin}(7)$-invariant $\mathrm{G}_2$-structure $\varphi$ on $S^7$ which is *nearly parallel* in the sense that
$$d\varphi = 4\kappa\psi \tag{4.16}$$
where $\psi = *\varphi$ as usual and $\kappa \neq 0$ is constant. Note that the metric determined by $\varphi$ has constant curvature $\kappa^2$. We clearly see that all the torsion forms vanish, except $\tau_0 = 4\kappa \neq 0$ and is constant. Hence, $\varphi$ is an integrable $\mathrm{G}_2$-structure but since
$$H = \frac{2}{3}\kappa\varphi, \tag{4.17}$$
we see that
$$dH = \frac{8}{3}\kappa^2\psi \neq 0 \tag{4.18}$$
by (4.16) and thus $\varphi$ is not strong.

Take $P$ to be the $\mathrm{G}_2$-frame bundle of $S^7$, and let $\theta$ be the connection on $P$ determined by $\nabla^+$. It is observed in [II05, §6] that $\nabla^+ H_\varphi = 0$, and consequently, cf. (2.41),
$$g(R_{\nabla^+}(X,Y)Z, W) = g(R_{\nabla^+}(Z,W)X, Y).$$
Arguing as in the proof of Theorem 2.32, it follows that $\theta$ is a $\mathrm{G}_2$-instanton. Furthermore, it is shown in [II05, §6] that the curvature $F_\theta$ of $\theta$ satisfies
$$\mathrm{tr}\, F_\theta \wedge F_\theta = -\frac{32\kappa^4}{27}\psi. \tag{4.19}$$
Combining (4.18) and (4.19), we see that the heterotic Bianchi identity (1.1) is satisfied, for a suitable choice of scaling of the Killing form on the Lie algebra of $\mathrm{G}_2$. Overall, we see that $(\varphi, H, \theta)$ defines a coupled $\mathrm{G}_2$-instanton on $S^7$. △

*Remark* 4.18. An interesting example in six dimensions where Theorem 2.32 applies is the 6-sphere with the standard nearly Kähler structure inherited from the imaginary octonions. According to [II05, §6], this provides a solution of the gravitino equation and the Bianchi identity with instanton connection $\nabla^+$ and non-closed torsion given by the Nijenhuis tensor of the $\mathrm{SU}(3)$-structure, with a structure very similar to Example 4.17. ◯

## 5 Approximate solutions

As we have seen, there is a connection between solutions of the heterotic $\mathrm{G}_2$ system, solutions of the coupled $\mathrm{G}_2$-instanton equations and the vanishing of generalized Ricci curvature. In [LSE23], "approximate" solutions to the heterotic $\mathrm{G}_2$ system were given in the sense that the connections involved were only "approximate" $\mathrm{G}_2$-instantons: here the "approximate" pertains to dependence on the non-zero constant $\alpha'$ which appears in the heterotic Bianchi identity as $\alpha' \to 0$. Motivated by this and our results thus far, in this section we propose a new definition of $\alpha'$-approximate $\mathrm{G}_2$-instantons and show that it not only leads to approximate solutions to the coupled $\mathrm{G}_2$-instanton equations, but also to generalized Ricci curvature which is approximately zero in a quantitative sense as $\alpha' \to 0$. We also demonstrate that the aforementioned examples from [LSE23] provide $\alpha'$-approximate $\mathrm{G}_2$-instantons and thus lead to approximate coupled $\mathrm{G}_2$-instantons and approximate generalized Ricci-flatness.

### 5.1 Motivation: contact Calabi–Yau 7-manifolds

In [LSE23], the heterotic $\mathrm{G}_2$ system was studied in the context of contact Calabi–Yau 7-manifolds, which admit a natural 1-parameter family of $\mathrm{G}_2$-structures that we now recall.

**Definition 5.1.** Let $(V, \omega, \Omega)$ be a Calabi–Yau 3-orbifold, i.e. a Kähler 3-orbifold with Kähler form $\omega$ and holomorphic volume form $\Omega$ satisfying
$$\mathrm{vol}_V = \frac{\omega^3}{3!} = \frac{1}{4}\mathrm{Re}\,\Omega \wedge \mathrm{Im}\,\Omega,$$
where $\mathrm{vol}_V$ is the volume form associated with the Kähler metric $g_V$ on $V$. Suppose that the total space of an $S^1$-(orbi)bundle $\pi : M^7 \to V$ is a *contact Calabi–Yau 7-manifold*, i.e. $M$ is endowed with a connection 1-form $\eta$ such that $d\eta = \omega$. For every $\varepsilon > 0$, we define an $S^1$-invariant $\mathrm{G}_2$-structure $\varphi_\varepsilon$ on $M^7$, with dual 4-form $\psi_\varepsilon$, by
$$\varphi_\varepsilon = \varepsilon\eta \wedge \omega + \mathrm{Re}\,\Omega \quad \text{and} \quad \psi_\varepsilon = \frac{1}{2}\omega^2 - \varepsilon\eta \wedge \mathrm{Im}\,\Omega. \tag{5.1}$$
The metric induced from this $\mathrm{G}_2$-structure and its corresponding volume form on $M$ are
$$g_\varepsilon = \varepsilon^2 \eta \otimes \eta + g_V, \qquad \mathrm{vol}_\varepsilon = \varepsilon\eta \wedge \mathrm{vol}_V. \tag{5.2}$$

Note that varying $\varepsilon$ in (5.2) amounts to rescaling the $\mathbb{S}^1$ fibres of $\pi: M \to V$, so that $\varepsilon \to 0$ corresponds to collapsing the fibres to zero size.

We now recall some basic observations about the family of $G_2$-structures $\varphi_\varepsilon$ in (5.1) on the contact Calabi–Yau 7-manifold $M$ from [LSE23, Lemmas 2.4 & 2.5]. We see that

$$d\varphi_\varepsilon = \varepsilon\omega^2, \qquad d\psi_\varepsilon = 0,$$

so the $G_2$-structures are co-closed. The torsion forms of $\varphi_\varepsilon$ are:

$$\tau_0 = \frac{6}{7}\varepsilon; \quad \tau_1 = 0; \quad \tau_2 = 0; \quad \tau_3 = \frac{8}{7}\varepsilon^2 \eta \wedge \omega - \frac{6}{7}\varepsilon \mathrm{Re}\Omega.$$

In particular, we observe that the structures are integrable (i.e. $\tau_2 = 0$) and admit a connection with totally skew-symmetric torsion

$$H_\varepsilon = -\varepsilon^2 \eta \wedge \omega + \varepsilon \mathrm{Re}\Omega, \tag{5.3}$$

which satisfies

$$dH_\varepsilon = -\varepsilon^2 \omega^2.$$

The above facts are important for showing that one can build approximate solutions to the heterotic $G_2$-system on $M$, using the $G_2$-structures $\varphi_\varepsilon$. To describe these approximate solutions, it is necessary to introduce a useful (and natural) local coframe which is adapted to the geometry of $M$.

**Definition 5.2.** *Given $\varepsilon > 0$, let $(M^7, \varphi_\varepsilon)$ be as in Definition 5.1. We choose a local Sasakian real orthonormal coframe on $M$:*

$$e^0 = \varepsilon\eta, \quad e^1, \quad e^2, \quad e^3, \quad Je^1, \quad Je^2, \quad Je^3, \tag{5.4}$$

*where $J$ is the transverse complex structure (from the Calabi–Yau $V$) so that $\{e^1, e^2, e^3, Je^1, Je^2, Je^3\}$ is a basic $SU(3)$-frame for $V$. In this frame, the Kähler and holomorphic volume forms are given by:*

$$\omega = e^1 \wedge Je^1 + e^2 \wedge Je^2 + e^3 \wedge Je^3 \quad \text{and} \quad \Omega = (e^1 + iJe^1) \wedge (e^2 + iJe^2) \wedge (e^3 + iJe^3).$$

By [LSE23, Proposition 3.2], we know that if we write $e = (e_1 \; e_2 \; e_3)^\mathrm{T}$ then the following structure equations hold:

$$d \begin{pmatrix} e_0 \\ e \\ Je \end{pmatrix} = - \begin{pmatrix} 0 & \frac{\varepsilon}{2}Je^\mathrm{T} & -\frac{\varepsilon}{2}e^\mathrm{T} \\ -\frac{\varepsilon}{2}Je & a & b - \frac{\varepsilon}{2}e_0 I \\ \frac{\varepsilon}{2}e & -b + \frac{\varepsilon}{2}e_0 I & a \end{pmatrix} \wedge \begin{pmatrix} e_0 \\ e \\ Je \end{pmatrix}, \tag{5.5}$$

where $a$ is a skew-symmetric $3 \times 3$ matrix of 1-forms, $b$ is a symmetric traceless $3 \times 3$ matrix of 1-forms, and $I$ is the $3 \times 3$ identity matrix. Therefore, if we define

$$\mathbf{A} = \begin{pmatrix} 0 & 0 & 0 \\ 0 & a & b \\ 0 & -b & a \end{pmatrix} \quad \text{and} \quad \mathbf{B} = \begin{pmatrix} 0 & Je^\mathrm{T} & -e^\mathrm{T} \\ -Je & 0 & -e_0 I \\ e & e_0 I & 0 \end{pmatrix} \tag{5.6}$$

we see that $\mathbf{A} + \frac{\varepsilon}{2}\mathbf{B}$ is the local matrix representing the Levi-Civita connection of $g_\varepsilon$ with respect to the local orthonormal coframe introduced in Definition 5.2. If we then let

$$\mathbf{I} = \begin{pmatrix} 0 & 0 & 0 \\ 0 & 0 & -I \\ 0 & I & 0 \end{pmatrix} \quad \text{and} \quad \mathbf{C} = \begin{pmatrix} 0 & Je^\mathrm{T} & -e^\mathrm{T} \\ -Je & -[e] & [Je] \\ e & [Je] & [e] \end{pmatrix} - e_0 \mathbf{I}, \tag{5.7}$$

where

$$\left[ \begin{pmatrix} e_1 \\ e_2 \\ e_3 \end{pmatrix} \right] = \begin{pmatrix} 0 & e_3 & -e_2 \\ -e_3 & 0 & e_1 \\ e_2 & -e_1 & 0 \end{pmatrix}, \tag{5.8}$$

one can define a family of connections on $TM$ as follows [LSE23, Proposition 3.21].

**Definition 5.3.** *Let $(M^7, \varphi_\varepsilon)$ be as in Definition 5.1 for some $\varepsilon > 0$. Recall the local coframe on $M$ in Definition 5.2 and the matrices $\mathbf{A}, \mathbf{B}, \mathbf{C}, \mathbf{I}$ defined with respect to this coframe in (5.6)–(5.7). For $k \in \mathbb{R} \setminus \{0\}$ and $\delta, m \in \mathbb{R}$ we define a connection $\theta^{\delta,k}_{\varepsilon,m}$ on $TM$ by the formula*

$$\theta^{\delta,k}_{\varepsilon,m} = \mathbf{A} + \frac{k\varepsilon}{2}\mathbf{B} + \frac{k\varepsilon\delta}{2}\mathbf{C} + \frac{km\varepsilon}{2}e_0 \mathbf{I}. \tag{5.9}$$

*Note that this local expression determines a globally defined connection on $TM$, that taking $\delta = m = 0$ and $k = 1$ in (5.9) yields the Levi-Civita connection $\nabla^{g_\varepsilon}$ of the metric $g_\varepsilon$ on $M$, and taking $\delta = k = 1$ and $m = 0$ in (5.9) yields the Bismut connection $\nabla^+$ associated with $g_\varepsilon$ and torsion $H_\varepsilon$.*

*Remark* 5.4. We can interpret the various parameters appearing in Definition 5.3 as follows. First, the parameter $k$ can be viewed as a "squashing" parameter, allowing us to rescale the connection along the fibres of $\pi : M \to V$ independently of the parameter $\varepsilon$. The matrix $\mathbf{C}$ is equivalent (up to a factor of $\varepsilon$) to the torsion $H_\varepsilon$ in (5.3) by [LSE23, Proposition 3.10], so the parameter $\delta$ varies the torsion of the connection along a canonical line, which contains the Bismut, Hull and Levi-Civita connections when $k = 1$ and $m = 0$. Finally, the parameter $m$ can be viewed as an additional "twist" parameter acting in the transverse directions for the fibration of $M$ over $V$. ◯

In [LSE23, Corollary 3.27], it was described how $\theta^{\delta,k}_{\varepsilon,m}$ in Definition 5.3 fails to be a $G_2$-instanton.

**Proposition 5.5.** *Using the notation of Definition 5.3, (5.7) and (5.8), the curvature $R^{\delta,k}_{\varepsilon,m}$ of the connection $\theta^{\delta,k}_{\varepsilon,m}$ on $TM$ satisfies:*

$$R^{\delta,k}_{\varepsilon,m} \wedge \psi_\varepsilon = \frac{k\varepsilon^2(6(1-\delta+m) + k(1-\delta)(1+3\delta))}{4} \cdot \frac{\omega^3}{3!}\mathbf{I} + \frac{k^2\varepsilon^2}{4}\eta \wedge \frac{\omega^2}{2!} \wedge \mathbf{M}^\delta_m, \tag{5.10}$$

*where*

$$\mathbf{M}^\delta_m = \begin{pmatrix} 0 & (1+m-5\delta)(1+\delta)e^T & (1+m-5\delta)(1+\delta)Je^T \\ (5\delta-1-m)(1+\delta)e & (\delta^2 - 2(2+m)\delta - 1)[Je] & (\delta^2 - 2(2+m)\delta - 1)[e] \\ (5\delta-1-m)(1+\delta)Je & (\delta^2 - 2(2+m)\delta - 1)[e] & -(\delta^2 - 2(2+m)\delta - 1)[Je] \end{pmatrix}.$$

*In particular, $\theta^{\delta,k}_{\varepsilon,m}$ is never a $G_2$-instanton.*

The main result on the heterotic $G_2$ system in this contact Calabi–Yau setting is the following [LSE23, cf. Theorem 1].

**Theorem 5.6.** *Let $\pi : M^7 \to V$ as in Definition 5.1 be a contact Calabi–Yau 7-manifold. Let $A$ be the pullback of the Levi-Civita connection of the Calabi–Yau metric on $V$, defined on $E = \pi^*TV$.*

*For all $\alpha' > 0$ there exist $\varepsilon = \varepsilon(\alpha') > 0$, $k = k(\alpha') > 0$, both tending to zero as $\alpha' \to 0$, and $\delta, m \in \mathbb{R}$ so that if $M$ is endowed with the $G_2$-structure $\varphi_\varepsilon$ as in (5.1), the connection $\theta^{\delta,k}_{\varepsilon,m}$ in Definition 5.3 on $TM$ and the connection $A$ on $E$, then we have a solution to the heterotic $G_2$ system, except that $\theta^{\delta,k}_{\varepsilon,m}$ is never a $G_2$-instanton but instead satisfies*

$$|R^{\delta,k}_{\varepsilon,m} \wedge \psi_\varepsilon|_{g_\varepsilon} = \mathcal{O}(\alpha')^2 \quad \text{as } \alpha' \to 0. \tag{5.11}$$

Concretely, three separate regimes are presented in [LSE23] of choices of the parameters $\varepsilon, k, \delta, m$ so that the conclusion of Theorem 5.6 holds for any positive $\alpha'$ sufficiently close to 0.

**Case 1.** $\delta \in \mathbb{R} \setminus \{0, -1\}$, $m = \delta - 1$, $k^2 = (\alpha')^{-3}$, $\varepsilon^2 = \dfrac{8}{\delta^2(1+\delta)^2}(\alpha')^5$.

**Case 2.** $\delta = 0$, $m < -1$, $k = (\alpha')^{-3}$, $\varepsilon^2 = -\dfrac{8}{(1+m)(1+3(\alpha')^3)}(\alpha')^8$.

**Case 3.** $\delta = -1$, $m > -2$, $k = (\alpha')^{-3}$, $\varepsilon^2 = \dfrac{8}{(2+m)(4-3(\alpha')^3)}(\alpha')^8$.

We shall return to the examples in Theorem 5.6 at the end of this section to understand in what sense the condition (5.11) gives "approximate" $G_2$-instantons and thus approximate solutions to the heterotic $G_2$ system.

## 5.2 Approximate $G_2$-instantons, coupled $G_2$-instantons and generalized Ricci curvature

We now return to the general setting of 7-manifolds with integrable $G_2$-structures. Given the relationship between solutions of the heterotic $G_2$ system, coupled $G_2$-instantons and the vanishing of the generalized Ricci curvature as seen in Section 4, and based on the results in Theorem 5.6, we are motivated to define a suitable notion of approximate $G_2$-instantons, and then to show that this leads to an appropriate sense of both approximate coupled $G_2$-instantons and approximate generalized Ricci-flatness.

Given this goal, we propose the following definition of approximate $G_2$-instantons in our context.

**Definition 5.7.** *Suppose that for a sequence of non-zero real numbers $\alpha' \to 0$ we have the following data.*

*Let $(M^7, \varphi)$ be a 7-manifold endowed with an integrable $G_2$-structure with induced metric $g$, dual 4-form $\psi$ and torsion 3-form $H$. Let $P \to M$ be a principal $K$-bundle over $M$, where the Lie algebra $\mathfrak{k}$ is endowed with an $\alpha'$-independent, non-degenerate, bilinear, symmetric pairing $\langle \cdot, \cdot \rangle : \mathfrak{k} \otimes \mathfrak{k} \to \mathbb{R}$. Let $\theta \in \Omega^1(P, \mathfrak{k})$ define a connection on $P$ with curvature $F_\theta$ and recall the induced connection $\nabla^{\theta,+}$. Suppose finally that the heterotic Bianchi identity is satisfied:*

$$dH = \alpha'\langle F_\theta \wedge F_\theta \rangle. \tag{5.12}$$

*We say that the connections $\theta$ are $\alpha'$-approximate $G_2$-instantons if*

$$\left|F_\theta \wedge \psi\right|_g = \mathcal{O}(\alpha')^2 \quad \text{and} \quad \left|\nabla^{\theta,+} F_\theta \wedge \psi\right|_g = \mathcal{O}(\alpha')^2 \quad \text{as } \alpha' \to 0. \tag{5.13}$$

*Remark* 5.8. Note that if $\varphi$ is integrable, then $\nabla^+$ preserves $\varphi$ and hence $\psi$, and so

$$\nabla^{\theta,+}(F_\theta \wedge \psi) = \nabla^{\theta,+} F_\theta \wedge \psi. \tag{5.14}$$

In particular, if $\theta$ is a $G_2$-instanton then $\nabla^{\theta,+} F_\theta \wedge \psi = 0$. Hence $G_2$-instantons give trivial examples of $\alpha'$-approximate $G_2$-instantons. ◯

In general, the first condition in (5.13), which is the one considered in [LSE23] (see Theorem 5.6), does not imply the second. Definition 5.7 therefore gives a stronger notion of approximate $G_2$-instanton which we shall see appears to be more natural, at least in our context.

*Remark* 5.9. As is well-known, $G_2$-instantons can bubble, which means that in a family their curvature can blow-up pointwise. To avoid this, it is natural to impose that their curvature stays bounded (pointwise), and so we can ask the same of our approximate $G_2$-instantons. It is in this setting that we can achieve our main results concerning approximate solutions. ◯

We now show that $\alpha'$-approximate $G_2$-instantons yield approximate coupled $G_2$-instantons in the following sense.

**Theorem 5.10.** *Suppose that we have $\alpha'$-approximate $G_2$-instantons on a principal $K$-bundle over $(M^7, \varphi)$ with integrable $G_2$-structure $\varphi$ and torsion $H$ satisfying (5.12) as in Definition 5.7. Recall $\rho(\varphi, H)$ given in Definition 4.1 and $\mathbb{F} \wedge \mathbb{F}^\dagger$ given in Lemma 2.28.*

*If the curvature $F_\theta$ of $\theta$ is bounded as $\alpha' \to 0$, then $(\varphi, H, \theta)$ give approximate solutions to the coupled $G_2$-instanton equation (4.2) in Definition 4.4 in the following sense as $\alpha' \to 0$:*

$$\begin{aligned}
|\rho(\varphi, H) + \langle F_\theta, (F_\theta \lrcorner \varphi)^\sharp \rangle|_g &= \mathcal{O}(\alpha')^2, \\
|(\nabla^{\theta,+} F_\theta) \lrcorner \varphi|_g &= \mathcal{O}(\alpha')^2, \\
|[F_\theta \lrcorner \varphi,] - \mathbb{F} \wedge \mathbb{F}^\dagger \lrcorner \varphi|_g &= \mathcal{O}(\alpha')^2, \\
dH - \alpha' \langle F_\theta \wedge F_\theta \rangle &= 0.
\end{aligned} \tag{5.15}$$

*Proof.* Since $\varphi$ is integrable, $\rho(\varphi, H) = 0$ by Lemma 4.3. The first equation in (5.15) is then an immediate consequence of the boundedness of $F_\theta$ and the first condition in (5.13) of $\alpha'$-approximate $G_2$-instantons. The second equation in (5.15) is precisely the second condition in (5.13). The fourth equation in (5.15) is satisfied by assumption. We are therefore only left with the third equation in (5.15).

The first term in the third equation is of order $\mathcal{O}(\alpha')^2$ by the first condition in (5.13). In (2.46) in the proof of Theorem 2.32, we saw locally we can write $\mathbb{F} \wedge \mathbb{F}^\dagger$ as:

$$- \sum_j \left( e_j \lrcorner \langle \cdot, \cdot \rangle^{-1}(\zeta_l \lrcorner F_\theta) \right) \wedge \left( e_j \lrcorner (\zeta_k \lrcorner F_\theta) \right) \otimes \zeta^k \otimes \zeta_l,$$

where $\{e_j\}$ form a local orthonormal frame on $M^7$ and $\{\zeta_j\}$ give an orthonormal basis for the Lie algebra of $K$. By Lemma 2.31 and the proof of Theorem 2.32, we deduce that $\pi_7(F_\theta) = 0$ forces $\pi_7(\mathbb{F} \wedge \mathbb{F}^\dagger) = 0$ and, moreover, there is a universal constant $C > 0$ so that

$$|\pi_7(\mathbb{F} \wedge \mathbb{F}^\dagger)|_g \leq C |F_\theta|_g |\pi_7(F_\theta)|_g.$$

The third equation in (5.15) now follows from the boundedness of $F_\theta$ and the first condition in (5.13). □

Now, as we have seen, taking $G_2$-instantons $\theta$ in Definition 5.7 leads to generalized Ricci-flatness because in this case, the following two terms, which are the components of the generalized Ricci curvature, must vanish:

$$d_\theta^* F_\theta + 4\tau_1 \lrcorner F_\theta - F_\theta \lrcorner H, \qquad \text{Ric}^+ + \alpha' \sum_j \langle i_{v_j} F_\theta, i_{v_j} F_\theta \rangle + \nabla^+ \tau_1. \tag{5.16}$$

We now examine each of these terms in turn for $\alpha'$-approximate $G_2$-instantons.

**Proposition 5.11.** *Let $\theta$ be $\alpha'$-approximate $G_2$-instantons over $(M^7, \varphi)$ as in Definition 5.7. Then the curvature $F_\theta$ satisfies*

$$\left| d_\theta^* F_\theta + 4\tau_1 \lrcorner F_\theta - F_\theta \lrcorner H \right|_g = \mathcal{O}(\alpha')^2 \quad \text{as } \alpha' \to 0.$$

*Proof.* Recall that, by Lemma 4.11, we have

$$d_\theta^* F_\theta + 4\tau_1 \lrcorner F_\theta - F_\theta \lrcorner H = 6\tau_1 \lrcorner \pi_7 F_\theta + \frac{1}{3}\tau_0 \pi_7 F_\theta \lrcorner \varphi - 3\pi_7 F_\theta \lrcorner \tau_3 + 3\sum_{j=1}^7 e_j \lrcorner \pi_7 \nabla^{\theta,+}_{e_j} F_\theta. \tag{5.17}$$

The result now immediately follows from the $\alpha'$-approximate $G_2$-instanton condition (5.13). □

We now turn to the second term in (5.16).

**Proposition 5.12.** *Let $(M^7, \varphi)$ be a 7-manifold with an integrable $G_2$-structure and $\theta$ be a connection on a principal $K$-bundle over $M$ as in Definition 5.7 so that the heterotic Bianchi identity (5.12) is satisfied. Then*

$$\mathrm{Ric}^+ + \alpha' \sum_j \langle i_{v_j} F_\theta, i_{v_j} F_\theta \rangle + 4\nabla^+ \tau_1 = 3\alpha' \sum_j \langle i_{v_j} F_\theta, i_{v_j} \pi_7 F_\theta \rangle. \quad (5.18)$$

*In particular, if $\theta$ are $\alpha'$-approximate $G_2$-instantons as in Definition 5.7 with bounded curvature $F_\theta$ as $\alpha' \to 0$, then*

$$\left| \mathrm{Ric}^+ + \alpha' \sum_j \langle i_{v_j} F_\theta, i_{v_j} F_\theta \rangle + 4\nabla^+ \tau_1 \right|_g = \mathcal{O}(\alpha')^3 \quad \text{as } \alpha' \to 0. \quad (5.19)$$

*Proof.* In (4.8) we saw that

$$\mathrm{Ric}^+_{ij} - \frac{1}{12}(dH)_{ab\mu i}\psi_{ab\mu j} + 4\nabla^+_i(\tau_1)_j = 0, \quad (5.20)$$

so it suffices to study the second term in (5.20) to obtain (5.18).

Note that

$$(dH)_i \lrcorner \psi_j = \frac{1}{6}(dH)_{ab\mu i}\psi_{ab\mu j}.$$

This observation, together with the heterotic Bianchi identity (5.12), then implies that

$$\frac{1}{12}(dH)_{ab\mu i}\psi_{ab\mu j} = \frac{1}{2}(dH)_i \lrcorner \psi_j = \frac{\alpha'}{2}(F_\theta)^\beta{}_i{}^\mu (F_\theta)_\beta{}^{\nu\rho}\psi_{j\mu\nu\rho}.$$

Using the decomposition

$$F_\theta \lrcorner \psi = 2\pi_7 F_\theta - \pi_{14} F_\theta = 3\pi_7 F_\theta - F_\theta. \quad (5.21)$$

we deduce that

$$\frac{1}{12}(dH)_{ab\mu i}\psi_{ab\mu j} = \alpha'\left(-(F_\theta)^\beta{}_{\mu i}(F_\theta)_{\beta\mu j} + 3(F_\theta)^\beta{}_i{}^\mu \pi_7 F_{\theta\,\beta j\mu}\right). \quad (5.22)$$

Inserting (5.22) in (5.20) gives (5.18). The final result then follows from the condition (5.13) in Definition 5.7 of $\alpha'$-approximate $G_2$-instantons, together with the assumption that $F_\theta$ is bounded as $\alpha' \to 0$. □

Combining Propositions 5.11 and 5.12, we immediately obtain the following.

**Theorem 5.13.** *Let $\theta$ be $\alpha'$-approximate $G_2$-instantons on a principal $K$-bundle over $(M^7, \varphi)$, endowed with an integrable $G_2$-structure with torsion $H$ satisfying the heterotic Bianchi identity (5.12), as in Definition 5.7. Suppose further that the curvature $F_\theta$ is bounded as $\alpha' \to 0$.*

*Let $Q = TM \oplus \mathrm{ad}P \oplus T^*M$ have the Courant algebroid structure defined by the pair $(H, \theta)$. If the divergence is given by*

$$\mathrm{div} = \mathrm{div}^{V_+} - 2\langle 4\tau_1, \cdot \rangle, \quad (5.23)$$

*then a torsion-free $V_+$-compatible generalized connection $D \in \mathcal{D}_0(V_+, \mathrm{div})$ is an approximate generalized Ricci-flat connection, in the sense that the generalized Ricci curvature has norm of order $\mathcal{O}(\alpha')^2$ as $\alpha' \to 0$.*

## 5.3 Approximate solutions on contact Calabi–Yau 7-manifolds

In this subsection, we revisit the setting of contact Calabi–Yau 7-manifolds $\pi : M^7 \to V$ endowed with the $G_2$-structures $\varphi_\varepsilon$ as in Definition 5.1. Recall that $\varphi_\varepsilon$ are integrable $G_2$-structures with torsion $H_\varepsilon$. Recall also the connections $\theta_{\varepsilon,m}^{\delta,k}$ on $TM$ in Definition 5.3 and that we can define the bundle $E = \pi^* TV$ and endow it with the pullback $A$ of the Levi-Civita connection from the Calabi–Yau 3-orbifold $V$. Given this leads us to the following definition.

**Definition 5.14.** *Let $K = G_2 \times \mathrm{SU}(3)$ and identify $G_2$ and $\mathrm{SU}(3)$ with their standard matrix representations (acting on $\mathbb{R}^7$ and $\mathbb{C}^3$ respectively). We can define a principal $K$-bundle $P$ over $M$ whose natural associated vector bundle is $TM \oplus E$. We can then define a connection $\theta$ on $P$ using the pair of connections $\theta_{\varepsilon,m}^{\delta,k}$ and $A$.*

*We also endow the Lie algebra $\mathfrak{k}$ of $K$ with the pairing $\langle \cdot, \cdot \rangle : \mathfrak{k} \otimes \mathfrak{k} \to \mathbb{R}$ with respect to the splitting $\mathfrak{k} = \mathfrak{g}_2 \oplus \mathfrak{su}(3)$:*

$$\langle (X_1, Y_1), (X_2, Y_2) \rangle = -\mathrm{tr}(X_1 X_2) + \mathrm{tr}(Y_1 Y_2). \quad (5.24)$$

*Clearly, $\langle \cdot, \cdot \rangle$ is non-degenerate, bilinear and symmetric.*

Theorem 5.6 then states that if we are given *any* sequence $\alpha' \to 0$, then we can choose positive parameters $\varepsilon = \varepsilon(\alpha')$, $k = k(\alpha')$ and real parameters $\delta, m$ *independent of* $\alpha'$ so that $H = H_\epsilon$ and $\theta$ given in Definition 5.14 satisfy the heterotic Bianchi identity (5.12) with $\langle \cdot, \cdot \rangle$ as in (5.24) (which we notice is $\alpha'$-independent). Moreover, since $A$ is a $G_2$-instanton by [LSE23, Lemma 3.1], (5.11) in Theorem 5.6 also gives that the curvature $F_\theta$ of $\theta$ satisfies the first condition in (5.13).

Theorems 5.10 and 5.13 imply that if $\theta$ also satisfies the second condition in (5.13) and has bounded curvature as $\alpha' \to 0$, then the "approximate" solutions to the heterotic $G_2$ system given by Theorem 5.6 give rise to approximate coupled $G_2$-instantons and an approximate generalized Ricci-flat connection on the associated Courant algebroid. This is what we now show.

**Theorem 5.15.** *Let $M^7$ be a contact Calabi–Yau 7-manifold as in Definition 5.1. Suppose we are given any sequence of positive numbers $\alpha' \to 0$. Let $\varepsilon = \varepsilon(\alpha') > 0$, $k = k(\alpha') > 0$, $\delta, m \in \mathbb{R}$ be the associated parameters given by Theorem 5.6 and let $M$ be endowed the integrable $G_2$-structure $\varphi_\varepsilon$ given in (5.1) with torsion $H = H_\varepsilon$. Let $P, \theta, \langle \cdot, \cdot \rangle$ be the principal $K$-bundle, connection and pairing on $\mathfrak{k}$ given in Definition 5.14.*

*Then the heterotic Bianchi identity (5.12) is satisfied and $\theta$ are $\alpha'$-approximate $G_2$-instantons in the sense of Definition 5.7 with bounded curvature as $\alpha' \to 0$. Hence, $(\varphi_\varepsilon, H_\varepsilon, \theta)$ are $\alpha'$-approximate coupled $G_2$-instantons in the sense of (5.15) and the Courant algebroid $Q = TM \oplus \mathrm{ad}P \oplus T^*M$ with structure $(H, \theta)$ and divergence as in (5.23) has a torsion-free $V_+$-compatible generalized connection with generalized Ricci curvature with norm of order $\mathcal{O}(\alpha')^2$ as $\alpha' \to 0$.*

*Proof.* As explained before the statement, we need only show that the curvature $F_\theta$ is bounded and that the second condition in (5.13) holds. By (5.14), we see that this second condition is equivalent to

$$\left|\nabla^{\theta,+}(F_\theta \wedge \psi_\varepsilon)\right|_{g_\varepsilon} = \mathcal{O}(\alpha')^2 \quad \text{as } \alpha' \to 0. \tag{5.25}$$

We already remarked that the connection $A$ on $E$ is a $G_2$-instanton and that it is pulled-back from $V$. Hence, its curvature $F_A$, and the norm of $F_A$, are $\alpha'$-independent since the metric $g_\varepsilon$ on $M$ is $\alpha'$-independent when restricted to basic forms by (5.2). Moreover, $F_A \wedge \psi_\varepsilon = 0$ and so (5.25) is trivially satisfied for $F_\theta = F_A$.

Given this discussion and the definition of $\theta$, it now suffices to show that the curvature $R^{\delta,k}_{\varepsilon,m}$ of $\theta^{\delta,k}_{\varepsilon,m}$ has bounded norm as $\alpha' \to 0$ and satisfies

$$\left|\nabla^{\theta^{\delta,k}_{\varepsilon,m},+}(R^{\delta,k}_{\varepsilon,m} \wedge \psi_\varepsilon)\right|_{g_\varepsilon} = \mathcal{O}(\alpha')^2 \quad \text{as } \alpha' \to 0. \tag{5.26}$$

In [LSE23, Proposition 3.17], the curvature $R^{\delta,k}_{\varepsilon,m}$ was written in terms of the local orthonormal coframe given in Definition 5.2 as:

$$R^{\delta,k}_{\varepsilon,m} = F_A + \frac{1}{2}k\varepsilon^2(1 - \delta + m)\omega\mathbf{I} + \frac{k^2\varepsilon^2}{4}\mathbf{Q}^\delta_m,$$

where $F_A$ is the curvature of the connection $A$ as above, $\mathbf{I}$ is given in (5.7) and $\mathbf{Q}^\delta_m$ depends only on $\delta$ and $m$ (and the local coframe), so is independent of $\alpha'$. Since all the terms except $F_A$ involve at least a factor of $k\varepsilon$, which tends to zero as $\alpha' \to 0$, we deduce that

$$|R^{\delta,k}_{\varepsilon,m} - F_A|_{g_\varepsilon} \to 0 \quad \text{as } \alpha' \to 0.$$

Since we already established that $F_A$ has bounded norm as $\alpha' \to 0$, the same must be true for $R^{\delta,k}_{\varepsilon,m}$.

We already saw the expression for $R^{\delta,k}_{\varepsilon,m} \wedge \psi_\varepsilon$ in (5.10). Note there that the matrix $\mathbf{M}^\delta_m$ is again independent of $\alpha'$. Recall that $\theta^{\delta,k}_{\varepsilon,m}$ is given in (5.9) and note that taking $\delta = 1$, $m = 0$ and $k = 1$ in this expression leads to the Bismut connection $\nabla^+$, and instead taking $\delta = m = 0$ and $k = 1$ yields the Levi-Civita connection. Altogether, we see that taking derivatives using $\nabla^{\theta^{\delta,k}_{\varepsilon,m},+}$ cannot decrease the powers of $k$ and $\varepsilon$ that already appear in (5.10). Therefore, the norm of $\nabla^{\theta^{\delta,k}_{\varepsilon,m},+}(R^{\delta,k}_{\varepsilon,m} \wedge \psi_\varepsilon)$ must have at least the same order as $\alpha' \to 0$ as the norm of $R^{\delta,k}_{\varepsilon,m} \wedge \psi_\varepsilon$. Since we are already given that this latter quantity is of order $\mathcal{O}(\alpha')^2$ as $\alpha' \to 0$ by (5.11), we deduce that (5.26) holds as desired. □

*Remark* 5.16. Theorem 5.15 shows that it is justified to say that the results in [LSE23], summarised in Theorem 5.6, truly lead to "approximate" solutions to the heterotic $G_2$ system. ○

# A  Generalised Ricci tensor

In this section, we briefly recall the construction in [GF14, GF19] of the generalized Ricci tensor using torsion-free generalized connections. We fix a string algebroid $E$ over a smooth manifold $M$, as in Section 2.2. Recall that a generalized connection on $E$ is given by a differential operator

$$D \colon \Omega^0(E) \to \Omega^0(E^* \otimes E)$$

satisfying the *anchored* Leibniz rule and compatibility with the pairing on $E$, that is,

$$D_a(fb) = \pi(a)(f)b + fD_ab, \qquad \pi(a)\langle b, c\rangle = \langle D_a b, c\rangle + \langle b, D_a c\rangle.$$

Associated to a generalized connection $D$ there are two natural quantities, given by the torsion tensor $T_D \in \Omega^0(\Lambda^3 E)$

$$T_D(a, b, c) = \langle D_a b - D_b a - [a, c], c\rangle + \langle D_c a, b\rangle$$

and the divergence operator $\operatorname{div}_D \colon \Omega^0(E) \to C^\infty(M)$

$$\operatorname{div}_D(a) = \operatorname{tr} Da.$$

Given a generalized metric $\mathbf{G}$ on $E$, there are many Levi-Civita generalized connections [CSCW11, GF19], that is, torsion-free and compatible, in the sense that $D\mathbf{G} = \mathbf{0}$ or, equivalently,

$$D(\Omega^0(V_\pm)) \subset \Omega^0(V_\pm),$$

and this phenomenon persists even if we fix the associated divergence operator. To see this, note that the compatibility with $\mathbf{G}$ implies that we can decompose $D$ in terms of four differential operators

$$\begin{aligned} D_+^+ \colon \Omega^0(V_+) \to \Omega^0(V_+^* \otimes V_+), & \qquad D_-^+ \colon \Omega^0(V_+) \to \Omega^0(V_+^* \otimes V_+), \\ D_+^- \colon \Omega^0(V_-) \to \Omega^0(V_+^* \otimes V_-), & \qquad D_-^- \colon \Omega^0(V_-) \to \Omega^0(V_-^* \otimes V_-). \end{aligned} \tag{A.1}$$

We refer to [GF19, Section 3.1] for a more precise statement of the next result.

**Lemma A.1** ([GF19]). *Let $(\mathbf{G}, \operatorname{div})$ be a pair given by a generalized metric $\mathbf{G}$ and a divergence operator $\operatorname{div}$ on a string algebroid $E$. Then, any torsion-free generalized connection compatible with $\mathbf{G}$ has mixed type operators $D_-^+$ and $D_+^-$ uniquely given by* (2.10). *Furthermore, there exists a non-empty affine space of torsion-free generalized connections compatible with $\mathbf{G}$ and with divergence $\operatorname{div}$, modelled on the sections of*

$$\Sigma_0^+ \oplus \Sigma_0^-,$$

*where $S_0^2 V_\pm$ the space of trace-free symmetric two-tensors and $\Sigma_0^\pm = (S_0^2 V_\pm \otimes V_\pm)/S^3 V_\pm$.*

With the previous classification in hand, we can define the generalized Ricci tensors of a pair $(\mathbf{G}, \operatorname{div})$, as follows. Firstly, given a generalized connection $D$ on $E$, we can define its *curvature operator*

$$R_D(a, b)c := D_a D_b c - D_b D_a c - D_{[a,b]} c$$

acting on a triple of sections $a, b, c \in \Omega^0(E)$. This is not a tensorial object, unlike in ordinary geometry [Gua10]. Recall also from Section 2.2 that a pair $(\mathbf{G}, \operatorname{div})$ determines uniquely an isomorphism $E \cong T \oplus \operatorname{ad} P \oplus T^*$ and a tuple $(g, H, \theta, \zeta_+\zeta_-)$ satisfying the Bianchi identity (2.4) (cf. Proposition 2.9).

**Proposition A.2** ([GF14, GF19]). *Let $(\mathbf{G}, \operatorname{div})$ be a pair given by a generalized metric $\mathbf{G}$ and a divergence operator $\operatorname{div}$ on a string algebroid $E$. Then, there are well-defined* generalized Ricci tensors

$$\operatorname{Rc}_{\mathbf{G},\operatorname{div}}^+ \in V_- \otimes V_+ \quad \text{and} \quad \operatorname{Rc}_{\mathbf{G},\operatorname{div}}^- \in V_+ \otimes V_-,$$

*uniquely determined by $(\mathbf{G}, \operatorname{div})$ via the formula*

$$\operatorname{Rc}_{\mathbf{G},\operatorname{div}}^+(a_-, b_+) = \operatorname{tr}_{V_+} \left( d_+ \to R_D(d_+, a_-) b_+ \right), \qquad \operatorname{Rc}_{\mathbf{G},\operatorname{div}}^-(b_+, a_-) = \operatorname{tr}_{V_-} \left( c_- \to R_D(c_-, b_+) a_- \right),$$

*where $R_D(a.b)c := [D_a, D_b]c - D_{[a,b]}c$ for any choice of torsion-free generalized connection $D$ compatible with $\mathbf{G}$ and with divergence $\operatorname{div}$.*

*Furthermore, in terms of the isomorphism $E \cong T \oplus \operatorname{ad} P \oplus T^*$ and the tuple $(g, H, \theta, \zeta_+\zeta_-)$ determined by $(\mathbf{G}, \operatorname{div})$, the generalized Ricci tensor $\operatorname{Rc}_{\mathbf{G},\operatorname{div}}^+$ (resp. $\operatorname{Rc}_{\mathbf{G},\operatorname{div}}^-$) is explicitly given by* (2.16) *(resp.* (2.17)*).*

An explicit formula for the generalized Ricci tensors in Proposition 2.9 was first computed in [GF14] via a specific choice of torsion-free generalized connection. Here we recall the expression for such a generalized connection, following [GFRT16]. We will need this expression for the proof of Proposition 2.18.

**Lemma A.3** ([GFRT16])**.** *Let $(\mathbf{G}, \mathrm{div})$ be a pair as before, on a string algebroid $E$ over a smooth manifold of dimension $n \geq 2$. Then, there exists a $\mathbf{G}$-compatible torsion-free generalized connection $D$ with $\mathrm{div}_D = \mathrm{div}$, with pure-type operators $D_+^+$ and $D_-^-$ constructed as follows: define $\varepsilon \in \Omega^0(E)$ by $\langle \varepsilon, \cdot \rangle := \mathrm{div}^G - \mathrm{div}$. Via the isomorphism $E \cong T \oplus \mathrm{ad}P \oplus T^*$ provided by $\mathbf{G}$, we can uniquely write*

$$\varepsilon = \sigma_+(\zeta_+^\sharp) + z + \sigma_-(\zeta_-^\sharp)$$

*for $\zeta_\pm \in \Omega^0(T^*)$ and $z \in \Omega^0(\mathrm{ad}P)$. Then, $D$ is defined by*

$$\begin{aligned}
D_{b_+} d_+ &= \sigma_+(\nabla_Y^{1/3} W) + \tfrac{1}{n-1}\sigma_+(g(Y,W)\zeta_+^\sharp - \zeta_+(W)Y), \\
D_{a_-} c_- &= \sigma_-(\nabla_X^{-1/3} Z - \tfrac{2}{3}g^{-1}\langle i_X F_\theta, t\rangle - \tfrac{1}{3}g^{-1}\langle i_Z F_\theta, r\rangle) \\
&\quad + d_X^\theta t - \tfrac{2}{3}F_\theta(X,Z) - \tfrac{1}{3}[r,t] \\
&\quad + \tfrac{1}{\dim \mathfrak{k}+n-1}((\langle r,t\rangle - g(X,Z))(z + \sigma_-(\zeta_-^\sharp)) - (\langle z,t\rangle - \zeta_-(Z))a_-),
\end{aligned}\tag{A.2}$$

*where $\nabla_X^{\pm 1/3} Y = \nabla_X^g Y \pm \tfrac{1}{6} g^{-1} H(X,Y,\cdot)$ and*

$$\begin{aligned}
a_- &= \sigma_-(X) + r = X + r - gX, \\
b_+ &= \sigma_+(Y) = Y + gY, \\
c_- &= \sigma_-(Z) + t, \\
d_+ &= \sigma_+(W).
\end{aligned}\tag{A.3}$$

# B  Special linear algebra

## B.1  Flat $G_2$-structure and decomposition of forms

On Euclidean space $\mathbb{R}^7$ endowed with the canonical basis $e_1, \ldots, e_7$ and the orientation $\mathrm{vol}_0 = e^1 \wedge \ldots \wedge e^7$, the standard flat $G_2$-structure is the three-form $\varphi_0 \in \Lambda^3(\mathbb{R}^7)^*$ explicitly given by (3.1). It induces the (flat) Euclidean metric defined by

$$g_0(X,Y) = \frac{1}{6\mathrm{vol}_0} i_X \varphi_0 \wedge i_Y \varphi_0 \wedge \varphi_0, \tag{B.1}$$

and has Hodge dual (3.3). The Lie group $G_2 \leq \mathrm{SO}(7)$ is defined as the stabiliser of $\varphi_0$ in $\mathrm{SO}(7)$, acting on the space of three-forms $\Lambda^3(\mathbb{R}^7)^*$.

Consider the space of degree $k$, skew-symmetric, multi-linear forms on $\mathbb{R}^7$, which we denote by

$$\Lambda^k := \Lambda^k(\mathbb{R}^7)^*.$$

Following [Kar08], we have the following natural decomposition of the spaces of forms $\Lambda^k$ into irreducible $G_2$-representations (as it is standard, the subscript denotes the dimension of the subspace, and irreducible representations are omitted):

$$\begin{aligned}
\Lambda^2 &= \Lambda_7^2 \oplus \Lambda_{14}^2 \cong \Lambda^1 \oplus \mathfrak{g}_2 \\
\Lambda^3 &= \Lambda_1^3 \oplus \Lambda_7^3 \oplus \Lambda_{27}^3 \cong \Lambda^0 \oplus \Lambda^1 \oplus S_0^2(\mathbb{R}^7) \\
\Lambda^4 &= \Lambda_1^3 \oplus \Lambda_7^4 \oplus \Lambda_{27}^4 \\
\Lambda^5 &= \Lambda_7^5 \oplus \Lambda_{14}^5
\end{aligned}$$

where $S_0^2(\mathbb{R}^7)$ denotes the space of trace-free symmetric two-tensors. In particular, taking the wedge product with $\varphi_0$ yields an isomorphism

$$\Lambda_7^2 \to \Lambda^6 \colon \beta \to \beta \wedge \psi_0$$

and hence, for $\beta \in \Lambda^2$, one has that

$$\beta \in \Lambda_{14}^2 \cong \mathfrak{g}_2 \iff \beta \wedge \psi_0 = 0.$$

## B.2 Spin linear algebra

Let $\operatorname{Cl}(\mathbb{R}^7)$ denote the Clifford algebra of $(\mathbb{R}^7, g_0)$, given by the tensorial algebra of $\mathbb{R}^7$ modulo the relation

$$X \cdot X = -g_0(X, X),$$

for $X \in \mathbb{R}^7$. The algebra $\operatorname{Cl}(\mathbb{R}^7)$ is isomorphic to $\operatorname{End}(\mathbb{R}^8) \oplus \operatorname{End}(\mathbb{R}^8)$, and it admits a real representation $\Delta_7 \simeq \mathbb{R}^8$ with generators (see [FKMS97, p. 261]):

$$\begin{array}{rclrcl}
e_1 &=& E_{18} + E_{27} - E_{36} - E_{45}, & e_2 &=& -E_{17} + E_{28} + E_{35} - E_{46}, \\
e_3 &=& -E_{16} + E_{25} - E_{38} + E_{47}, & e_4 &=& -E_{15} - E_{26} - E_{37} - E_{48}, \\
e_5 &=& -E_{13} - E_{24} + E_{57} + E_{68}, & e_6 &=& E_{14} - E_{23} - E_{58} + E_{67}, \\
e_7 &=& E_{12} - E_{34} - E_{56} + E_{78}, & & &
\end{array} \tag{B.2}$$

where $E_{ij}$ is the standard basis of the Lie algebra $\mathfrak{so}(8)$: it is $-1$ in the position $i,j$ and skew-symmetric. Upon restriction of this representation to $\operatorname{Spin}(7) \subset \operatorname{Cl}(\mathbb{R}^7)$ we obtain the (irreducible) real spin representation

$$\kappa : \operatorname{Spin}(7) \to \operatorname{SO}(\Delta_7).$$

The group $\operatorname{Spin}(7)$ acts transitively on the sphere and $\operatorname{G}_2$ can be identified with the subgroup of $\operatorname{Spin}(7)$ preserving a spinor (see e.g. [FKMS97]).

**Proposition B.1.** *The Lie group $\operatorname{G}_2$ is canonically isomorphic to the subgroup of $\operatorname{Spin}(7)$ preserving the spinor $\eta_0 := (1, 0, \cdots, 0) \in \Delta_7$:*

$$\operatorname{G}_2 \cong \{g \in \operatorname{Spin}(7) : g \cdot \eta_0 = \eta_0\}. \tag{B.3}$$

As a $\operatorname{G}_2$-representation the space of spinors $\Delta_7$ splits into irreducible components as

$$\Delta_7 \simeq \mathbb{R}^8 \simeq \mathbb{R} \oplus \mathbb{R}^7 \simeq \langle \eta_0 \rangle \oplus \Lambda^1,$$

corresponding to the real and purely imaginary octonions. The identification of $\Lambda^1$ inside $\Delta_7$ is simply the embedding $\alpha \in \Lambda^1 \mapsto \alpha^\sharp \cdot \eta_0 \in \Delta_7$, cf. [FKMS97, p. 262]. The relation between the descriptions of $\operatorname{G}_2$ as stabiliser of a 3-form $\varphi_0$ and as a spinor $\eta_0$ are related by [ACFH15, p. 545]:

$$\varphi_0(X, Y, Z) := \langle X \cdot Y \cdot Z \cdot \eta_0, \eta_0 \rangle. \tag{B.4}$$

Following [FKMS97], there is a natural $\operatorname{G}_2$-equivariant map

$$\mu : \Lambda^2 \to \Delta_7 : \beta \mapsto \beta \cdot \eta_0.$$

Using the isomorphisms $\Lambda^2 \simeq \Lambda^1 \oplus \mathfrak{g}_2$ and $\Delta_7 \simeq \langle \eta_0 \rangle \oplus \Lambda^1$, one can easily see that $\mu|_{\mathfrak{g}_2} \equiv 0$ by invariance. Furthermore, as demonstrated in [FKMS97, p. 262], $\mu|_{\Lambda^1}$ is an embedding. This leads to a characterization of the Lie algebra $\mathfrak{g}_2$:

$$\mathfrak{g}_2 = \{\beta \in \Lambda^2 : \beta \cdot \eta_0 = 0\} \subset \Lambda^2. \tag{B.5}$$

Consider the action of 3-forms on the spinor $\eta_0$

$$\nu : \Lambda^3 \to \Delta_7 : \gamma \mapsto \gamma \cdot \eta_0.$$

Using $\Lambda^3 \simeq \langle \varphi_0 \rangle \oplus \Lambda^1 \oplus \Lambda^3_{27}$, the different pieces in this decomposition act on $\eta_0$ via the following formulae (see [FI03])

$$\varphi_0 \cdot \eta_0 = -7\eta_0, \qquad (X \lrcorner \psi_0) \cdot \eta_0 = 4X \cdot \eta_0, \qquad \gamma \cdot \eta_0 = 0, \tag{B.6}$$

for all $X \in \mathbb{R}^7$ and $\gamma \in \Lambda^3_{27}$. Consider now the induced 'Dirac-type' map $\psi : \Lambda^3 \to \Delta_7$, defined by

$$\psi : \Lambda^3 \to \Delta_7 : \gamma \mapsto \slashed{\gamma} \cdot \eta_0 = \sum_j e^j \cdot (e_j \lrcorner \gamma) \cdot \eta_0$$

Using $\Delta_7 \simeq \langle \eta_0 \rangle \oplus \Lambda^1$, we can describe how the irreducible components of $\Lambda^3$ behave under $\psi$:

**Lemma B.2.** *Under the decomposition $\Lambda^3 \simeq \langle \varphi_0 \rangle \oplus \Lambda^1 \oplus \Lambda^3_{27}$, the irreducible components of $\Lambda^3$ under the map $\psi$ acts as*

$$\slashed{\varphi}_0 \cdot \eta_0 = -\frac{21}{2}\eta_0, \qquad \slashed{i_X \psi_0} \cdot \eta_0 = 6X \cdot \eta_0, \qquad \slashed{\gamma} \cdot \eta_0 = 0, \tag{B.7}$$

*for all $X \in \mathbb{R}^7$ and $\gamma \in \Lambda^3_{27}$.*

*Proof.* First, let $\gamma = \frac{1}{3!}\gamma_{ijk}e^{ijk}$ an arbitrary 3-form, then $e_\mu \lrcorner \gamma = \frac{1}{2!}\gamma_{\mu jk}e^{jk}$ and consequently

$$\slashed{\gamma} \cdot \eta_0 = e^\mu \cdot (e_\mu \lrcorner \gamma) \cdot \eta_0 = \frac{1}{2!}\gamma_{\mu jk}e^\mu \cdot (e^j \wedge e^k) \cdot \eta_0.$$

Using the canonical embedding $\mathfrak{spin}(7) \simeq \mathfrak{so}(7) = \Lambda^2 \to \mathrm{Cl}(\mathbb{R}^7)$ given by $e^j \wedge e^k \in \Lambda^2 \mapsto \frac{1}{2}e^j \cdot e^k \in \mathrm{Cl}(\mathbb{R}^7)$, in accordance with the convention $v \cdot v = -|v|^2$, cf. [LM90, Prop. 6.2], we then have

$$\slashed{\gamma} \cdot \eta_0 = \frac{1}{4}\gamma_{\mu jk}e^\mu \cdot e^j \cdot e^k \cdot \eta_0.$$

Using the explicit representation (B.2) and $\eta_0 = (1, 0, \cdots, 0) \in \Delta_7 \simeq \mathbb{R}^8$, the last equation $\slashed{\gamma} \cdot \eta_0 = 0$ for $\gamma \in \Lambda^3_{27}$ holds by invariance. The first equation holds, by a direct computation. For the second one, it is enough to prove for $X = e_1 \in \mathbb{R}^7$, hence by invariance for the whole irreducible representation. We have in this case

$$i_X \psi_0 = e^{256} + e^{234} + e^{457} + e^{367},$$

and

$$e^1 \cdot \eta_0 = (E_{18} + E_{27} - E_{36} - E_{45}) \cdot \eta_0 = E_{18} \cdot \eta_0 = (0, \cdots, 0, 1).$$

On the other hand,

$$\slashed{e_1 \lrcorner \psi_0} \cdot \eta_0 = (0, \cdots, 0, 6)$$

and the result follows. $\square$